\documentclass[reqno]{amsart}

\usepackage[margin=1.5in]{geometry}
\usepackage[usenames, dvipsnames]{color}
\usepackage{amsthm, amsmath, amssymb, mathrsfs, cleveref, enumerate, array}

\newtheorem{thm}{Theorem}[section]
\newtheorem{prop}[thm]{Proposition}
\newtheorem{lem}[thm]{Lemma}
\newtheorem{cor}[thm]{Corollary}

\newcommand{\bthm}{\begin{thm}}
\newcommand{\ethm}{\end{thm}}
\newcommand{\blem}{\begin{lem}}
\newcommand{\elem}{\end{lem}}
\newcommand{\bprop}{\begin{prop}}
\newcommand{\eprop}{\end{prop}}
\newcommand{\bcor}{\begin{cor}}
\newcommand{\ecor}{\end{cor}}
\newcommand{\bdf}{\begin{df}}
\newcommand{\edf}{\end{df}}
\newcommand{\bpf}{\begin{proof}}
\newcommand{\epf}{\end{proof}}

\newcommand{\bfl}{\begin{flalign*}}
\newcommand{\efl}{\end{flalign*}}

\newcommand{\ben}{\begin{enumerate}}
\newcommand{\een}{\end{enumerate}}

\newcommand{\N}{\mathbb{N}}
\newcommand{\Z}{\mathbb{Z}}

\newcommand{\R}{\mathbb{R}}

\newcommand{\E}{\mathcal{X}}

\newcommand{\sgn}{\text{sgn}}
\newcommand{\var}{\mathrm{Var}}

\newcommand{\lt}{\langle}
\newcommand{\rt}{\rangle}

\begin{document}

\title[Hydrodynamics for long-range asymmetric systems]{Hydrodynamic limits for long-range asymmetric interacting particle systems}
\author{Sunder Sethuraman and Doron Shahar }
\date{}

\address{\noindent Department of Mathematics, University of Arizona,
  Tucson, AZ  85721
\newline
e-mail:  \rm \texttt{sethuram@math.arizona.edu}
}

\address{\noindent Department of Mathematics, University of Arizona,
  Tucson, AZ  85721
\newline
e-mail:  \rm \texttt{dshahar@math.arizona.edu}
}

\begin{abstract}
We consider the hydrodynamic scaling behavior of the mass density with respect to a general class of mass conservative interacting particle systems on $\Z^n$, where the jump rates are asymmetric and long-range of order $\|x\|^{-(n+\alpha)}$ for a particle displacement of order $\|x\|$.  Two types of evolution equations are identified depending on the strength of the long-range asymmetry.  When $0<\alpha<1$, we find a new integro-partial differential hydrodynamic equation, in an anomalous space-time scale.  On the other hand, when $\alpha\geq 1$, we derive a Burgers hydrodynamic equation, as in the finite-range setting, in Euler scale.

\end{abstract}



\maketitle

\section{Introduction}

In this paper, we consider hydrodynamic limits in a class of mass conserving 
particle systems in several dimensions $n\geq 1$ on $\Z^n$ with certain asymmetric long-range interactions.  These limits, when they exist, capture the space-time scaling limit of the microscopic empirical mass density field of the particles as the solution of a `hydrodynamic' equation governing a macroscopic flow.  
When the interactions are symmetric and finite-range, such limits have been shown in a variety of stochastic particle systems (cf. \cite{Demasi}, \cite{Kipnis}, \cite{Spohn}).   Also, when the interactions are asymmetric and finite-range, for systems, such as `simple exclusion' and `zero-range', as well as other processes, hydrodynamics has been proved (cf. \cite{baha},  \cite{baha1}, \cite{Saada}, \cite{eymar}, \cite{fritz0}, \cite{fritz}, Chapter 8 in \cite{Kipnis}, \cite{Rezakhanlou}, and reference therein).

However, less is known about hydrodynamics when the dynamics is of long-range type, although such processes are natural in applications, for instance with respect to wireless communications.  
The only works, to our knowledge, which considers `long-range' limits are \cite{cedric} and \cite{Jara}, where hydrodynamics of types of symmetric, long-range exclusion and zero-range processes was shown.

In this context, 
our main purpose is to derive the hydrodynamic equation in a general class of asymmetric long-range `misanthrope' models, which includes simple exclusion and zero-range systems. 
Another motivation was to understand if there is a `mode coupling' basis for certain `stationary' fluctuation results in asymmetric long-range models seen in \cite{seth}.  There, the fluctuations of the empirical mass density field, translated by characteristic speeds, was shown to obey in a sense either a stochastic heat or Burgers equation, depending on the strength of the long-range interactions.  One may ask whether such fluctuations could be inferred from associated hydrodynamics through mode coupling analysis (cf. \cite{Spohn}), as is the case with respect to asymmetric finite-range systems.

Informally, the particle systems studied follow a collection of dependent random walks which interact in various ways.  For instance, in the exclusion and zero-range particle systems, the random walks interact infinitesimally in time and space respectively.  In the exclusion process, particles move freely except in that jumps, according to a jump probability $p(\cdot)$, to already occupied locations are suppressed.  Whereas, in the zero-range process, the jump rate of a particle at a site depends on the number of particles at that site, but the location of the jump is freely chosen according to $p(\cdot)$.  In this article, we will consider the general `misanthrope' process, for which features of exclusion and zero-range interactions are combined, so that both the jump rate and location of jump may depend infinitesimally on the other particles.

In such dynamics, as mass is
preserved, that is no birth or death allowed, there is a 
family of product
invariant measures $\nu_\rho$ indexed by density $\rho$.   Let
 $\eta_t(x)$ denote the number of particles at location $x$ at time $t$.  
 
 By `long-range', to be concrete, we mean, for $\alpha>0$ and $d\in \Z^n$, that $p(\cdot)$ takes the form
 $$p(d) = \frac{{\bf 1}(d>0)}{\|d\|^{n+\alpha}},$$
where $d>0$ means $d_i\geq 0$ for $1\leq i\leq n$ and $d\neq 0$.
The form we have chosen may be generalized
as discussed in Subsection \ref{remarks}.

We will start the process in certain `local equilibrium' nonstationary states $\mu^{N}$, that is
when initially
particles are put independently on lattice sites, according to a varying mass density $\rho_0$, where the marginal at vertex $x$ has mean $\rho_0(x/N)$, and $N$ is a scaling parameter.  We will restrict attention to initial densities $\rho_0$ such that the relative entropy of $\mu^{(N)}$ with respect to an invariant measure $\nu_{\rho^*}$ for $\rho^*>0$ is of order $N^n$.  In effect, this means $\rho_0=\rho_0(u)$ is a function which equals a constant $\rho^*$ for all $u$ large.  This restriction is further discussed in Subsection \ref{remarks}.

Consider, the `hydrodynamic' density, where space is scaled by $N$ and time is speeded up by $N^\theta$,
$$\rho(t,u) \ = \ \lim_{\epsilon\downarrow 0}\lim_{N\uparrow\infty} \frac{1}{(2N\epsilon)^n} \sum_{|y/N - u|\leq \epsilon} \eta_{N^\theta t}(y).$$
Initially, by the law of large numbers, $\rho(0,\cdot) = \rho_0(\cdot)$.  Our goal will be to derive, choosing $\theta = \theta(\alpha)$ appropriately, a `hydrodynamic' partial differential equation for $\rho(t,\cdot)$.  

The choice of $\theta$ is usually determined by the time scale needed in order for a single particle to travel a microscopic distance of order $N^n$, or a nonzero macroscopic distance.  When $\alpha>1$, as $p(\cdot)$ has a mean, the travel time is of the same order as in the finite-range asymmetric case, namely of order $N$, indicating $\theta =1$, the `Euler' scale.  While, when $0<\alpha<1$, because of the heavier tail in $p(\cdot)$, the travel time is of shorter duration, and it turns out $\theta$ should be taken as $\theta = \alpha$, an anomalous scale, interestingly the same as in \cite{Jara} when the jumps are symmetric.

Our main results are as follows.  When $0<\alpha<1$ (Theorem \ref{alpha<1}), we derive that the hydrodynamic equation is the weak form of
$$\partial_t \rho(t,u) = \int_{[0,\infty)^n} \frac{F(\rho(t, u-v), \rho(t,u)) - F(\rho(t,u), \rho(t, u+v))}{
\|v\|^{n+\alpha}} dv,$$
where the function $F$ reflects a homogenization of the microscopic rates in the system.

This equation appears novel in the PDE literature.  It is, in a sense, a `long-range'  integro-differential form of a Burgers equation.  Although the step function $\rho_0(u) = 1_{[a,\infty)}$, in $n=1$, is an invariant solution, as no particle moves, one suspects if motion is allowed in $\R^n$, because of the long-range character of the jump probability, the solution may be more regular than in the Burgers equation, where shocks may form in finite time from smooth initial densities.   Although we show existence of weak solutions of the hydrodynamic equation, uniqueness of these solutions, in any particular class of solutions, is not yet known.  As a consequence, the result we show is that all limit points of the mass density field satisfy weakly this hydrodynamic equation.  

However, when $\alpha>1$ (Theorem \ref{alpha>1}), under an additional assumption that the misanthrope system is `attractive', that is a monotonicity condition on the rates (cf. definition in Section \ref{models}), we show that the hydrodynamic equation is a Burgers equation
$$\partial_t \rho(t,x) + \gamma_\alpha\partial_{{\bf 1}(n)} F(\rho(t,x)) = 0,$$
where  $\gamma_\alpha$ is a specified constant, $F$ again depends on particle interactions,
and $\partial_{{\bf 1}(n)}$ is the directional derivative in $x$ in direction $\langle 1, 1, \ldots, 1\rangle$.  For the boundary value $\alpha=1$, one recovers the same hydrodynamic equation as when $\alpha>1$, however with an extra `log' scaling factor in the scaling, as remarked in Subsection \ref{remarks}.

The $\alpha\geq 1$ hydrodynamic equation may be understood in terms of results say in \cite{Saada}, \cite{Rezakhanlou}, for finite-range asymmetric systems.  When $\alpha>1$, the mean of the jump probability $p(\cdot)$ is bounded.  In particular, long jumps are not so likely, and it is perhaps expected in this case that a Burgers equation would be derived.

  In both settings, these hydrodynamic limit results are the first for long-range asymmetric misanthrope processes on $\Z^n$.  In Subsection \ref{remarks}, further remarks on these limits, their assumptions, and extensions are discussed.

The general scheme of proof, as is customary, is to obtain the hydrodynamic equation by an application of an It\^o formula with respect to the evolving empirical mass density.  In this computation, the generator action gives an average of nonlinear rates of interaction, in terms of the occupation variables $\eta_t$.  The main point is to replace this average, which immediately cannot be written in terms of the empirical density itself, by a homogenized or averaged function of the empirical mass density, thereby allowing one to close the equation.

When $0<\alpha<1$, we follow the `entropy' method strategy of Guo-Papanicolaou-Varadhan (cf. \cite{Kipnis}) as invoked in \cite{Jara} for the symmetric long-range zero-range model.  There are however important differences, especially with respect to the `$1$ and $2$-block' estimates, where the general long-range asymmetric misanthrope structure complicates the analysis.  Notably, here the `attractiveness' condition is not used.

When $\alpha\geq 1$, we follow the scheme in \cite{Rezakhanlou} and Chapter 8 \cite{Kipnis} for finite-range systems, although several steps in the infinite-volume long-range setting take on a different character.  The technique is to show a `$1$-block' estimate, and then to close equations, through use of Young measures, by invoking a uniqueness result for measure-valued solutions in \cite{DiPerna}.  Though the `attractiveness' condition on the process is used in two important places, many estimates we make do not rely on this condition. 
 Moreover, verification of several of the conditions in \cite{DiPerna} seems novel, and may be of independent interest.

Finally, returning to part of our motivation discussed above, in light of the form of the hydrodynamics shown for $\alpha>1$, there is no difference in the type of hydrodynamic equation when $\alpha\geq 3/2$ or when $1<\alpha<3/2$, and so it would seem the fluctuation results seen in \cite{seth}, when $\alpha\geq 3/2$, may not have a `mode coupling' interpretation.

The structure of the article is as follows.  In Section \ref{models}, we introduce the processes studied and, in Section \ref{results}, we state our main results, Theorems \ref{alpha<1} and \ref{alpha>1}, and related remarks.  After some preliminaries in Section \ref{preliminaries}, we prove Theorem \ref{alpha<1} in Section \ref{proofsection_alpha<1}, relying on $1$ and $2$-block estimates shown in Section \ref{blocks}.  In Section \ref{proof_section_alpha>1}, we prove Theorem \ref{alpha>1}, stating key inputs, Theorems \ref{msw_thm}, \ref{ent_cond_thm}, \ref{mass_thm}, and \ref{init_cond_thm}, which are then proved in 
Sections \ref{msw_sect}, \ref{ent_cond_section}, \ref{L^1mass_sect}, \ref{init_sect}, with the aid of estimates in Section \ref{couplingsection} and the Appendix.

\tableofcontents

\section{Models}
\label{models}

Let $\N_0 = \N \cup \{0\}$. We will consider a class of $n\geq 1$ dimensional `misanthrope' particle systems evolving on the state space $\E=\N_0^{\Z^n}$, which includes simple exclusion and zero-range systems.  The configuration $\eta_t = \{\eta_t(x): x\in \Z^n\}$ gives the number of particles $\eta_t(x)$ at locations $x \in \Z^n$ at time $t$.  Let $p: \Z^n \to [0, \infty)$ be a single particle transition rate such that $\sum_{d} p(d) < \infty$.

In the simple exclusion process, at most one particle may occupy each site, $\eta(x)=0$ or $1$ for all $x \in \Z^n$.  Informally, each particle carries an exponential rate $1$ clock.  When a clock rings, the particle may displace by $d$ with probability proportional to $p(d)$.  If the destination site is empty, this jump is made, however, if the proposed destination is occupied, the jump is suppressed and the clock resets, hence the name `simple exclusion'.
Formally, the system is the Markov process on $\{0,1\}^\Z \subset \E$ with generator
$$Lf(\eta) = \sum_{x,d} p(d)\eta(x)(1-\eta(x+d))(f(\eta^{x,x+d})-f(\eta))$$
where $\eta^{x,y}$ is the configuration obtained from $\eta$ by moving a particle from $x$ to $y$:
$$\eta^{x,y}(z)= \left\{\begin{array}{cl} \eta(x)-1& \text{ if } z =x\\ \eta(y)+1 & \text{ if } z=y \\ \eta(z) & \text{ if } z\neq x,y. \\  \end{array}\right.$$
See \cite{Liggett} for the construction and further details of the simple exclusion process. 

In the zero-range process, however, any number of particles may occupy a site. Informally, each site $x$ holds an exponential clock with rate $g(\eta(x))$, where $g:\N_0\rightarrow \R_+$ is a fixed function, such that $g(0)=0$ and $g(k)>0$ for $k\geq 1$. When a clock rings, from that site a particle at random displaces by $d$ with chance proportional to $p(d)$.  The name `zero-range' comes from the observation that, infinitesimally, particles interact only with those on the same site.  Formally, the zero range process is a Markov process on $\E$ with generator 
$$Lf(\eta) = \sum_{x,d} p(d)g(\eta(x))(f(\eta^{x,x+d})-f(\eta)).$$
To construct the zero range process, we require that $g$ be Lipschitz.   
See \cite{Andjel} for the construction and further details of the zero range process.

On the other hand, the misanthrope process is a more general 
Markov process on $\E$ with generator 
$$Lf(\eta) = \sum_{x,d} p(d)b(\eta(x), \eta(x+d))(f(\eta^{x, x+d})-f(\eta)),$$
Here, to avoid degeneracies,
$b: \N_0 \times \N_0 \to [0, \infty)$ is such that $b(0,m)=0$ for $m\geq 0$.  We now describe two cases:
(i) If $b(l,m)=0$ for some $l,m\geq 1$, then $b(l,m)=0$ for all $l\geq 0$ and $m\geq M_0$, where $M_0\geq1$ is the first such $m$, and in this case, $b(l,m)>0$ for $1\leq l\leq M_0$, $0\leq m<M_0$.  (ii) If $b(l,m)>0$ for $l,m\geq 1$, we denote $M_0=\infty$, and set $b(l,m)>0$ for $l\geq 1$ and $m\geq 0$.

To get the simple  exclusion process, let $b(l,m)= {\bf 1}(l=1, m=0)$ so that $b(\eta(x), \eta(y)) = \eta(x)(1-\eta(y))$. To recover the zero-range process, let $b(l,m)= g(l)$ so then $b(\eta(x), \eta(y)) = g(\eta(x))$. The name `misanthrope' refers to the observation in \cite{Cocozza}, where the process was introduced, that particles tend to avoid crowded sites, if $b(l,m)$ is increasing in $l$ and decreasing in $m$.

In this paper, we concentrate on `decomposable' misanthrope systems, where $b(l,m) = g(l)h(m)$, in terms of functions $g$ and $h$, satisfying the restrictions on $b(\cdot,\cdot)$ above.
The associated generator reduces to the form
$$Lf(\eta) = \sum_{x,d} p(d)g(\eta(x))h(\eta(x+d))(f(\eta^{x, x+d})-f(\eta))$$
To aid in construction of the process and for other estimates, we will impose that (i) $g$ is Lipschitz:  $|g(k+1)-g(k)|\leq \kappa$ for $k\geq 0$, (ii) $h$ be bounded, in which case, $h$ is also Lipschitz, $|h(k+1)-h(k)|\leq \kappa_1:=2\|h\|_\infty$ for $k\geq 0$, and (iii) $|g(a)h(b) - g(u)h(v)|\leq \kappa_2\big[|a-u| + |b-v|\big]$ for $a,b,u,v\geq 0$.   The last condition (iii) is a sufficient ingredient to construct the process, and forces $g$ to be bounded if $h$ is nontrivial (cf. equation (7.1) in \cite{Saada_construction}).  However, it is not a necessary condition, and will not be used in the main body of the paper.

   Since $g(0)=0$, we have $g(l)\leq \kappa l$. We also have $h(0) >0$.  If $h$ has a zero, and $M_0<\infty$ is the first root, then $h(m)=0$ for $m\geq M_0$; in this case, the process, starting with less than $M_0$ particles per site, remains so in the evolution.  
	
	We refer to \cite{Saada_construction} for further discussion and the construction of the process on a complete, separable space $\mathcal{X}_0 = \{\eta: \|\eta\|_{\mathcal{X}_0}<\infty\} \subset \N_0^{\Z^n}$ with metric $\|\eta-\xi\|_{\mathcal{X}_0} = \sum_{x\in \Z^n} \beta(x)|\eta(x)-\xi(x)|$, where $\beta(\cdot)$ is a suitable positive function on $\Z^n$ such that $\sum_{x\in \Z^n}\beta(x)<\infty$.

\subsection{Long range asymmetric transitions}
 In this article, we concentrate on  `long range' totally asymmetric processes, where $p(\cdot)$ is in form
 $$p(d)= \frac{{\bf 1}(d>0)}{\|d\|^{n+\alpha}}.$$ 
Here, $\|\cdot\|$ is the Euclidean norm, and $d=(d_1, ..., d_n)>0$ means $d_i\geq 0$ for all $i$, but $d\neq 0$. We require $\alpha > 0$ so that $\sum_d p(d) < \infty$.  Although more general transition rates can be treated, as discussed in Subsection \ref{remarks}, the form of $p$ chosen allows for simplified notation and encapsulates the complexity of the more general situation.

We will distinguish three cases $\alpha<1$, $\alpha =1$, and $\alpha >1$. The transition rate $p(\cdot)$ has a finite mean exactly when $\alpha>1$, and the corresponding model shares some of the properties of the finite-range situation where $p$ is compactly supported.  However, when $\alpha<1$, the behavior of the associated process does reflect that long jumps are more likely.  The $\alpha=1$ case, although borderline, turns out in some ways to be similar to the $\alpha >1$ case.

In particular, a random walk with transition rate $p(\cdot)$ will take an order $\gamma_N$ steps to travel an order $N$ distance on $\Z^n$ where
$$
\gamma_N = \left\{\begin{array}{rl}
N & \ {\rm when \ }\alpha>1\\
N/\log(N) & \ {\rm when \ }\alpha =1\\
N^\alpha & \ {\rm when \ } \alpha<1.
\end{array}\right.$$
These orders will be relevant when discussing hydrodynamic space-time scaling of the process.

\subsection{Invariant Measures}
As the decomposable misanthrope system is mass conservative, one expects a family of invariant measures $\nu_\rho$ indexed by particle density $\rho\geq 0$.  In fact, there is a family of translation-invariant product measures, for a general class of misanthrope processes, including the long-range asymmetric decomposable models, when 
$$g(i)h(j)-g(j)h(i)= g(i)h(0)-g(j)h(0)$$
for $0\leq i,j \leq M_0$ if $M_0<\infty$ and $i,j\geq 0$ otherwise, which we will also assume (cf. \cite{Saada_construction}).  In the case $h$ is nontrivial, this implies a linear relation between $g$ and $h$.

To specify the marginal $\Theta_\rho$ of the measure $\nu_\rho = \prod_{x\in \Z^n}\Theta_\rho$, consider the measure $\bar\Theta_\lambda$ on $\N_0$ given by
$$\bar\Theta_{\lambda}(k) = \frac{1}{Z(\lambda)} \frac{\lambda^k\prod_{j=0}^{k-1} h(j)}{\prod_{j=1}^k g(j)}$$
where $Z(\lambda)$ is the normalizing constant.
Let $\rho(\lambda) = \sum_{k\geq 0} k\bar\Theta_\lambda(k)$ be the mean of $\bar\Theta_\lambda$.  Both $Z(\lambda)$ and $\rho(\lambda)$ are well defined for $0\leq \lambda <\lambda_c$, where $\lambda_c = \infty$ if $M_0<\infty$ and $\lambda_c=\liminf_{k\uparrow\infty} g(k)/h(k)$ otherwise.  One can see they are strictly increasing on this range, and so invertible.  Let $\rho_c = \lim_{\lambda \to \lambda_c}\rho(\lambda)$.
Now, for density $\rho\in [0,\rho_c)$, define $\Theta_\rho = \bar\Theta_{\lambda(\rho)}$ where $\lambda =\lambda(\rho)$ is such that $\rho(\lambda) = \rho$.

The functions,
$$\Phi(\rho)= E_{\nu_\rho}[g(\eta(x))], \quad {\rm and } \quad \Psi(\rho)= E_{\nu_\rho}[h(\eta(x))],$$
will play important roles in the sequel.   One can observe that $\Phi,\Psi$ are $C^1$ on their domains. 

Moreover, as $g$ and $h$ are Lipschitz, both
$\Phi$ and $\Psi$ will be Lipschitz. (following the method of Corollary 2.3.6 \cite{Kipnis}).  Also, note by boundedness of $h$ that $\|\Psi\| \leq \|h\|_\infty$.  Hence, we have the following inequalities, 
\begin{eqnarray}
&& g(\eta(x))h(\eta(y))\leq \kappa \|h\|_\infty \eta(x), \  \  \Phi(\rho)\Psi(\rho) \leq \kappa \|h\|_\infty\rho , \ \ {\rm and }\nonumber\\
&& \ \ \ |\Phi(b)\Psi(b) - \Phi(a)\Psi(a)| \leq \kappa \|h\|_\infty |b-a| + \kappa_1\kappa b|b-a|.
\label{LB}
\end{eqnarray}

In later calculations, we will need finite exponential moments of $\eta(x)$ and $g(\eta(x))$ with respect to $\nu_\rho$ for $\rho \in [0, \rho^c)$. Since $g(\eta(x)) \leq \kappa \eta(x)$ for some constant $\kappa > 0$, we note $g(\eta(x))$ will have a finite exponential moment if $\eta(x)$ does. 
We say that FEM satisfied if
$$E_{\nu_\rho}\big[e^{\gamma\eta(x)}\big] = \frac{1}{Z(\lambda)} \sum_{k=0}^\infty \frac{(\lambda e^\gamma)^k\prod_{j=0}^{k-1} h(j)}{\prod_{j=1}^k g(j)} = \frac{Z(\lambda e^\gamma)}{Z(\lambda)}<\infty$$ 
for $\gamma\geq 0$ and $\rho\in [0,\rho_c)$.   FEM is a condition on the rates $g$ and $h$, which we will assume holds throughout.  For instance, if $\lim_{k\rightarrow\infty} h(k)/g(k+1)=0$ or $M_0<\infty$, then FEM holds.

We will also assume that $\rho_c=\infty$ when $h(m)>0$ for $m\geq 0$, noting of course, if $M_0<\infty$, then $\rho_c=M_0$.  This rules out processes for which there is a `max' invariant state, but no bound on the occupation number per site, for instance.

To relate with zero-range and simple exclusion, if
we set $h\equiv 1$, the measure $\nu_\rho$ reduces to the well known family of invariant measures for the zero-range process. However, when $h(1)=0$, we recover that $\nu_\rho$ is the Bernoulli product measure with parameter $\rho\in [0,1]$.

Finally, we remark, with respect to $\nu_{\rho^*}$, one may construct an $L^2(\nu_{\rho^*})$ Markov process as in \cite{seth_ext}.  The associated adjoint $L^*$ may be computed as the generator of the process with reversed jump rates $p^*(d) = p(-d)$ for $d\in \Z^n$.

\subsection{Initial, empirical and process measures}
We will examine the scaling behavior of the process as seen when time is speeded up by $\gamma_N$ and space is scaled by parameter $N\geq 1$.  Let $L_N := \gamma_N L$ and $\eta^N_t: = \eta_{\gamma_Nt}$ for $t\geq 0$.  Let $T^N_t$ be the associated semigroup.  Often, we will drop the superscript `$N$' when the context is clear.

We will focus on the cases $\alpha\neq 1$, discussing the case $\alpha =1$ in Subsection \ref{remarks}.  Then, for $\alpha\neq 1$, we have $\gamma_N = N^{\alpha\wedge 1}$.  The space-time scaling is the `Euler' scaling when $\alpha>1$, but is an anomalous scale when $\alpha<1$.

Define, for $t\geq 0$, the empirical measure $$\pi^N_t = \frac{1}{N^n} \sum_{x \in \Z^n} \eta^N_t(x) \delta_{\frac{x}{N}}.$$
  We will 
use the following notation for spatial integration against test functions $G$:
$$\lt \pi^N_t, G(t,u) \rt = \frac{1}{N^n} \sum_{x \in \Z^n} \eta^N_t(x) G\big(t, \frac{x}{N}\big).$$
For $T>0$ fixed, the measure-valued trajectories $\{\pi^N_t: 0 \leq t \leq T\}$ are in the Skorohod space $D([0, T], M_+(\R^n))$, where  
$M_+(\R^n)$ is the set of positive Radon measures on $\R^n$ endowed with the vague topology.

Suppose that we start the process at level $N$ according to an initial measure $\mu^N$.  We denote the distribution at times $t\geq 0$ by $\mu^N_t := \mu^NT^N_t$.  The initial measures 
 that we will use are such that the law of large numbers holds in probability with respect to an initial density $\rho_0$ as $N\rightarrow\infty$:
$$\langle \pi^N_0, G(u)\rangle = \frac{1}{N}\sum_{x\in \Z^n} G(\frac{x}{N}) \eta_0^N \rightarrow \int G(u) \rho_0(u)du.$$
We will assume that 
$\rho_0:\R^n \rightarrow \R$ is continuous, and the range of $\rho_0$ lies in $[0,\rho_c)$. In this respect, we will take that $\mu_N$ is a product measure, whose marginal at $x\in \Z^n$ is $\Theta_{\rho_0(x/N)}$ with mean $\rho_0(x/N)$.
Moreover, we will assume that the relative entropy of $\mu^N$ with respect to an invariant measure $\nu_{\rho^*}$ for $0<\rho^*<\rho_c$ is of order $N^n$.  
Then,  
for $|x|$ large, 
the marginals of $\mu^N$, in this case would be very close to those of $\nu_{\rho^*}$ and, in particular, $\rho_0(x)\sim \rho^*$.  For convenience, throughout, we will assume that $\rho_0$ is such that it equals $\rho^*$ outside a compact set.  We also remark the measures $\{\mu^N\}$, by their definition, are stochastically bounded by $\nu_{\rho^\#}$ where $\rho^\# = \|\rho_0\|_\infty$.

Define also $\{P^N \}_{N\geq 1}$ to be the sequence of probability measures on the Skorohod space $D([0, T], M_+(\R^n))$, governing the processes $\{\pi^N_\cdot\}_{N\geq 1}$ starting from $\{\mu^N\}_{N\geq 1}$.   Expectation with respect to $P^N$ will be denoted as $E^N$.

\subsection{Additional assumption when $\alpha>1$}

We will assume further the following condition, in force only when $\alpha>1$.
Namely, we will assume the misanthrope process is `attractive', that is $b(n,m)$ is increasing in $n$ and decreasing in $m$.  In other words, the `decomposable' process is `attractive' when $g$ is increasing and $h$ is decreasing.

\section{Results}
\label{results}

We will split the main results according to the settings $\alpha<1$ and $\alpha>1$.

\medskip

Suppose $\alpha<1$. and consider the operator $\mathcal{L}$ acting on smooth, compactly supported test functions $G: \R^n \rightarrow \R$ by
 \begin{equation}
 \label{hyd_op_alpha<1}
 \mathcal{L}(G)= \int_{[0, \infty)^n} \frac{\Phi(G(u-v))\Psi(G(u))-\Phi(G(u))\Psi(G(u+v))}{\|v\|^{n+\alpha}} dv.
 \end{equation}
 Let also 
 $$\nabla^{\alpha,v} G(u) = \frac{G(u+v)-G(u)}{\|v\|^{n+\alpha}}.$$

\bthm \label{alpha<1} 
Suppose $\alpha<1$.  Then, the sequence $\{P^N\}_{N\geq 1}$ is tight, and every limit point $P^*$ 
 is supported on absolutely continuous measures $\pi_t=\rho(t, u)du$ whose densities are weak solutions of the hydrodynamic equation $\partial_t \rho = \mathcal{L}(\rho)$  
with initial condition $\rho(0, u) = \rho_0(u)$.

That is, for test functions $G$ with compact support in $[0,T)\times \R^n$, we have 
\begin{eqnarray*}
&& \int_{\R^n} \rho_0(u) G_0(u)du + \int_0^\infty \int_{\R^n} \rho(s,u) \partial_sG_s(u)  ds \\
 &&\ \ \ \ \ +  \int_0^T\int_{\R^n} \int_{[0, \infty)^n} \Phi\left(\rho(s,u)\right)\Psi\left(\rho(s,u+v)\right) \nabla^{\alpha, v}G_s(u) dv du ds =0.
 \end{eqnarray*}
 \ethm

\medskip

We now assume $\alpha>1$, and state the hydrodynamic limit in this setting.

\bthm \label{alpha>1}
Suppose $\alpha >1$, and in addition that the process is `attractive'.
Then, $\{P^N\}_{N\geq 1}$ converges weakly to the point mass supported on the absolutely continuous measure $\pi_t=\rho(t, u)du$ whose density is the weak entropy solution (cf. \eqref{entropy_solution_defn}) of the hydrodynamic equation
\begin{equation}
\label{hyd_alpha>1}\partial_t \rho + \gamma_\alpha \partial_{\textbf{1}(n)} [\Phi(\rho)\Psi(\rho)]=0,
\end{equation}
with initial condition $\rho(0, u) = \rho_0(u)$.  Here, $\textbf{1}(n)$ is the unit vector in the direction $\lt 1, ..., 1\rt$, and $\gamma_\alpha$ is the constant defined by $\gamma_\alpha = \sum_{\|d\|=1}^\infty d_1/\|d\|^{n+\alpha}$.
\ethm

We comment, as is well known, the `Burgers' equation \eqref{hyd_alpha>1} may not have a strong solution for all times.  However, 
a weak solution $\rho(t,u)$ is an `entropy' solution if in the weak sense,
\begin{eqnarray}\label{entropy_solution_defn}
&& \partial_t |\rho-c|+\gamma_\alpha \partial_{\textbf{1}(n)}[\sgn(\rho-c)(\Phi\Psi(\rho)-\Phi\Psi(c))] \leq 0, \ \ {\rm where \ }c\in \R \ {\rm and} \\
&&\exists {\rm \ a \ null \ set \ }\mathcal{E}\subset [0,T] {\rm \ such \ that \ } \lim\limits_{\stackrel{t \searrow 0}{t\not\in \mathcal{E}}} \int_{-R}^R |\rho_t(u)- \rho_0(u)| du =0, \ \ {\rm  for \ all \ }R>0.\nonumber
\end{eqnarray}
In this case, Kru$\breve{\text{z}}$kov proved that there is a unique bounded entropy solution if $\rho_0$ is bounded, which is implied by our assumptions \cite{Kruz}. 

\subsection{Remarks}
\label{remarks}

We now make several remarks about Theorems \ref{alpha<1} and \ref{alpha>1}.
\vskip .1cm

{\bf 1.} {\it Uniqueness of solution}. When $\alpha<1$, an open question is to understand in what sense the weak solution is unique.  If there is an unique weak solution $\rho(t,u)$ of the hydrodynamic equation, then $P^N$ would converge weakly to $\delta_{\rho(t,u)du}$.  However, it is not clear what additional criteria, as in the finite-range or $\alpha>1$ setting, needs to be imposed to ensure a unique weak solution.

In this context, we note, in \cite{Jara}, for certain attractive long-range symmetric zero-range evolutions, with symmetric jump rate $p^{sym}(d) = p(d) + p(-d)$, the hydrodynamic equation derived is $\partial_t \rho = \mathcal{L}^{sym} \rho$ where 
$$\mathcal{L}^{sym}(G)= \int_{\R^n} \frac{\Phi(G(u-v))\Psi(G(u))-\Phi(G(u))\Psi(G(u+v))}{\|v\|^{n+\alpha}} dv$$ 
is isotropic, as opposed to \eqref{hyd_op_alpha<1}.  Uniqueness of weak solution, under an `energy' condition, is shown there.  Symmetric long-range exclusion processes are also considered in \cite{Jara}.   However in such models, as is well known, the hydrodynamic equation is linear, and so uniqueness of solution is more immediate.

\vskip .1cm
{\bf 2.} {\it General jump rates.} The jump rate $p(\cdot)$ may be generalized to a larger class, in which jumps are allowed in all directions.  When $\alpha<1$, the jump rate can be in form, say 
\begin{equation*}
\label{gen_p}
p^{gen}(d) = \beta(d)/\|d\|^{n+\alpha} \ {\rm where \ } \beta(y) = \sum _{i=1}^n \big[b_i^+{\bf 1}(y\cdot e_i\geq 0) + b_i^-{\bf 1}(y\cdot e_i\leq 0)\big]{\bf 1}(y\neq 0),
\end{equation*}
 in terms of constants $\{b_i^+, b_i^-\}_{i=1}^n$, and $\{e_i\}_{i=1}^n$ is the standard basis.  We note, in this case, $p^{gen}$ may even be symmetric as in \cite{Jara}.

However, when $\alpha>1$, the same generalization is allowed, except the jump rate must have a drift, $\sum d p^{gen}(d)\neq 0$.

Under these generalizations, the form of Theorems \ref{alpha<1} and \ref{alpha>1} remain the same except that the hydrodynamic equation now involves straightforwardly the constants $\{b_i^+, b_i^-\}$: When $\alpha<1$, weakly
$$
\partial_t \rho =   \sum_{i=1}^n \sum_{\sigma=+, -}b_i^\sigma \int_{\sigma\big[e_i\cdot v\big]>0} \frac{\Phi(\rho(u-v))\Psi(\rho(u))-\Phi(\rho(u))\Psi(\rho(u+v))}{\|v\|^{n+\alpha}} dv,$$
and when $\alpha>1$, in \eqref{hyd_alpha>1} and \eqref{entropy_solution_defn}, $\gamma_\alpha\partial_{\textbf{1}(n)}$ is replaced by $\sum d p^{gen}(d)\cdot \nabla$.
The proofs are the same, albeit with more notation.
\vskip .1cm

{\bf 3.} {\it Case $\alpha =1$.}  When $\alpha=1$, a log correction is needed in the definition of the empirical measure since the jump rate $p$ does not have mean, but just `barely' so in that $\sum_{\|d\|\leq N} dp(d) = O(\log(N))$.  In this case, instead of $\pi^N_t$, we should use the measures $\frac{1}{N^n} \sum_{x \in \Z^n} \eta_{(N/\log(N))t}(x) \delta_{\frac{x}{N}}$.  The arguments, when $\alpha>1$, are straightforwardly adapted to yield the equation $\partial_t \rho +  \partial_{\textbf{1}(n)} [\Phi(\rho)\Psi(\rho)]=0$, the constant $\gamma_\alpha$ being replaced by $1$.

\vskip .1cm

{\bf 4.} {\it Long-range communication $\alpha>2$ versus $1\leq \alpha\leq 2$}.  In the Euler scale, when $\alpha>2$, as opposed to when $1<\alpha\leq 2$, the influence from long distances to the origin, say, is minimal.  From considering the single particle displacement rates, the chance a particle displaces by order $N$ is of order $N^{1-\alpha}$.  So, the likelihood of particles a distance of order $N$ away from the origin to pass by is minimal when $\alpha>2$, but this chance it appears is nontrivial when $1<\alpha\leq 2$.  

It seems not clear how to use the method in \cite{Rezakhanlou} to approximate $L^1$ initial densities by arbitrary ones. 
 Hence, rather than start in an $L^1$ density $\rho_0$, under which the system would have only a finite number of particles at each scaling level $N$ as in \cite{Rezakhanlou}, we have tried to understand infinite volume effects, using the `entropy' method, by starting in a non-integrable density $\rho_0$.  That $\rho_0(u)=\rho^*$ for large $u$ is a consequence of this method.

\vskip .1cm
{\bf 5.} {\it Use of `attractiveness'  
when $\alpha>1$.} Only for the proof of Theorem \ref{alpha>1} is `attractiveness'
used.  This condition allows to show in Step 1 of Section \ref{proof_section_alpha>1} that solutions $\rho$ are in $L^\infty$ when $M_0=\infty$.
 However, when $M_0<\infty$, we have a priori that $\rho\in L^\infty$ and `attractiveness' is not needed for this point.
 
On the other hand, `attractiveness' is used to rewrite the generator of a coupled process in \eqref{alpha>1_generator}, and then to bound it in Lemma \ref{alpha>1martbound}.  These are important ingredients for the `ordering' Lemma \ref{orderinglemma}, which is used to show a `measure weak' formulation of the entropy condition in Theorem \ref{ent_cond_thm}, proved in Section \ref{ent_cond_section}.  

\vskip .1cm
{\bf 6.} {\it Initial conditions.}  Only in the proof of Theorem \ref{alpha>1} is the full description of the initial measures $\mu_N$ used. In particular, the full structure is employed in Step 3a in Subsection \ref{coupmicro_sect}, for the proof of the entropy condition inequality. However, with respect to the proof of Theorem \ref{alpha<1}, we note only the fact that the marginals of $\mu_N$ at $x\in \Z^n$ have mean $\rho_0(x/N)$ is used.

\section{Preliminaries}
\label{preliminaries}
Throughout this paper, 
a test function will be a smooth $C^{1,2}$ function $G: [0, T) \times \R^n \to \R$ with compact support.  Typically, the letter $R$ will be such that  
the support of $G$ will be contained in $[0, T)\times [-R, R]^n$. 
Define $\|G\| = \sup_{t,u} |G(t,u)|$, and similarly $\|\nabla G\|$, $\|\nabla^2 G\|$ and $\|\partial_s G\|$.

For $y\in \Z^n$, let $\tau_y$ represent the shift operator: 
$\tau_y(\eta(x))= \eta(x+y)$ and $\tau_y(f(\eta))= f(\tau_y\eta)$.
Define, for $d\in \Z^n$, $h_d(\eta(x))= \tau_dh(\eta(x))= h(\eta(x+d))$.

Define also $|y|= \max\{y_1, ..., y_n\}$ for $y= (y_1, ..., y_n) \in \Z^n$.  In later calculations, we will use the notion of an `$l$-block' average of a function $f=f(\eta)$:
That is, define
$$f^l(\eta) = \frac{1}{(2l+1)^n}\sum_{|y| \leq l} \tau_y f(\eta).$$
In particular, $\eta^l(x) = \frac{1}{(2l+1)^n}\sum_{|y| \leq l} \eta(x+y)$.

Form now the mean-zero martingale with respect to $\langle G, \pi^N_t\rangle$:
$$
M^{N,G}_t  :=\lt \pi^N_t, G_t \rt -\lt \pi^N_0, G_0 \rt - \int_0^t \lt \pi^N_s, \partial_sG_s \rt ds -\int_0^t N^{1\wedge \alpha} L_N\lt \pi^N_s , G_s \rt ds.
$$
Also, with respect to its quadratic variation, 
$$\lt M^{N,G}\rt_t := \int_0^t N^{1\wedge \alpha} L_N[(\lt \pi^N_s, G_s \rt)^2]-2N^{1\wedge \alpha} \lt \pi^N_s, G_s \rt L_N \lt \pi^N_s, G_s \rt ds,$$
we have that $(M^{N,G}_t)^2 - \lt M^{N,G}\rt_t$ is a mean-zero martingale.

Explicitly, we may compute 
$$N^{1 \wedge \alpha} L_N\lt \pi^N_s , G_s \rt = \frac{N^{1\wedge \alpha}}{N^n}\sum\limits_{x \in \Z^n} \sum\limits_{\|d\|=1}^\infty \frac{1}{\|d\|^{n+\alpha}}  gh_d(\eta_s(x)) \big[G_s\big(\frac{x+d}{N}\big)-G_s\big(\frac{x}{N}\big)\big]$$
and
\begin{equation}
\label{gen_quad}
\lt M^{N,G}\rt_t = \int_0^t \frac{N^{1 \wedge \alpha}}{N^{2n}}\sum\limits_{x \in \Z^n} \sum\limits_{\|d\|=1}^\infty \frac{1}{\|d\|^{n+\alpha}}gh_d(\eta_s(x))\big[G_s\big(\frac{x+d}{N}\big)-G_s\big(\frac{x}{N}\big)\big]^2 ds.
\end{equation}

Here, and in the body of the paper, our convention will be that the sums over $d$ implicitly contain the restriction that $d=(d_1,\ldots,d_n)>0$, that is $d_i>0$ for $1\leq i\leq n$, as $p$ is supported on such $d$, to reduce notation.

\subsection{Entropy and Dirichlet forms}
\label{entropy_section}

Recall the distribution of the process at the $N$th level at time $t\geq 0$, $\mu^N_t = \mu^N T^N_t$.   Consider the relative entropy $H(\mu^N_t| \nu_{\rho^*})$ of $\mu^N_t$ with respect to the invariant measure $\nu_{\rho^*}$.  One may show that $H(\mu^N_t|\nu_{\rho^*})$ is finite, and hence $\mu^N_t$ is absolutely continuous with respect to $\nu_{\rho^*}$.  In terms of the Radon-Nikodym derivative $f^N_t = d\mu^N_t/d\nu_{\rho^*}$, we have $H(\mu^N_t|\nu_{\rho^*}) = H(f^N_t)$ where $H(f)=\int f \log f d\nu_{\rho^*}$.

Define now the Dirichlet form of a density $f$ by
$D(f)=-\int \sqrt{f}L^{sym}\sqrt{f}d\nu_{\rho^*}$, where we define $L^{sym}= (L+L^*)/2$ as the symmetric part of $L$. 
We will on occasion 
define new Dirichlet forms in terms of pieces of the above Dirichlet form.
For $x,y\in \Z^n$, define the bond Dirichlet form as
$$
  D^{x,y}(f) 
= \frac{1}{2} \int p^{sym}(y-x)g(\eta(x))h(\eta(x+y))(\sqrt{f(\eta^{x,y})}- \sqrt{f(\eta)})^2 d\nu_{\rho^*},
 $$
where
 $p^{sym}(d)= (p(d)+p(-d))/2$. Roughly speaking, $D^{x,y}(f)$ is a measure of how much $f(\eta)$ can vary as one particle is moved from $x$ to $y$ or vice versa. In particular, if $D^{x,y}(f)=0$, then $f(\eta)= f(\eta^{x,y})$ when $p(y-x)g(\eta(x))h(\eta(x+y)) \neq 0$.  In terms of these pieces, the `full' Dirichlet form may be written as $D(f)= (1/2)\sum_{x,y} D^{x,y}(f)$. 

It is a calculation to relate the entropy and Dirichlet form as follows:
\begin{equation}
\label{ent_dir}
H(\mu^N_t|\nu_{\rho^*})  +2\gamma_N\int_0^t D(f^N_s)ds\leq H(\mu^N|\nu_{\rho^*}).
\end{equation}
By convexity 
of the Dirichlet form, we have the bound
$D\big(\bar{f}^N_t\big) \leq H(\mu^N_0| \nu_{\rho^*})/(2\gamma_N t)$,
where $\bar{f}^N_t = \frac{1}{t}\int_0^t f^N_s ds$.
Moreover, we have
$H(\mu^N_t|\nu_{\rho^*})\leq H(\mu^N| \nu_{\rho^*}) \leq CN^n$, by our entropy assumption on the initial distributions $\{\mu^N\}$, and so with $C_0=C/(2t)$,
\begin{equation}
\label{dirichlet_bound}
D\big(\bar{f}^N_t\big) \leq \frac{C_0N^n}{\gamma_N}.
\end{equation}

In the finite volume, \eqref{ent_dir} and \eqref{dirichlet_bound} are well-known (cf. \cite{Kipnis}).  In the infinite volume, to obtain finiteness of the relative entropy, \eqref{ent_dir} and \eqref{dirichlet_bound}, we may approximate $\mu^N_t$ by distributions $\mu^{N,R}_t=\mu^NT^{N,R}_t$ of processes with dynamics localized to boxes of growing width $R$ and sites frozen outside, with semigroup $T^{N,R}_t$ starting from the product measure $\mu^N$, for which the relative entropy and localized Dirichlet form of $d\mu^{N,R}_t/d\nu_{\rho^*}$ satisfy \eqref{ent_dir} and \eqref{dirichlet_bound}.  
In particular, by the construction (cf. \cite{Saada_construction}), for Lipschitz functions $u$ on the complete, separable metric space $\mathcal{X}_0$, we have $T^{N,R}_t u \rightarrow T^N_t u$ as $R\uparrow\infty$; also, as $|T^{N,R}_t u(\eta)|, |T^N_t u(\eta)| \leq c_ue^{c_rt}\|\eta\|_{\mathcal{X}_0} + |u|({\bf 0})\in L^1(\mu_N)$, where $c_u$ is the Lipschitz constant with respect to $u$, ${\bf 0}$ is the empty configuration, and $c_r$ is a constant depending on process parameters, we have the `convergence', $E_{\mu^{N,R}_t}[u]\rightarrow E_{\mu^N_t}[u]$, as $R\uparrow\infty$.  Therefore, $\mu^{N,R}_t$ converges weakly to $\mu^N_t$ by the Portmanteau theorem (cf. Section 3.9 in \cite{Durrett}).  Hence, by lower semi-continuity of the relative entropy (cf. \cite{dupuis}), $H(\mu^N_t|\nu_{\rho^*})\leq CN^n$ is finite.  Now, note that the localized Dirichlet form is greater than the Dirichlet form $D_K$ involving only bonds in a fixed box with width $K$ for all large $R$, and that such fixed forms increase as $K$ grows to the full one.  The form $D_K$ is lower semi-continuous, $\liminf_{R\uparrow\infty} D_K(d\mu^{N,R}_s/d\nu_{\rho^*}) \geq D_K(d\mu^N_s/d\nu_{\rho^*})$, by use of the `convergence', and observations
$-2\sqrt{f(\eta)f(\eta^{x,x+y})} = \sup_\epsilon -\epsilon f(\eta) -\epsilon^{-1}f(\eta^{x,x+y})$, $g(\cdot)h(\cdot)$ is Lipschitz in $\mathcal{X}_0$ by use of the earlier construction assumption (iii), and $E_{\nu_{\rho^*}}[f(\eta^{x,x+y})g(\eta(x))h(\eta(x+y))] = E_{\nu_{\rho^*}}[f(\eta)g(\eta(x+y))h(\eta(x))]$ (applied with $f = d\mu^{N,R}_s/d\nu_{\rho^*}$ and $d\mu^N_s/d\nu_{\rho^*}$).   With these ingredients, it is straightforward to conclude \eqref{ent_dir} and \eqref{dirichlet_bound}.  See also \cite{Landim_Mourragui} and references therein for related approaches.

Recall now the `entropy inequality':  For $\gamma>0$,
$$ E_\mu[f] \leq \frac{1}{\gamma} \Big( E_\nu[e^{\gamma f}] + H(\mu|\nu)\Big).$$
Then, with respect to a function $f$ on the configuration space, we have
\begin{equation}
\label{entropy_ineq}
E^N[f(\eta_s)] \leq \frac{1}{\gamma N^n} \big\{\ln  E_{\nu_{\rho^*}}[e^{\gamma N^n f(\eta)}] + CN^n  \big\}.
\end{equation}

A common application of the entropy inequality is
to bound the numbers of particles in various sets.   
\blem\label{particlebound} For $N\geq 1$, $0\leq s\leq T$, and sets $A_N\subset \Z^n$ with ${\rm Card}(A_N) \leq C_1 N^n$, we have
$$\frac{1}{N^n}E^N\big[\sum_{x\in A_N} \eta_s(x)\big] \leq C_1K,$$
where $K=(\ln E_{\nu_{\rho^*}} \big[e^{\gamma\eta(0)}\big] + C/C_1)/\gamma$.
\elem

\bpf
By the entropy inequality \eqref{entropy_ineq}, and finite exponential moments FEM, the left-side of the display, is bounded by 
$$\frac{1}{\gamma N^n} \big\{\ln E_{\nu_{\rho^*}}\big[\exp\big(\gamma\sum_{x\in A_N} \eta(x)\big)\big]+ CN^n \big\}$$
for $\gamma =1$ say.  Since $\nu_{\rho^*}$ is a translation invariant product measure, we have the further bound $(C_1/\gamma)  \ln E_{\nu_{\rho^*}}\big[e^{\gamma\eta(0)}\big]+ (C/\gamma)$ to finish.
 \epf

For later reference, we state the following `truncation' bounds, which holds under FEM, using also the entropy inequality; see p 90-91 in \cite{Kipnis}.

\blem \label{particlebound2} For $R<\infty$ and $0\leq s\leq T$, we have
\begin{eqnarray*}
&&\limsup_{A \to \infty }  \limsup_{N \to \infty }  E^N \big[ \frac{1}{N^n} \sum\limits_{|x| \leq RN} \eta_s(x) {\bf 1}(\eta_s(x) > A)\big] = 0 \ \ {\rm and \ }\\
&&\limsup_{A \to \infty } \limsup_{l \to \infty } \limsup_{N \to \infty }  E^N \big[ \frac{1}{N^n} \sum\limits_{|x| \leq RN} \eta_s^l(x) {\bf 1}(\eta_s^l(x) > A)\big] = 0.
\end{eqnarray*}
\elem

\subsection{Generator and martingale bounds}

We now collect a few useful bounds.

\blem \label{alpha<1martbound} 
Let $G$ be a test function supported on $[0,T)\times [-R,R]^n$.  We have
$$E^N[|N^{1 \wedge \alpha} L_N\lt \pi^N_s , G_s \rt|] \leq C_G$$ where $C_G = \kappa \|h\|_\infty(n+\frac{\sigma_n}{|\alpha-1|}\|\nabla G\|+ \frac{\sigma_n}{\alpha}2\|G\|)K'$ and $K'=2(R+1)^nK$. \elem

\bpf
First, as $h$ is bounded and $g$ is Lipschitz, by \eqref{LB},
we have 
\begin{equation*}
|N^{1 \wedge \alpha} L_N\lt \pi^N_s , G_s \rt| \leq \kappa \|h\|_\infty N^{1\wedge \alpha} \sum\limits_{\|d\|=1}^\infty \frac{1}{\|d\|^{n+\alpha}} \frac{1}{N^n}\sum\limits_{x \in \Z^n} \eta_s(x) \big|G_s\big(\frac{x+d}{N}\big)-G_s\big(\frac{x}{N}\big)\big|.
\end{equation*}
The sum over $d$ can be divided into a sums over $1\leq \|d\| \leq N$ and $\|d\| >N$. We may bound
$\big|G_s\big(\frac{x+d}{N}\big)-G_s\big(\frac{x}{N}\big)\big| $ by $ \|\nabla G\|\frac{\|d\|}{N}{\bf 1}(|x|\leq (R+1)N)$ when $1\leq \|d\|\leq N$, and by $\|G\|\big({\bf 1}(|x|\leq RN)+{\bf 1}(|x+d|\leq RN)\big)$ when $\|d\|>N$.  Hence, we have the further bound
\begin{eqnarray}
|N^{1 \wedge \alpha} L_N\lt \pi^N_s , G_s \rt|&\leq&
\kappa \|h\|_\infty \Big(\frac{N^{1\wedge \alpha}}{N}\sum\limits_{\|d\|=1}^N \frac{1}{\|d\|^{n-1+\alpha}}\|\nabla G\| \frac{1}{N^n}\sum_{|x| \leq (R+1)N} \eta_s(x) \nonumber \\
&&+ N^{1\wedge \alpha}\sum\limits_{\|d\|=N+1}^\infty \frac{1}{\|d\|^{n+\alpha}} \|G\| \frac{1}{N^n}\sum_{\substack {|x| \leq RN \\ |x+d|\leq RN}} \eta_s(x)\Big). 
\label{alpha<1_eq0}
\end{eqnarray}

Both sums over $x$ add over at most $2((R+1)N)^n$ sized regions. 
Hence, by Lemma \ref{particlebound}, the expected value of both sums are less than $2K(R+1)^nN^n$.

Also, the sums over $d$ can be bounded as follows: 
\begin{eqnarray}
\sum\limits_{\|d\|=1}^{aN} \frac{1}{\|d\|^{n-1+\alpha}} &\leq & n+\sigma_n \int_{1}^{aN} \frac{1}{r^\alpha} dr = n+\sigma_n\big(\frac{1}{\alpha-1}-\frac{1}{\alpha-1}\frac{aN}{(aN)^\alpha}\big) \nonumber\\
\sum\limits_{\|d\|=bN+1}^\infty \frac{1}{\|d\|^{n+\alpha}} &\leq & \sigma_n \int_{bN}^\infty \frac{1}{r^{1+\alpha}} dr =  \frac{\sigma_n}{\alpha (bN)^\alpha},
\label{sumsond}
\end{eqnarray}
where $\sigma_n$ is the surface area of the part of an $n$-sphere, of radius $1$ centered at the origin, contained in the first orthant. We note also, an alternate bound, $\sum_{\|d\|=1}^{aN} \|d\|^{-(n-1+\alpha)} \leq \sigma_n N^{1-\alpha}\int_0^a r^{-\alpha}dr$ can be used when $\alpha<1$.

Then, 
\begin{eqnarray*}
&& E^N |N^{1 \wedge \alpha} L_N\lt \pi^N_s , G_s \rt| \\
&&\ \  \leq 
  \kappa \|h\|_\infty K'\Big(\frac{N^{1\wedge \alpha}}{N} \Big(n+ \sigma_n\Big[\frac{1}{\alpha-1}-\frac{1}{\alpha-1}\frac{N}{N^\alpha}\Big]\Big)\|\nabla G\| + N^{1\wedge \alpha} \frac{\sigma_n}{\alpha N^\alpha} 2\|G\|\Big) \leq C_G,
  \end{eqnarray*}
as desired.
\epf

We state here straightforward corollaries of the proof of Lemma \ref{alpha<1martbound}, adjusting the values of $a$ and $b$ in the sums over $d$ near \eqref{sumsond}.

\begin{lem}
\label{lemma23_cor}
We have, when $\alpha<1$, that
$$\lim_{\epsilon\downarrow 0} \lim_{D\uparrow\infty}\lim_{N\uparrow\infty}
E^N \Big|\int_0^t \frac{N^{\alpha}}{N^n} \sum\limits_{x \in \Z^n} \sum\limits_{\substack{\|d\|< \epsilon N \\ \|d\|> DN }} \frac{1}{\|d\|^{n+\alpha}} \eta_s(x) \big[G_s\big(\frac{x+d}{N}\big)-G_s\big(\frac{x}{N}\big)\big] ds \Big | =0.$$
\end{lem}

\begin{lem}
\label{lemma23_cor1}
We have, when $\alpha>1$, that
$$\lim_{\epsilon\downarrow 0} \lim_{N\uparrow\infty}
E^N \Big|\int_0^t \frac{N}{N^n} \sum\limits_{x \in \Z^n} \sum\limits_{\|d\|> \epsilon N} \frac{1}{\|d\|^{n+\alpha}} \eta_s(x) \big[G_s\big(\frac{x+d}{N}\big)-G_s\big(\frac{x}{N}\big)\big] ds \Big | =0.$$
\end{lem}

The difference of quadratic variations can be bounded as follows:

\blem \label{alpha<1quadbound} For $0\leq t_1\leq t_2\leq T$, we have that
$$E^N\big|\lt M^{N,G}\rt_{t_2} - \lt M^{N,G}\rt_{t_1}\big| \leq K_G\frac{|t_2-t_1|}{N^n}$$
 where $K_G= 2\|G\|C_G$. \elem

\bpf Recall the formula for $\langle M^{N,G}\rangle_t$ in the beginning of the section.  By \eqref{LB},
it is enough to show that  
\begin{equation}
\label{lemma4.6_eq}
\kappa \|h\|_\infty \frac{N^{1 \wedge \alpha}}{N^{2n}}E^N \sum\limits_{x \in \Z^n} \sum\limits_{\|d\|=1}^\infty \frac{1}{\|d\|^{n+\alpha}}\eta_s(x) \big| G_s\big(\frac{x+d}{N}\big)-G_s\big(\frac{x}{N}\big)\big|^2  \leq \frac{K_G}{N^n}.
\end{equation}
We can bound one factor $\big| G_s\big(\frac{x+d}{N}\big)-G_s\big(\frac{x}{N}\big)\big|$ by $2\|G\|$. The left-side of the display is then bounded by $2\|G\|/N^n$ times
$$ \kappa \|h\|_\infty N^{1 \wedge \alpha} E^N \sum\limits_{\|d\|=1}^\infty \frac{1}{\|d\|^{n+\alpha}}\frac{1}{N^{n}}\sum\limits_{x \in \Z^n}\eta_s(x) \big| G_s\big(\frac{x+d}{N}\big)-G_s\big(\frac{x}{N}\big)\big|. $$
However, we have already bounded this expression in the proof of Lemma \ref{alpha<1martbound} by $C_G$.  
\epf

\subsection{Tightness of $\{P^N\}$}

We now show, when $\alpha\neq 1$, that the sequence $\{P^N\}$ is tight and therefore weakly relatively compact. For smooth $G$ with compact support, let $P^N_G$ be the induced distribution of $\{\langle \pi^N_t, G\rangle: t\in [0,T]\}$.
To prove that $\{P^N\}$ is tight, it is enough to show that $\{P^N_G\}$ is tight for all such $G$ (cf. Proposition 1.7, Chapter 4 in \cite{Kipnis}).  We will in fact show tightness estimates with respect to the uniform topology, stronger than the Skorohod topology.

\begin{prop}
\label{PNtightness}
The sequence
$\{P^N_G\}$ is tight with respect to the uniform  topology:  For smooth $G$ with compact support in $\R^n$, the following holds.

\ben
\item For every $\epsilon >0$, there is a compact $K \subseteq \R$ such that $\sup_N P^N\big(\lt \pi^N_0, G\rt \notin K\big) \leq \epsilon$.
\item For every $\epsilon>0$, $$\limsup\limits_{\delta \to 0^+} \ \limsup\limits_{N \to \infty} P^N\Big(\sup_{\stackrel{|t-s|< \delta}{0\leq s,t\leq T}} |\lt \pi^N_t, G\rt -\lt \pi^N_s, G\rt|> \epsilon\Big)=0$$
\een
\end{prop}

\bpf To prove the first condition, it is enough to show that 
$\sup_N E^N[|\lt \pi^N_0 , G\rt|]$ is finite.  But, by Lemma \ref{particlebound},
$$E^N |\lt \pi^N_0 , G\rt| \leq \|G\|\frac{1}{N^n}E_{\mu^N}\sum_{|x|\leq RN} \eta_0(x) <\infty.$$

To prove the second condition,
 for $t>s$, we may write
\begin{eqnarray}
\label{tightness1}
&&\sup_{|t-s|< \delta} |\lt \pi^N_t, G\rt -\lt \pi^N_s, G\rt| \\ 
&&\ \ \ \leq  \sup_{|t-s|< \delta} \int_s^t \big|N^{1\wedge \alpha} L_N \lt \pi^N_r, G \rt \big| dr + \sup_{|t-s|< \delta} \big| M^{N,G}_t- M^{N,G}_s\big|. \nonumber
\end{eqnarray}
The second term on the right-side of \eqref{tightness1} is bounded through the triangle inequality, Doob's inequality, and the quadratic variation estimate Lemma \ref{alpha<1quadbound}:
$$E^N  \sup_{|t-s|< \delta} \big| M^{N,G}_t- M^{N,G}_s\big|^2  \leq 4 E^N \sup_{0\leq t \leq T} \big(M^{N, G}_t\big)^2 \leq 16E^N\langle M^{N,G}\rangle_T = O\big({N^{-n}}\big).$$

For the first term on the right-side of \eqref{tightness1}, as is done in the proof of  Lemma \ref{alpha<1martbound}, we bound the integrand by  \eqref{alpha<1_eq0}.  We now analyze the first term in \eqref{alpha<1_eq0}; the other term is similarly handled.  Write the first term as $I_1+I_2$, in terms of a parameter $A$, where
\begin{eqnarray*}
I_1= \kappa \|h\|_\infty\frac{N^{1\wedge \alpha}}{N}\sum\limits_{\|d\|=1}^N \frac{1}{\|d\|^{n-1+\alpha}}\|\nabla G\| \frac{1}{N^n}\sum_{|x| \leq RN} \eta_r(x){\bf 1}(\eta_r(x)\leq A)\\
I_2 = \kappa \|h\|_\infty\frac{N^{1\wedge \alpha}}{N}\sum\limits_{\|d\|=1}^N \frac{1}{\|d\|^{n-1+\alpha}}\|\nabla G\| \frac{1}{N^n}\sum_{|x| \leq RN} \eta_r(x){\bf 1}(\eta_r(x)> A).
\end{eqnarray*}

We may bound $I_1$, as in the proof of Lemma \ref{alpha<1martbound}, by $I_1\leq [C_G/K'](2R+1)^nA$.  Correspondingly, $\sup_{|t-s|<\delta}\int_s^t I_1 dr \leq \delta [C_G/K'](2R+1)^nA$, which vanishes as $\delta\downarrow\infty$.

For the term $I_2$, we use the following approach.  
For each $\delta$, partition $[0,T]$ into $n=\lceil T/\delta \rceil$ intervals $[t_i, t_{i+1}]$ for $i=0, 1, \ldots, n-1$, of length $T/n$. 
Then, 
$$\sup_{|t-s| < \delta} \int_s^t I_2 dr \leq 3 \max_{i} \int_{t_i}^{t_{i+1}} I_2 dr.$$
It follows that
$$P^N\Big(\sup_{|t-s| < \delta} \int_s^t I_2 dr > \epsilon\Big) \leq \sum_{i} P^N\left(\int_{t_i}^{t_{i+1}} I_2 dr > \epsilon/3\right) \leq \frac{3}{\epsilon} \int_0^T E^N[I_2] dr. $$
Since the sum $(N^{1\wedge \alpha}/N)\sum_{1\leq \|d\| \leq N} \|d\|^{-(n-1+\alpha)}$ is bounded (cf. \eqref{sumsond}), by Lemma \ref{particlebound2},
we have that $\lim_{A\uparrow\infty}\lim_{N\uparrow\infty} E^N[I_2] =0$, to finish the proof.
\epf

\section{Proof outline:  Hydrodynamic limits when $\alpha<1$}
\label{proofsection_alpha<1}

  We outline the proof of Theorem \ref{alpha<1}, refering to `$1$ and $2$-block' estimates later proved in Section \ref{blocks}.
\vskip .2cm

{\it Step 1.} First, by Doob's inequality and the quadratic variation bound Lemma \ref{alpha<1quadbound}, for $\epsilon_0>0$,
$$\limsup\limits_{N \to \infty} P^N\Big(\sup_{0\leq t\leq T}\big| M^{N,G}_t  \big|> \epsilon_0 \Big)  \leq \frac{4}{\epsilon_0^2}\limsup\limits_{N \to \infty} E^N \big(M^{N,G}_T \big)^2=0.$$
As  $G$ has compact support, we may choose $t<T$ large enough so that $G_t$, and hence $\lt \pi^N_t, G_t\rt$ vanishes.  Therefore, for such $t$,
$$\limsup\limits_{N \to \infty} P^N\Big(\Big|\lt \pi^N_0 , G_0\rt + \int_0^t  \lt \pi^N_s, \partial_sG_s \rt ds + \int_0^t N^{\alpha} L_N\lt \pi^N_s , G_s\rt ds \Big|> \epsilon_0 \Big)=0.$$

\vskip .2cm

{\it Step 2.} Next, in order that
$\lt \pi^N_0 , G_0\rt + \int_0^t  \lt \pi^N_s, \partial_sG_s \rt ds + \int_0^t N^{\alpha} L_N\lt \pi^N_s , G_s\rt ds$ looks like the weak formulation of a hydrodynamic equation, we will work to replace $\int_0^t N^{\alpha} L_N\lt \pi^N_s , G_s\rt ds$ by appropriate terms.
Noting the generator expression near \eqref{gen_quad},
$$\int_0^t N^{\alpha} L_N\lt \pi^N_s , G_s\rt ds = \int_0^t \frac{N^{\alpha}}{N^n} \sum\limits_{x \in \Z^n} \sum\limits_{\|d\|=1}^\infty \frac{1}{\|d\|^{n+\alpha}} gh_d(\eta_s(x)) \big[G_s\big(\frac{x+d}{N}\big)-G_s\big(\frac{x}{N}\big)\big] ds.$$ 
   We now truncate the sum over $d$ to when  $\|d\|$ is at least $\epsilon N$ and at most $DN$.  By Lemma \ref{lemma23_cor}, as $gh_d(\eta(x)) \leq \kappa \|h\|_\infty\eta(x)$ (cf. \eqref{LB}), the excess vanishes,
  where $\epsilon\downarrow 0$ and $D\uparrow \infty$, after $N\uparrow\infty$.
Therefore, after limits on $N$, $\epsilon$ and $D$ are taken, 
\begin{eqnarray*}
&&\lt \pi^N_0 , G_0\rt + \int_0^t  \lt \pi^N_s, \partial_sG_s \rt ds \\
&&\ \ \ + \int_0^t \frac{N^{\alpha}}{N^n} \sum\limits_{x \in \Z^n} \sum\limits_{\|d\|=\epsilon N}^{DN} \frac{1}{\|d\|^{n+\alpha}} gh_d(\eta_s(x)) \big[G_s\big(\frac{x+d}{N}\big)-G_s\big(\frac{x}{N}\big)\big] ds
\end{eqnarray*}
vanishes in probability. Here, and elsewhere, we write $\epsilon N$ and $DN$ for $\lceil \epsilon N\rceil$ and $\lfloor DN\rfloor$.

\vskip .2cm

{\flushleft \it Step 3a.}  We will now like to replace the nonlinear terms `$gh_d(\eta_s(x))$' by effective linear ones in terms of $\pi^N_s$.

The first replacement involves substituting $gh_d(\eta_s(x))$, with its average over $l$-blocks: $(gh_d)^l(\eta_s(x))$, where $l$ diverges after $N$ diverges, but before the limits on $\epsilon$ and $D$. By a discrete integration-by-parts, smoothness and compact support of $G$, the error introduced is of the expected order
$$\kappa\|h\|_\infty \|\nabla G\|\frac{ N^\alpha}{l N^n}E^N \int_0^t \sum_{\|d\|=\epsilon N}^{DN} \frac{1}{\|d\|^{n+\alpha}}  \sum_{|x|\leq (R+D)N}\eta_s(x) ds,$$
  which vanishes, noting Lemma \ref{particlebound}.

Therefore, we have, as the various parameters tend to their limits, that
\begin{eqnarray*}
&&\lt \pi^N_0 , G_0\rt + \int_0^t  \lt \pi^N_s, \partial_sG_s \rt ds   \\
&&\ \ \ + \int_0^t \frac{N^{\alpha}}{N^n} \sum\limits_{x \in \Z^n} \sum\limits_{\|d\|=\epsilon N}^{DN} \frac{1}{\|d\|^{n+\alpha}} (gh_d)^l\big(\eta_s(x))\big) \big[G_s\big(\frac{x+d}{N}\big)-G_s\big(\frac{x}{N}\big)\big] ds
\end{eqnarray*}
vanishes in probability.

\vskip .1cm

{\flushleft \it Step 3b.} Next, we perform what is usually called the `$1$-block' replacement. 
We would like to replace $(gh_d)^l(\eta_s(x))$ by $\Phi(\eta^l_s(x))\Psi(\eta^l_s(x+d))$, the `averaged' function of the local mass density.  That is, we wish to show
$$\limsup\limits_{l \to \infty}\limsup\limits_{N \to \infty}E^N\int_0^t\frac{N^\alpha}{N^n} \sum_{\|d\|=\epsilon N}^{DN} \frac{1}{\|d\|^{n+\alpha}}\Big |\sum_{x \in \Z^n}  H^1_{x,d,l}(\eta_s)\big[ G_s\big(\frac{x+d}{N}\big)-G_s\big(\frac{x}{N}\big)\big]\Big |ds=0 $$
where $H^1_{x,d,l}(\eta)=(gh_d(\eta(x)))^l - \Phi(\eta^l(x))\Psi(\eta^l(x+d))$. By discrete integration-by-parts and bounding $G(x/N)$ by ${\bf 1}(|x|\leq RN)\|G\|$, it will be enough to show that 
\begin{eqnarray*}
&&E^N\int_0^t \frac{N^\alpha}{N^n} \sum\limits_{\|d\|=\epsilon N}^{DN} \frac{1}{\|d\|^{n+\alpha}} \sum\limits_{|x| \leq (R+D_N} \big|H^1_{x,d,l}(\eta_s)\big| ds \ \ {\rm and \ }\\
&&\ \ \ \ \ \ \ \ \ \ \ \ \ \ E^N\int_0^t \frac{N^\alpha}{N^n} \sum\limits_{\|d\|=\epsilon N}^{DN} \frac{1}{\|d\|^{n+\alpha}} \sum\limits_{|x| \leq (R+D)N} \big| \tau_{-d}H^1_{x,d,l}(\eta_s)\big| ds] 
\end{eqnarray*}
both vanish.  This is proved as a consequence of Proposition \ref{1block_alpha<1} in Subsection \ref{1blocks}.

After this $1$-block replacement, we have $$\lt \pi^N_0 , G_0\rt + \int_0^t  \lt \pi^N_s, \partial_sG_s \rt ds +  \int_0^t \frac{1}{N^{2n}} \sum\limits_{x \in \Z} \sum\limits_{\|d\|=\epsilon N}^{DN} \Phi\big(\eta^l_s(x)\big)\Psi\big(\eta^l_s(x+d)\big) \nabla^{\alpha,d} G_s\big(\frac{x}{N}\big)   ds$$
vanishes in the limit, 
where 
$\nabla^{\alpha,d} G_s(x/N) = \|(d/N)\|^{-n-\alpha}\big[G_s\big((x+d)/N\big)-G_s\big(x/N\big)\big]$.

\vskip .1cm

{\flushleft \it Step 3c.} The final estimate is the so-called `$2$-blocks' replacement,
where $\eta^l_s(x)$ is replaced by $\eta^{\epsilon'N}_s(x)$ in terms of a parameter $\epsilon'$.  We will write $\epsilon'N$ instead of $\lfloor \epsilon'N\rfloor$ throughout. 
That is, we will like to show, as in order $N\uparrow\infty$, $\epsilon'\downarrow 0$ and $l\uparrow\infty$, that
\begin{eqnarray*}
&&E^N \int_0^t \frac{1}{N^{2n}} \sum\limits_{x \in \Z^n}\sum\limits_{\|d\|=\epsilon N}^{DN}   \big|\Phi(\eta^l_s(x))\Psi(\eta^l_s(x+d))  \\
&&\ \ \ \ \ \ \ \ \ \ \ \ \ - \Phi(\eta^{\epsilon'N}_s(x))\Psi(\eta^{\epsilon'N}_s(x+d))\big| \big|\nabla^{\alpha,d} G_s\big(\frac{x}{N}\big) \big| ds 
\end{eqnarray*}
vanishes.

We observe that $\nabla^{\alpha,d} G_s\left(\frac{x}{N}\right)$ is zero, unless $-(R+D)N \leq x \leq RN$, in which case it is bounded. Also, as $\Phi$ is Lipschitz and $\Psi$ is bounded by $\|h\|_\infty$, we have 
$$\big|\Phi(\eta^l_s(x))\Psi(\eta^l_s(x+d)) - \Phi(\eta^{\epsilon'N}_s(x))\Psi(\eta^{\epsilon'N}_s(x+d))\big| \leq H^2_{x,d,l, \epsilon' N}(\eta_s)$$ 
where 
$H^2_{x,d,l, \epsilon' N}(\eta) = \kappa \|h\|_\infty |\eta^l(x)-\eta^{\epsilon' N}(x)|+ \Phi(\eta^l(x))|\Psi(\eta^l(x+d))-\Psi(\eta^{\epsilon' N}(x+d))|$.
Then, to show the $2$-blocks replacement, it will enough to show
$$\lim_{l\uparrow\infty}\lim_{\epsilon'\downarrow 0}\lim_{N\uparrow\infty}
E^N\int_0^t \frac{1}{N^n} \sum\limits_{\|d\|=\epsilon N}^{DN} \frac{1}{N^n} \sum\limits_{|x| \leq (R+D)N } H^2_{x,d,l, \epsilon' N}(\eta_s)  ds  =0.$$
This is a consequence of Proposition \ref{2blocks_alpha<1} in Subsection \ref{2blocks}.

We now observe that an $\epsilon' N$-block is macroscopically small, and may written in terms of $\pi^N_s$ as follows:
$$\eta^{\epsilon' N}_s(x) = \left(\frac{2\epsilon' N}{2\epsilon' N+1}\right)^n \langle \pi^N_s, \iota_{\epsilon'}(\cdot - x/N)\rangle$$
where $\iota_{\epsilon'} = (2\epsilon')^{-n}{\bf 1}([-\epsilon', \epsilon']^n)$. 
Hence, after the $2$-blocks replacement, we have `closed the equation', that is
\begin{eqnarray*}
&&\lt \pi^N_0 , G_0\rt + \int_0^t  \lt \pi^N_s, \partial_sG_s \rt ds   \\
&&\ \ \ + \int_0^t \frac{1}{N^n} \sum\limits_{x \in \Z^n} \frac{1}{N^n} \sum\limits_{\|d\|=\epsilon N}^{DN} \Phi\Psi_d\left(\langle \pi^N_s,   \iota_{\epsilon '}(\cdot -x/N\rangle\right) \nabla^{\alpha, d}G_s\left(\frac{x}{N}\right) ds
\end{eqnarray*}
vanishes in probability as $N\uparrow\infty$ and $\epsilon'\downarrow 0$.

\vskip .2cm

{\it Step 4.}   
We may replace the Riemann sums with integrals limited by $\epsilon$ and $D$.  As $\Phi,\Psi$ are Lipschitz and $\Psi$ is bounded (cf. \eqref{LB}), and as $\nabla^{\alpha, d}G_s$ is smooth,
the error accrued
is of expected order 
$N^{-(n+1)} \int_0^t E^N \sum_{|x|\leq (R+D)N} \eta_s(x)ds$,
which vanishes 
by Lemma \ref{particlebound}.

Further, we may then replace the limits in the integrals by $0$ and $\infty$, respectively.  The error of this replacement, comparing to Riemann sums, vanishes by Lemma \ref{lemma23_cor}.

Hence, we obtain that
$$\lt \pi^N_0 , G_0\rt + \int_0^t  \lt \pi^N_s, \partial_sG_s \rt ds +  \int_0^t \int_{\R^n} \int_{[0, \infty)^n} \Phi\Psi_v\left(\langle \pi^N_s,\iota_{\epsilon '}(\cdot - u)\rangle \right) \nabla^{\alpha, v}G_s(u) dv du ds$$
converges to zero in probability as $N\uparrow\infty$ and $\epsilon'\downarrow 0$, where for $v\in \R^n$, $\Psi_v(f(u))= \Psi(f(u + v))$ and we recall
$\nabla^{\alpha, v}G_s(u) = \|v\|^{-n-\alpha}\big [G_s(u+v)-G_s(u)\big ]$.

\vskip .2cm

{\it Step 5.}  Now, according to Proposition \ref{PNtightness},
the measures $\{P^N_G\}$ are tight, with respect to uniform topology.  Let $\{N_k\}$ be a subsequence where the measures converge to a limit point $P^*$.  The function of $\pi$,
$$\lt \pi_0 , G_0\rt + \int_0^t  \lt \pi_s, \partial_sG_s \rt ds +  \int_0^t \int_{\R^n} \int_{[0, \infty)^n} \Phi\Psi_v\big(\langle \pi^N_s,\iota_{\epsilon '}(\cdot -u)\rangle\big) \nabla^{\alpha, v}G_s(u) dv du ds,$$
is continuous for each $\epsilon'>0$.  Then,
letting $N_k\uparrow\infty$, we recover that
\begin{eqnarray*}
&&\lim_{\epsilon'\downarrow 0}
P^*\Big(\Big|\lt \pi_0 , G_0\rt + \int_0^t  \lt \pi_s, \partial_sG_s \rt ds \nonumber \\
&&\ \ \ \ \ \  +\int_0^t \int_{\R^n} \int_{[0, \infty)^n} \Phi\Psi_v\big(\langle \pi_s,  \iota_{\epsilon '}(\cdot -u)\rangle\big) \nabla^{\alpha, v}G_s(u) dv du ds\Big| > \epsilon_0 \Big) = 0.
\end{eqnarray*}
\vskip .2cm

{\it Step 6.}
Now, we claim that $P^*$ is supported on 
on measures $\pi_s$ that are absolutely continuous with respect to Lebesgue measure, and so
$\pi_s= \rho(s,u)du$ for an $L^1_{loc}$ function $\rho(s,u)$.  Indeed, this follows,
under condition FEM, with the same proof given for zero-range processes on pages 73-75 of \cite{Kipnis}.  We also have $\langle \pi_0,G_0\rangle = \langle \rho_0, G_0\rangle$ from our initial conditions.
Hence, $\langle \pi_s,  \iota_{\epsilon'}(\cdot - u)\rangle =  (2\epsilon')^{-n}\int_{[-\epsilon',\epsilon']^n}\rho(s,u+v)dv$.   
 Note also that $P^*$-a.s. $$\limsup_{\epsilon' \to 0}  \Big| \frac{1}{(2\epsilon')^n}\int_{[-\epsilon',\epsilon']^n}\rho(s,u+v)dv-\rho(s,u)\Big|=0 \ \text{ $u$-a.e. and in $L^1_{loc}$.}$$
 
  By properties of $\Phi, \Psi$, we have
\begin{align*}
&|\Phi\Psi_v (\langle \pi_s,  \iota_{\epsilon'}(\cdot - u)\rangle) - \Phi\Psi_v(\rho(s,u))| \\
& \ \leq \kappa\|\psi\|_\infty|\langle \pi_s,  \iota_{\epsilon'}(\cdot - u)\rangle - \rho(s,u)| + \kappa \rho(s,u)|\Psi_v(\langle \pi_s,  \iota_{\epsilon'}(\cdot - u)\rangle) - \Psi_v(\rho(s,u))|.
\end{align*}
As $\Psi$ is bounded, the second term on the right-side is bounded by $2\kappa\|\psi\|_\infty \rho(s,u)$.
Note also that $\sup_{w\in \R^n} E_{P^*}\int_0^t |\langle \pi_s, {\bf 1}(|\cdot -w|\leq R)\rangle ds = \sup_{w\in \R^n} E_{P^*} \int_0^t \int_{|u-w|\leq R} \rho(s,u)duds<\infty$ by Lemma \ref{particlebound} and lower semi-continuity in $\pi$ of the associated mapping. 
Then, as $G$ has compact support, by the $L^1_{loc}$ convergence, and use of dominated convergence, 
we have, 
with respect to each limit point $P^*$, a.s.
\begin{eqnarray*}
&& \int_{\R^n} \rho_0(u)  G_0(u)du + \int_0^\infty \int_{\R^n} \rho(s,u) \partial_sG_s(u) \rt ds \\
 &&\ \ \ \ \ +  \int_0^t \int_{\R^n} \int_{[0, \infty)^n} \Phi\Psi_v\left(\rho(s,u)\right) \nabla^{\alpha, v}G_s(u) dv du ds =0.
 \end{eqnarray*}
Since $G$ has compact support in $[0,T)$ with respect to time, we may replace $t$ by $T$.

In other words, every limit point $P^*$ is supported on absolutely continuous measures, $\pi_s=\rho(s,u)du$, whose densities $\rho(s,u)$ are weak solutions of the hydrodynamic equation.  This concludes the proof of Theorem \ref{alpha<1}.
\qed

\section{$1$-block and $2$-block estimates}
\label{blocks}

We discuss the $1$ and $2$ block estimates when $\alpha<1$, and also a $1$-block estimate when $\alpha>1$ in the next three subsections.

\subsection{$1$-block estimate: $\alpha <1$}
\label{1blocks}

We now prove the $1$-block replacement used in Section \ref{proofsection_alpha<1}.  As a comment, in Step 3, due to the long range setting, we use a somewhat nonstandard estimate.

\begin{prop}[$1$-block estimate]
\label{1block_alpha<1}
When $\alpha<1$, we have
$$\limsup\limits_{l \to \infty}\limsup\limits_{N \to \infty} E^N\int_0^t \sum\limits_{\|d\|=\epsilon N}^{DN} \frac{ N^\alpha}{\|d\|^{n+\alpha}} \frac{1}{N^n} \sum\limits_{|x| \leq RN} \left|H^1_{x,d,l}(\eta_s)\right| ds = 0$$
where $H^1_{x,d,l}(\eta)=(gh_d)^l(\eta(x))) - \Phi\big(\eta^l(x)\big)\Psi\big(\eta^l(x+d)\big)$
\end{prop}

\bpf 
The proof goes through a few steps.

\vskip .1cm
{\it Step 1.}
We first introduce a truncation.  As both $|h|, |\Psi|\leq \|h\|_\infty$, and both $g, \Phi$ are Lipschitz, we can bound $H^1_{x, d, l}\leq 2\kappa\|h\|_\infty\eta^l(x)$.
Since $\limsup_{N\uparrow\infty}\sum_{\|d\|=\epsilon N}^{DN} N^\alpha/\|d\|^{n+\alpha}<\infty$,
 we can introduce the indicator function ${\bf 1}(\eta^l_s(x) \leq A)$ into the integrand in the display by Lemma \ref{particlebound2}.

We may introduce one more truncation.  As $|H^1_{x,d,l}(\eta_s)|{\bf 1}(\eta^l(x) \leq A)$ is bounded in terms of $A$, and again the sum $\sum_{\|d\|=\epsilon N}^{DN} N^\alpha/\|d\|^{n+\alpha}$ is uniformly bounded in $N$, we can introduce the indicator function ${\bf 1}(\eta^l(x+d) \leq A)$ in the integrand again by say Lemma \ref{particlebound2}.

  It will be enough to show, for each $A$, as $N\uparrow\infty$ and $l\uparrow\infty$, that the following vanishes,
$$
E^N\int_0^t \sum\limits_{\|d\|=\epsilon N}^{DN} \frac{ N^\alpha}{\|d\|^{n+\alpha}} \frac{1}{N^n} \sum\limits_{|x| \leq RN} \left|H^1_{x,d,l}(\eta_s)\right|{\bf 1}(\eta^l_s(x) \vee \eta^l_s(x+d)\leq A) ds.
$$

\vskip .2cm

{\it Step 2.}  Recall the density $\bar{f}^N_t$ in Subsection \ref{entropy_section}.
The expected value above equals
$$t \int \sum\limits_{\|d\|=\epsilon N}^{DN} \frac{ N^\alpha}{\|d\|^{n+\alpha}} \frac{1}{N^n} \sum\limits_{|x| \leq RN} \left|H^1_{x,d,l}(\eta)\right|{\bf 1}(\eta^l(x) \vee \eta^l(x+d) \leq A) \bar{f}^N_t(\eta) \nu_{\rho^*}(d\eta).$$

Given the Dirichlet bound on $\bar{f}^N_t$ in \eqref{dirichlet_bound} of order $N^n/N^\alpha$, we need only show that
$$\sup_{D(f)\leq \frac{C_0N^n}{N^\alpha}} \int \sum\limits_{\|d\|=\epsilon N}^{DN} \frac{ N^\alpha}{\|d\|^{n+\alpha}} \frac{1}{N^n} \sum\limits_{|x| \leq RN} \left|H^1_{x,d,l}(\eta)\right| {\bf 1}(\eta^l(x) \vee \eta^l(x+d) \leq A) f(\eta) \nu_{\rho^*}(d\eta)$$
vanishes as $N\uparrow\infty$ and $l\uparrow\infty$.

Let $f^{R,N}(\eta)= \frac{1}{(2RN+1)^n} \sum_{|x|\leq RN} \tau_x f(\eta)$.  By translation-invariance of $\nu_\rho$,
the above display equals
\begin{equation*}
\label{alpha<1_1block_step2}\sup_{D(f)\leq \frac{C_0N^n}{N^\alpha}} \int \sum\limits_{
\|d\|=\epsilon N}^{DN} \frac{ N^\alpha}{\|d\|^{n+\alpha}} \frac{(2RN+1)^n}{N^n}) \left|H^1_{0,d,l}(\eta)\right|{\bf 1}(\eta^l(0)\vee \eta^l(d) \leq A) f^{R,N}(\eta) \nu_{\rho^*}(d\eta).
\end{equation*}

 \vskip .2cm
 
{\it Step 3.}
At this stage, there is a trick that is not part of the standard $1$-block argument because, in $H^1_{0,d,l}$, we in fact have $2$ $l$-blocks, about $0$ and $d$. Let $\xi$ and $\zeta$ be configurations on $[-l,l]^n$ that equal $\eta$ and $\tau_d\eta$, respectively, on $[-l ,l]^n$. Then,
$$H^1_{0,d,l}(\eta) =  \frac{1}{(2l +1)^n}\sum\limits_{|y|\leq l} g(\xi(y)h(\zeta(y)) - \Phi(\xi^l(0))\Psi(\zeta^l(0))=: \mathcal{H}^1_l(\xi, \zeta).$$
Let $\nu^1_{\rho^*}(d\xi, d\zeta)$ be the product measure on pairs of configurations $(\xi, \zeta)$ induced by $\nu_{\rho^*}$, and let $\bar{f}_{l,d}(\xi, \zeta)$ be the conditional expectation of $f^{R,N}(\eta)$ given configurations $\eta$ that equal $\xi$ on $[-l,l]^n$ and $\zeta$ on $[-l+d, l+d]^n$. 

 Define now
 $$\bar{f}_{N,l, \epsilon, D} =  \Big(\sum\limits_{\|d\|=\epsilon N}^{DN} \frac{N^\alpha}{\|d\|^{n+\alpha}} \bar{f}_{l,d}\Big) \Big/ \sum\limits_{\|d\|=\epsilon N}^{DN} \frac{N^\alpha}{\|d\|^{n+\alpha}}.$$
  Given $(2RN+1)^n/N^n\leq (2R+1)^n$ and $\sum_{\|d\|=\epsilon N}^{DN} N^\alpha/\|d\|^{n+\alpha}$ is bounded in terms of $\epsilon$ and $D$, it will be sufficient to show that
$$\limsup\limits_{l\uparrow\infty}\limsup\limits_{N\uparrow\infty}
\sup\limits_{D(f)\leq \frac{C_0N^n}{N^\alpha}}  \int   \left|\mathcal{H}^1_l(\xi, \zeta)\right|{\bf 1}(\xi^l(0)\vee \zeta^l(0) \leq A) \bar{f}_{N, l, \epsilon, D}(\xi, \zeta) \nu^1_{\rho^*}(d\xi, d\zeta) = 0.$$

\vskip .2cm

{\it Step 4.}
  Let $D_1^{x,y}(f)=D^{x,y}(f)$ and $D_2^{x,y}(f)=D^{x,y}(\tau_d f)$ be the bond Dirichlet forms with respect to configurations $\xi$ and $\zeta$ respectively.  Define now a new Dirichlet form,
$$D^2_l(f):= \sum_{|y| \leq l} \left(D_1^{0,y}(f)+D_2^{0,y}(f)\right).$$ 
In Lemma \ref{1blockalpha<1} below, we prove the following bound 
$D^2_l(\bar{f}_{N,l, \epsilon, D}) \leq {C_1}/{N^\alpha}$.
Therefore, it will be enough to show 
$$\limsup\limits_{l\uparrow\infty}\limsup\limits_{N\uparrow\infty}\sup\limits_{D^2_l(f)\leq \frac{C_1}{N^\alpha}}  \int   \big|\mathcal{H}^1_l(\xi, \zeta)\big|{\bf 1}(\xi^l(0)\vee \zeta^l(0) \leq A) f(\xi, \zeta) \nu^1_{\rho^*}(d\xi, d\zeta)=0.$$

Because of the truncation, ${\bf 1}(\xi^l(0)\vee \zeta^l(0) \leq A)$, the domain of $f$ is a finite set of configurations $(\xi,\zeta)$. 
Hence, the condition $D^2_l(f)\leq C_1/N^\alpha$ specifies a compact set of densities, and 
the 
supremum in the display, for each $N$ and $l$, is attained at some density denoted $f_{N, l}$. As $N\uparrow\infty$, any convergent subsequence of $\{f_{N, l}\}$ approaches a function whose Dirichlet form vanishes. Therefore, after $N\uparrow\infty$, the supremum is bounded by 
$$ \sup\limits_{D^2_l(f)=0} \int  \left| \mathcal{H}^1_l(\xi, \zeta)\right|{\bf 1}(\xi^l(0)\vee \zeta^l(0) \leq A) f(\xi, \zeta) \nu^1_{\rho^*}(d\xi, d\zeta) .$$

\vskip .2cm

{\it Step 5.} If $D^2_l(f)=0$, then $f$ is a constant, namely $f=C_{j_1, j_2}$ along the hyperplanes 
$$H^2_{j_1, j_2} =\Big\{(\xi, \zeta): \sum_{|y|\leq l}\xi(y) = j_1, \sum_{|y|\leq l}\zeta(y) = j_2 \Big\},$$ 
where $j_1, j_2= 0, 1, ..., (2l+1)^nA$. Because $f$ is a probability density, the constants $C_{j_1, j_2}$ are nonnegative and $\sum C_{j_1, j_2} \nu^1_{\rho^*}\big(H^2_{j_1,j_2}\big)=1$. Therefore, the last expression can be written as a supremum over all possible weighted averages of the integral over hyperplanes:
$$\sup\limits_{\sum c_{j_1, j_2} =1}\sum\limits_{j_1, j_2 =0}^{(2l+1)^nA} c_{j_1, j_2} \int_{H^2_{j_1, j_2}}  \big| \mathcal{H}^1_l(\xi, \zeta)\big| \nu^{l,j_1,j_2}(d\xi, d\zeta),$$
where $\nu^{l,j_1,j_2}$ is the canonical product measure $\nu^1_{\rho^*}$ conditioned to $H^2_{j_1,j_2}$, no longer depending on $\rho^*$.

In particular, it will be enough to show that
$$ \limsup\limits_{l\uparrow\infty}\sup\limits_{0 \leq j_1, j_2 \leq (2l+1)^nA} \int  \big| \mathcal{H}^1_l(\xi, \zeta) \big| \nu^{l, j_1, j_2}(d\xi, d\zeta)=0.$$

\vskip .2cm
{\it Step 6.}
We try to make the integrand independent of $l$.  We may choose $l$ 
so that $(2l+1)^n = q(2k+1)^n$, that is, an $l$-block is partitioned into $k$-blocks.  
Let $B_1, .., B_q$ denote the $k$-blocks.  Then,
$$\mathcal{H}^1_l(\xi, \zeta) = \frac{1}{q} \sum\limits_{i=1}^q \frac{1}{(2k+1)^n}\sum\limits_{y \in B_i} \left(g(\xi(y)h(\zeta(y)))-\Phi(\xi^l(0))\Psi(\zeta^l(0))\right)$$
Under $\nu^{l,j_1, j_2}$, the distribution of $\sum_{y \in B_i} (g(\xi(y))h(\zeta(y))-\Phi(\xi^l(0)\Psi(\zeta^l(0)))/(2k+1)^n$ doesn't depend on $i$. Therefore, noting $\xi^l(0)=\frac{j_1}{(2l+1)^n}$ and $\zeta^l(0)=\frac{j_1}{(2l+1)^n}$, we can bound the previous supremum by 
\begin{eqnarray*}
&&\sup\limits_{0 \leq j_1, j_2 \leq (2l+1)^nA} \int \Big|\frac{1}{(2k+1)^n}\sum_{|y|\leq k}g(\xi(y))h(\zeta(y))\\
&&\ \ \ \ \ \ \ \ \ \ \ \ \ \ \ \  -\Phi\big(\frac{j_1}{(2l+1)^n}\big)\Psi\big(\frac{j_2}{(2l+1)^n}\big)\Big| \nu^{l, j_1, j_2}(d\xi, d\zeta).
\end{eqnarray*}

We can then take the limit as $l\uparrow\infty$, that is, as $q\uparrow\infty$ to arrive, by a local central limit theorem or equivalence of ensembles estimate as in Corollary 1.7 in Appendix 2 \cite{Kipnis}, at the expression
$$\sup\limits_{0\leq \rho_1, \rho_2 \leq A} \int  \big|\frac{1}{(2k+1)^n}\sum_{|y|\leq k}g(\xi(y))h(\zeta(y))-\Phi(\rho_1)\Psi(\rho_2 )\big| (\nu_{\rho_1} \times \nu_{\rho_2})(d\xi, d\zeta).$$
But, this quantity, say using a Chebychev bound, vanishes uniformly for $0\leq \rho\leq A$ as $k\uparrow\infty$ by the law of large numbers, since $\Phi(\rho_1)\Psi(\rho_2)=E_{\nu_{\rho_1}\times \nu_{\rho_2}}[g(\xi(0))h(\zeta(0))]$.

As a remark, this last step is rather interesting. Normally, the usual $1$-block estimate ends by showing that an average of a function of the $\xi(y)$ converges to its expected value. Here, in the $\alpha<1$ case, we end up with term that looks like a covariance. \epf

We now show a bound on $D^2_l(\bar{f}_{N,l, \epsilon, D})
 = \sum\limits_{|y|\leq l} \big( D_1^{0,y}(\bar{f}_{N, l, \epsilon, D}) + D_2^{0,y}(\bar{f}_{N, l, \epsilon, D})\big)$.

\blem \label{1blockalpha<1} Let $f$ be a density such that $D(f) \leq \frac{C_0N^n}{N^\alpha}$.  Then, $D^2_l(\bar{f}_{N,l,\epsilon, D}) \leq {C_1}/{N^\alpha}$. 
\elem 

\bpf By the convexity of the Dirichlet form, for $i=1,2$,
$$D_i^{0,y}(\bar{f}_{N, l, \epsilon, D}) \leq \frac{1}{\sum\limits_{\|d\|=\epsilon N}^{DN} \frac{N^\alpha}{\|d\|^{n+\alpha}}} \sum\limits_{\|d\|=\epsilon N}^{DN} \frac{N^\alpha}{\|d\|^{n+\alpha}}D_i^{0,y}(\bar{f}_{l, d}).$$
Moreover,
$D_1^{0,y}(\bar{f}_{l, d}) \leq D^{0,y}(f^{R,N})$ and $D_2^{0,y}(\bar{f}_{l, d}) \leq D^{d,y+d}(f^{R,N})$.
It follows
$D^2_l(\bar{f}_{N, l, \epsilon,D})$
is less than
$$\frac{1}{\sum\limits_{\|d\|=\epsilon N}^{DN} \frac{N^\alpha}{\|d\|^{n+\alpha}}} \sum\limits_{\|d\|=\epsilon N}^{DN} \frac{N^\alpha}{\|d\|^{n+\alpha}}  \frac{1}{(2RN+1)^n}  \sum\limits_{|x|\leq RN} \sum\limits_{|y|\leq l} \big( D^{0,y}(\tau_x f) + D^{d,y+d}(\tau_xf)\big).
$$

Noting $D^{0,y}(\tau_x f) =D^{x,y+x}(f)$, the last display is bounded by
$$\frac{(2RN+1)^{-n}}{\sum\limits_{\|d\|=\epsilon N}^{DN} \frac{N^\alpha}{\|d\|^{n+\alpha}}} \sum\limits_{\|d\|=\epsilon N}^{DN} \frac{N^\alpha}{\|d\|^{n+\alpha}} \sum\limits_{\substack{|y|\leq l\\ |x|\leq RN}} \big( D^{x,y+x}(f) + D^{d+x,y+d+x}(f)\big).$$

Note that  $\sum\limits_{\substack{|y|\leq l\\ |x|\leq RN}} D^{x,y+x}(f) + D^{d+x,y+d+x}(f) \leq 4D(f)$, as each bond is counted at most four times.  Hence, we may bound the last display further by
$$\frac{4D(f)}{(2RN+1)^n} \leq \frac{4}{(2RN+1)^n}\frac{C_0N^n}{N^\alpha} \leq \frac{C_1}{N^\alpha}$$
where $C_1= 2C_0/R^n$. \epf

\subsection{$2$-blocks estimate: $\alpha<1$}
\label{2blocks}

The proof of the $2$-blocks estimate is similar to the preceding $1$-block estimate, so we will give only a brief overview of the key differences.

\begin{prop}[$2$-blocks estimate]
\label{2blocks_alpha<1}
When $\alpha<1$, we have
 $$\limsup\limits_{l \to \infty} \limsup\limits_{\epsilon' \to 0^+} \limsup\limits_{N \to \infty}E^N\left[\int_0^t \frac{1}{N^n} \sum\limits_{\|d\|=\epsilon N}^{DN} \frac{1}{N^n}\sum\limits_{|x| \leq (R+D)N } H^2_{x,d,l, \epsilon' N}(\eta_s)  ds \right] =0, $$
where 
$$H^2_{x,d,l, \epsilon' N}(\eta) = \kappa\|h\|_\infty \big|\eta^l(x)-\eta^{\epsilon' N}(x)\big|+ \Phi(\eta^l(x))\big|\Psi(\eta^l(x+d))-\Psi(\eta^{\epsilon' N}(x+d))\big|
$$
\end{prop}

\bpf 
The proof uses several steps.
\vskip .2cm

{\it Step 1.}
Analogous to the $1$-block proof, we introduce a truncation. 
We can bound the second term of $H^2_{x,d,l,\epsilon'N}(\eta)$ as
$\Phi(\eta^l(x))|\Psi(\eta^l(x+d))-\Psi(\eta^{\epsilon' N}(x+d))| \leq 2\kappa\|h\|_\infty\eta^l(x)$, since  
$\Psi(\cdot)\leq \|h\|_\infty$ and $\Phi$ is Lipschitz.  
Given that $\sum_{\|d\|=\epsilon N}^{DN} N^{-n}$ remains bounded as $N\uparrow\infty$, 
we can introduce the indicator function ${\bf 1}(\eta^l_s(x) \leq A)$ onto the second term of $H^2_{x,d,l, \epsilon' N}(\eta)$ by Lemma \ref{particlebound2}. Since $\Psi$ is Lipschitz, the
truncated second term is less than $\kappa_1\Phi(A)\big|\eta^l(x+d)- \eta^{\epsilon' N}(x+d)\big|$.
The proposition will follow if we show, as $\epsilon'\downarrow 0$, $l\uparrow\infty$ and $N\uparrow\infty$, that
$$E^N\int_0^t \frac{1}{N^n} \sum\limits_{\|d\|=\epsilon N}^{DN} \frac{1}{N^n}\sum\limits_{|x| \leq (R+D)N } |\eta_s^l(x^*)-\eta_s^{\epsilon' N}(x^*)|  ds  =0, $$
for both $x^*=x$ and $x^*=x+d$.

As in the standard $2$-blocks estimate, we will replace an $\epsilon'N$ block, $\eta^{\epsilon' N}_s$, by an average of $l$-blocks, $\eta^l_s$.
Specifically, we will replace $\eta^{\epsilon' N}_s(x^*)$ by 
$$\frac{1}{(2\epsilon' N+1)^n} \sum_{2l < |y|\leq \epsilon' N } \eta^l_s(x^*+y).$$ 
The expected error introduced is of order $E^N\int_0^t N^{-2n}\sum_{|x|\leq R'}\eta_s(x)ds$, for some $R'$, which vanishes by say Lemma \ref{particlebound}.

By bounding the `average' over $y$ by a supremum, it will be enough to show that
$$\sup_{2l< |y| \leq \epsilon' N} E^N\int_0^t \frac{1}{N^n} \sum\limits_{\|d\|=\epsilon N}^{DN} \frac{1}{N^n} \sum\limits_{|x|\leq (R+D)N}\big|\eta^l_s(x^*)- \eta^l_s(x^*+y)\big| ds$$
vanishes.

\vskip .2cm

{\it Step 2.} From here, the proof of the 2-blocks estimates proceeds in the same way as for the 1-block estimate. We can write the expected value in terms of $\bar{f}^N_t$ and then majorize by a factor $t$ times
$$ \sup_{2l< |y| \leq \epsilon' N} \sup_{D(f)\leq \frac{C_0N^n}{N^\alpha}} \int \frac{1}{N^n} \sum\limits_{\|d\|=\epsilon N}^{DN} \frac{1}{N^n} \sum\limits_{|x|\leq (R+D)N}\left|\eta^l(x^*)- \eta^l(x^*+y)\right| f(\eta) \nu_{\rho^*}(d\eta).$$

Suppose $x^*=x+d$. As $\nu_{\rho^*}$ is translation invariant, we may replace $x+d$ by $x$ and $f(\eta)$ by $\tau_{-d}f(\eta)$. Now, $\tau_{-d}f(\eta)$ and therefore $N^{-n}\sum_{\|d\|=\epsilon N}^{DN} \tau_{-d}f$, by convexity, satisfy the same entropy and Dirichlet form bounds as $f$.  Hence, without loss of generality, we may assume $x^*=x$ in the sequel. 

In this case, $\sum_{\|d\|=\epsilon N}^{DN} N^{-n}$ can then be pulled out, bounded above by a constant. It will be enough to show that
$$ \sup_{2l< |y| \leq \epsilon' N} \sup_{D(f)\leq \frac{C_0N^n}{N^\alpha}} \int \frac{1}{N^n} \sum\limits_{|x|\leq (R+D)N}\left|\eta^l(x)- \eta^l(x+y)\right| f(\eta) \nu_{\rho^*}(d\eta)$$
which looks like the standard $2$-blocks estimate, say in \cite{Kipnis}.

\vskip .2cm

{\it Step 3.}
We may introduce the indicator function ${\bf 1}(\eta^l(x) \vee \eta^l(x+y) \leq A)$ to the integrand by Lemma \ref{particlebound2}.
By translation-invariance of $\nu_{\rho^*}$, we can shift the summand by $\tau_{-x}$.  Recall the averaged density $f^{R+D,N}$, introduced in Step 2 in Subsection \ref{1blocks}.
Multiplying and dividing by $2(R+D)N+1)^n$ and noting that the factor $(2(R+D)N+1)^n/N^n$ is bounded, by convexity of the Dirichlet form, it will be enough to show the following vanishes:
$$ \sup_{2l< |y| \leq \epsilon' N} \sup_{D(f)\leq \frac{C_0N^n}{N^\alpha}} \int \left|\eta^l(0)- \eta^l(y)\right|{\bf 1}(\eta^l(0) \vee \eta^l(y) \leq A) f^{R+D,N}(\eta) \nu_{\rho^*}(d\eta).$$

\vskip .2cm

{\it Step 4.} Let $\xi_1$ and $\xi_2$ be configurations on $[-l,l]^n$, equal to $\eta$ and $\tau_y\eta$, respectively. Let $\nu^2_{\rho^*}(d\xi_1, d\xi_2)$ be the associated induced measure with respect to $\nu_{\rho^*}$.
Let also $\bar{f}_{l, y}(\xi_1, \xi_2)$ be the conditional expectation of $f^{R+D,N}(\eta)$ given configurations $\eta$ that equal $\xi_1$ on $[-l,l]^n$ and $\xi_2$ on $[-l+y, l+y]^n$.  The last display in Step 3 equals
$$ \sup_{2l< |y| \leq \epsilon' N} \sup_{D(f)\leq \frac{C_0N^n}{N^\alpha}} \int \left|\xi^l_{\bf 1}(0)- \xi^l_2(0)\right|{\bf 1}(\xi^l_{\bf 1}(0) \vee \xi^l_2(0) \leq A) \bar{f}_{l,y}(\xi_1, \xi_2) \nu^2_{\rho^*}(d\xi_1, d\xi_2).$$

With $D_1^{w,z}(f) = D^{w,z}(f)$ and $D_2^{w,z}(f) = D^{w,z}(\tau_y f)$, we now introduce a Dirichlet form, 
$$D^*_{l}(f)=\sum_{|x| \leq l} \left(D_1^{0, x}(f)+ D_2^{0, x}(f)\right) + \frac{D^*(f)}{p^{sym}(y)}, \ \ \ {\rm where}$$
$$D^*(f)/p^{sym}(y) = D^{0,y}(f)/p^{sym}(y) = E_{\nu_{\rho^*}}\big[ g(\eta(0))h(\eta(y))(\sqrt{f}(\eta^{0,y}) - \sqrt{f}(\eta))^2\big]$$ 
is the `unit jump rate' Dirichlet form on the bond between the centers of the $l$-blocks involved.
Dividing by $p^{sym}(y)$ ensures, no matter the size of $y$, that a zero form $D_l^*(f)=0$ implies that $f$ is invariant to particle motion within each $l$-block and also motion between the centers.  In this case, $f$ takes a constant value along each of the hyperplanes 
$$H^2_j =\big\{(\xi_1, \xi_2): \sum_{|y|\leq l}(\xi_{\bf 1}(y)+\xi_2(y)) = j\big\}$$ 
for $j= 0, 1, ..., 2(2l+1)^n A$.

In Lemma \ref{2blockalpha<1} at the end of the Subsection, for $2l<|y|\leq \epsilon' N$, we prove the bound 
$D^*_l(\bar{f}_{l,y}) \leq {C_2}{N^{-\alpha}}+C_3(\epsilon')^\alpha$.
Therefore, it will be enough to show the following vanishes:
$$\sup_{D^*_l(f)\leq \frac{C_2}{N^\alpha}+C_3(\epsilon')^\alpha} \int \left|\xi^l_{\bf 1}(0)- \xi^l_2(0)\right|{\bf 1}(\xi^l_{\bf 1}(0) \vee \xi^l_2(0) \leq A) f(\xi_1, \xi_2) \nu^2_{\rho^*}(d\xi_1, d\xi_2).$$

As in the $1$-block proof, as particle numbers are bounded, we may take limits, as $N\uparrow\infty$ and $\epsilon'\downarrow 0$, to restrict the supremum above to densities $f$ such that $D^*_l(f)=0$.

\vskip .2cm
{\it Step 5.}
Hence, at this stage, $f$ equals a constant $C_j$ along each hyperplane $H^2_j$ for $j\leq 2(2l+1)^nA$.  Because $f$ is a probability density, these constants $C_j$ are non-negative and $\sum_j C_j \nu^2_{\rho^*}(H^2_j) = 1$. Therefore, 
we need only show
$$\sup_{0 \leq j \leq 2(2l+1)^nA} \int \left|\xi^l_{\bf 1}(0)- \xi^l_2(0)\right| \nu^{2, l, j}(d\xi, d\zeta)$$
vanishes, where $\nu^{2,l, j}$ is the canonical measure on configurations $(\xi_1, \xi_2)$ which distributes $j$ particles among the two $l$-blocks.
 
However, both $\xi^l_1(0)$ and $\xi^l_2(0)$ equal $\frac{j}{2(2l+1)^n}$.  Hence, adding and subtracting $\frac{j}{2(2l+1)^n}$ inside the absolute value, it will be enough to control
$\var_{\nu^{2, l,j}}(\xi_i^l(0)) \leq C E_{\nu^{2,l,j}}[l^{-n}(\bar\xi_1(0))^2 + 2\bar\xi_1(0)\bar\xi_1(e_1)]$, $\bar\xi(x) = \xi(x) - j/[2(2l+1)^n]$ and $e_1 = (1,0,\ldots, 0)$.  By the equivalence of ensembles as used in Step 6 of Proposition \ref{1block_alpha<1}, noting $\nu_{j/[2(2l+1)^n]}$ is a product measure with identical marginals, the variance vanishes as $l\uparrow\infty$. \epf

\subsubsection{Moving particle lemma}

We now prove the following bound on $D^*_{l, y}(\bar{f}_{l,y})$.  Part of the strategy is inspired by \cite{Jara} where a similar `moving particle' estimate was proved.  Recall that  $p^{sym}(y)= \big(p(y)+p(-y)\big)/2$, which equals $p(|y|)/2$ in our totally asymmetric setting.

\blem\label{2blockalpha<1} Suppose $2l < |y| \leq \epsilon' N$ and $D(f) \leq \frac{C_0N^n}{N^\alpha}$.   Then,  
$$D^*_l(\bar{f}_{l,y})\leq \frac{C_2}{N^\alpha}+ C_3 (\epsilon')^\alpha.$$ 
\elem

\bpf First, by the same argument as in Lemma \ref{1blockalpha<1}, the sum
$$\sum_{|x| \leq l} \big(D_1^{0, x}(\bar{f}_{l,y})+D_2^{0, x}(\bar{f}_{l,y})\big) \leq C_2/N^\alpha, \ \ {\rm where \ }C_2= C_0/(R+D)^n.$$

Therefore, we need to control the form $D^*(\bar{f}_{l,y})/p^{sym}(y)$, which reflects motion from $0$ to $y=(y_1,\ldots, y_n)$, where $y_i\geq 0$.  
By the definition of $D^*$ and convexity of the Dirichlet form,
$$\frac{D^*(\bar{f}_{l,y})}{p^{sym}(y)} \leq  \frac{D^{0, y}(f^{(R+D)N})}{p^{sym}(y)}  \leq \frac{1}{(2(R+D)N+1)^n} \sum\limits_{|z|\leq (R+D)N} \frac{D^{z, z+y}(f)}{p^{sym}(y)}.$$
The term $D^{z, z+y}(f)/p^{sym}(y)$, reflecting a displacement by $y$, is now split up into several jumps of various lengths using Jensen's inequality and the structure of $\nu_{\rho^*}$. We will take advantage that, by definition, $p^{sym}$ is supported on all displacements.  The decomposition here is simpler than that which was used in \cite{Jara}.

We will split $D^{z, z+y}(f)/p^{sym}(y)$ into two jumps, one displacing by $k=(k_1,\ldots , k_n)$ where $0\leq k_i$, $|k|\leq |y|$, and $0 \neq k \neq y$, and one displacing by $y-k$. In the following, we will call $|y|_n = (|y|,\ldots, |y|)$.

By Jensen's inequality, and $p^{sym}(\cdot)/p^{sym}(|\cdot|_n)\geq 1$,
$$\frac{D^{z, z+y}(f)}{p^{sym}(y)} \leq 2 \frac{D^{z, z+k}(f)}{p^{sym}(|k|_n)} + 2\frac{D^{z+k, z+y}(f)}{p^{sym}(|y-k|_n)} = 4  \frac{D^{z, z+k}(f)}{p(|k|_n)} + 4\frac{D^{z+k, z+y}(f)}{p(|y-k|_n)}.$$

Therefore, 
$$
D^*(\bar{f}_{l,y})/p^{sym}(y) \leq
\frac{1}{(2(R+D)N+1)^n} \sum\limits_{|z|\leq (R+D)N} 4 \Big( \frac{D^{z, z+k}(f)}{p(|k|_n)} + \frac{D^{z+k, z+y}(f)}{p(|y-k|_n)}\Big)$$
which is less than
\begin{eqnarray*}
&& \frac{2^{2-n}}{((R+D)N)^np(|k|_n)} \sum\limits_{z}  D^{z, z+k}(f)\\
&&\ \ \ \ \ \ \  + \frac{2^{2-n}}{((R+D)N)^np(|y-k|_n)} \sum\limits_{z}  D^{z, z + y-k}(f) \\
&&\ \ \leq  \frac{2^{2-n}}{((R+D)N)^n \min\{p(|k|_n), p(|y-k|_n)\}} \sum\limits_{z}  D^{z, z+k}(f) +D^{z,z+y-k}(f). 
\end{eqnarray*}

Let 
\begin{eqnarray*}
a_k &=& \sum_z D^{z,z+k}(f) + D^{z,z+y-k}(f)
\ \ {\rm and }\\ 
b_k &= & \frac{((R+D)N)^n \min\{p(|k|_n), p(|y-k|_n)\}}{2^{2-n}} \geq \frac{((R+D)N)^n}{2^{2-n}}p(|y|_n),
\end{eqnarray*}
as $|y|\geq |k|, |y-k|$ and $p$ is decreasing in $|\cdot|_n$.
Then, 
$b_k D^*(\bar{f}_{l,y})/{p^{sym}(y)} \leq a_k$.

When $a_k$ is summed over all $k$ such that $0\leq k_i$, $|k|\leq |y|$, and $0 \neq k \neq y$, each bond is counted at most twice, and so
$\sum_k a_k  \leq 2D(f)$.
In particular, $[D^*(\bar f_{l,y})/p^{sym}(y)]\sum_k b_k \leq 2D(f)$.  Calculating
$\sum_k b_k \geq (R+DN)^n2^{n-2}p(|y|_n)\big[(|y| + 1)^n -2\big]$, we have
$$\frac{D^{*}(\bar{f}_{l,y})}{p^{sym}(y)} \leq \frac{2D(f)}{\sum\limits_{k} b_k} \leq \frac{4D(f)/((R+D)N)^n}{(|y|+1)^n-2)p(|y|_n)} = \frac{4C_0}{(R+D)^nN^\alpha} \frac{\||y|_n\|^{n+\alpha}}{(|y|+1)^n-2}.$$

As $\||y|_n\|$ scales like $|y|$, we have
$\||y|_n\|^{n+\alpha}/{[(|y|+1)^n-2]} =O((\epsilon' N)^\alpha)$. Therefore,
$$\frac{D^{*}(\bar{f}_{l,y})}{p^{sym}(y)}\leq \frac{O\big((\epsilon'N)^\alpha\big)}
{(R+D)^n N^\alpha} = O\big((\epsilon')^\alpha\big),$$
which gives the desired estimate. \epf

\subsection{$1$-block estimate: $\alpha>1$}
\label{1blocks_alpha>1}

The proof of the 1-block estimate, when $\alpha>1$, is similar to that when $\alpha<1$, but with fewer complications.  The argument is also similar to that in the standard finite-range setting in \cite{Kipnis}.  For completeness,
 we summarize the proof.

\begin{prop}[$1$-block estimate]
\label{1block_prop_alpha>1}
When $\alpha>1$, for each $d\in \Z^n$,
$$\limsup_{l \to \infty} \limsup_{N \to \infty} E^N\int_0^t \frac{1}{N^n}\sum\limits_{|x| \leq RN}\big| (gh_d)^l(\eta_s(x))-\Phi\Psi(\eta^l_s(x))\big|ds = 0.$$
\end{prop}

\bpf 
Following the proof of Proposition \ref{1block_alpha<1}, for $\alpha<1$,
we may introduce the indicator function ${\bf 1}(\eta^l_s(x)\vee \eta^l_s(x+d) \leq A)$, and 
bound the expectation in the display
by
\begin{equation}
\label{alpha>1_1block}
\sup\limits_{D(f)\leq \frac{C_0N^n}{N}} \int \left|(gh_d)^l(\eta(0))- \Phi\Psi(\eta^l(0))\right|{\bf 1}(\eta^l(0)\vee \eta^l(d) \leq A)f^{R,N}(\eta) \nu_{\rho^*} (d\eta).
\end{equation}

Let $\nu^{l,d}_{\rho^*}(d\xi)$ be the induced distribution of configurations $\xi$ equal to $\eta$ on $[-l, l]^n \cup d+[-l, l]^n$.  Let $\bar{f}_{l,d}(\xi)$ be the conditional expectation of $f^{R,N}(\eta)$ given configurations $\eta$ that equal $\xi$ on $[-l,l]^n\cup d+[-l, l]^n$. Introduce a Dirichlet form, 
$$D_{l,d}(f)= \sum_{y \in [-l,l]^n\cup d+[-l, l]^n} D^{0,y}(f).$$ 
In Lemma \ref{1blockalpha>1} below, when $D(f) = O(N^{n-1})$, we show that
$D_{l, d}(\bar{f}_{l, d})\leq {C_1}/{N}$.

Therefore, we can replace the supremum in \eqref{alpha>1_1block} by that over densities $f$ such that
$D_{l,d}(f)\leq C_1/N$.
As the truncation enforces a compact configuration space, after $N\uparrow\infty$, the supremum may be further replaced by $D_{l,d}(f)=0$. In this case, $f$ will be a constant $C_j\geq 0$ on hyperplanes of the form
$$H_j =\big\{\xi: \sum_{y \in [-l,l]^n \cup d+[-l, l]^n}\xi(y) = j\big\}.$$ 
As
$|[-l,l]^n\cup d+[-l, l]^n|^{-1} \sum_{y \in [-l,l]^n \cup d+[-l, l]^n} \xi(y)\leq \xi^l(0)+\xi^l(d)\leq 2A$,
the index $j \leq 2A|[-l,l]^n\cup d+[-l, l]^n|$.
Moreover, as $\sum_j C_j \nu^{l,d}_{\rho^*}(H_j) = 1$, 
we may bound \eqref{alpha>1_1block}
by a supremum over hyperplanes:
$$\sup\limits_{j}  \int   \big|(gh_d)^l(\xi(0))- \Phi\Psi(\xi^l(0))\big| \nu^{l,d,j}(d\xi),$$
where $\nu^{l,d,j}$ is the canonical measure 
supported on the hyperplane $H_j$.

As before, in the proof of Proposition \ref{1block_alpha<1}, we can partition $[-l, l]^n$ into $k$-blocks assuming $(2l+1)^n = q(2k+1)^n$. Let $B_1, .., B_q$ be the $q$ number of $k$-blocks.  Then,
$$(gh_d)^l(\xi(0)) - \Phi\Psi(\xi^l(0)) = \frac{1}{q} \sum\limits_{i=1}^q \frac{1}{(2k+1)^n}\sum\limits_{y \in B_i} \big(g(\xi(y))h(\xi(y+d)))-\Phi\Psi(\xi^l(0))\big).$$
Under the measure $\nu^{l, d, j}$, the distributions of $\sum_{y \in B_i} \big(g(\xi(y))h(\xi(y+d)))-\Phi\Psi(\xi^l(0))\big)$ do not depend on $i$. Therefore, it is enough to show
$$\limsup\limits_{k\rightarrow\infty}\limsup\limits_{l\rightarrow\infty}\sup_j  E_{\nu^{l,d,j}}  \big|   (gh_d)^k(\xi(0))- \Phi(\xi^l(0))\big|.$$

Now, we would like to replace $\Phi\Psi(\xi^l(0))$ by $\Phi\Psi(\rho)$ where $\rho = j/|[-l,l]^n\cup d+[-l, l]^n|$.  This should hold because $[-l, l]^n$ and $d+[-l, l]^n$ will have sufficient overlap for large $l$. To make this precise, 
bound $|\Phi\Psi(\xi^l(0))-\Phi\Psi(\rho)|\leq C(A)|\xi^l(0)-\rho|$ since $\Phi, \Psi$ are Lipschitz and $\xi^l(0)\leq 2A$. As in Step 6 of the proof of Proposition \ref{1block_alpha<1}, we can then split the $l$-block into $k$-blocks to show $\sup_j E_{\nu^{l,d,j}}|\xi^l(0)-\rho|$ vanishes as $l$ and then $k$ go to infinity, by an equivalence of ensembles estimate. 

Therefore, we only need to show
$\sup_j  E_{\nu^{l,d,j}}  \big|   (gh_d)^k(\xi(0))- \Phi\Psi(\rho)\big|$
goes to zero as $l$ and then $k$ go to infinity where $\rho = j/|[-l,l]^n\cup d+[-l, l]^n|$. By an equivalence of ensembles estimate, as in Step 6 of Proposition \ref{1block_alpha<1}, we need only show $E_{\nu_{\rho}} |  (gh_d)^k(\xi(0))- \Phi\Psi(\rho)|$ vanishes uniformly over bounded $\rho$ as $k\uparrow\infty$.  But, as $\nu_\rho$ is a product measure with identical marginals, and $E_{\nu_{\rho}}[gh_d(\xi(0))] = \Phi\Psi(\rho)$, this follows by a law of large numbers.  
\epf

We now prove the bound on $D_{l,d}(\bar{f}_{l, d})$.  Although the argument is similar to a finite-range setting estimate in \cite{Kipnis}, as it is short, we include it for convenience of the reader.

\blem \label{1blockalpha>1} Suppose $D(f) \leq \frac{C_0N^n}{N}$.  Then, $D_{l,d}(\bar{f}_{l,D}) \leq \frac{C_1}{N}$.
  \elem 

\bpf 
By convexity of the Dirichlet form,
$$
D_{l,d}(\bar{f}_{l,d}) = \sum\limits_{y \in [-l,l]^n\cup d+[-l, l]^n} D^{0,y}_{l,d}(\bar{f}_{l,d}) \leq
  \sum\limits_{y \in [-l,l]^n\cup d+[-l, l]^n} D^{0,y}(f^{R,N})$$
  which is less than
\begin{eqnarray*}
&& \sum\limits_{y \in [-l,l]^n\cup d+[-l, l]^n} \frac{1}{(2RN+1)^n} \sum\limits_{|x|\leq RN} D^{0,y}(\tau_x f)\\
&&\ \  \leq 
\frac{1}{(2RN+1)^n} \sum\limits_{\substack{y \in [-l,l]^n\cup d+[-l, l]^n\\ |x|\leq RN}} D^{x,y+x}(f).
\end{eqnarray*}
As bonds may be repeated twice in the sum in the last expression, we may bound by twice the full Dirichlet form, and obtain as desired
$D_{l,d}(\bar f_{l,d}) \leq \frac{2}{(2RN+1)^n} D(f) \leq {C_1}/{N}$,
where $C_1 = {2C_0}/{R^n}$. \epf

\section{Proof outline:  Hydrodynamic limits when $\alpha>1$}
\label{proof_section_alpha>1}

When $\alpha >1$, as the expected jump size $\sum_{d}dp(d)$ is finite, one expects in Euler scale to recover a similar hydrodynamic equation as when the jumps have finite range.  The strategy employed here is to follow the scheme of arguments in \cite{Rezakhanlou} and Chapter 8 in \cite{Kipnis} for finite-range processes.

However, in the long-range setting, several important steps are different.  In particular, we have worked to remove reliance on `attractiveness', a monotonicity condition on the rates, although it is still used in two, albeit, important places, namely to bound the hydrodynamic density as an $L^\infty$ object in Step 1 below, and to show the `Ordering' Lemma \ref{orderinglemma}, which is used to prove a measure weak entropy formulation.  On the other hand, the proof includes new arguments to bound uniformly the `mass difference from $\rho^*$' in the system, and to handle the `initial boundary layer' estimate, needed to apply a form of DiPerna's uniqueness characterization.

The first step in the argument is to use a $1$-block replacement estimate.  Here, we do not rely on `attractiveness' as in \cite{Rezakhanlou}, but the `entropy' method.  Part of the reason for this choice, as discussed in Subsection \ref{remarks}, is that, when $1<\alpha<2$, it is not clear how to use the `$L^1$-initial density' method in \cite{Rezakhanlou}. However, an artifact of using the `entropy' method is that we need to start from initial profiles $\rho_0$, which are close to $\rho^*$ at large distances.  

Since a `$2$-blocks' estimate is not available in the general asymmetric model, as also in \cite{Rezakhanlou} and \cite{Kipnis}, we use the concept of Young measures and DiPerna's characterization of measure-valued weak entropy solutions of the hydrodynamic equation to finish.

In terms of the process $\eta_t$, define a collection of Young measures as 
$$\pi^{N, l}_t(du, d\lambda) = \frac{1}{N^n}\sum_{x\in \Z^n} \delta_{\frac{x}{N}}(du)\delta_{\eta^l_t(x)}(d\lambda).$$
Integration with respect to $\pi^{N,l}_t$ against test functions is as follows:
$$\lt \pi^{N,l}_t, H(t,u, \lambda) \rt = \frac{1}{N^n} \sum_{x \in \Z^n}  H\big(t, \frac{x}{N}, \eta^l_t(x)\big)$$  
Denote by $\{Q^{N, l}\}$ the induced measures for the process $\{\pi^{N, l}_t: t\in [0,T]\}$.

Let $M_+(\R^n \times [0, \infty))$ be the set of positive Radon measures on $\R^n \times [0, \infty)$. Define $L^\infty([0, T], M_+(\R^n \times [0, \infty)))$ to be the space of functions $\pi_t:[0, T] \to M_+(\R^n \times [0, \infty))$ such that $\lt \pi_t, F \rt$ is essentially bounded in time for every continuous function $F$ with compact support in $\R^n \times [0, \infty)$.  The topology on $L^\infty([0, T], M_+(\R^n \times [0, \infty)))$ is such that elements $\pi$ and $\bar\pi$
are close if they give similar values upon integrating against a dense collection of test functions over space, $\lambda$, and time, that is if
$$ \int_0^T \lt \pi_s(du, d\lambda), F(u, \lambda) \rt H(s) ds \sim \int_0^T \lt \bar\pi_s(du, d\lambda), F(u, \lambda) \rt H(s) ds.$$
More precisely, the distance between $\pi$ and $\bar\pi$ is
$$d_{MV}(\pi,\bar{\pi}) = \sum_{k\geq 1} \frac{1}{2^k} \frac{d\big(\langle \pi, F_k\rangle, \langle \bar\pi, F_k\rangle\big)}
{1+ d\big(\langle \pi, F_k\rangle, \langle \bar\pi, F_k\rangle\big)},$$
where $\{F_k\}_{k\geq 1}$ is a dense sequence in the space of compactly supported functions in $\R^n\times [0,\infty)$, with respect to the uniform topology.  Here,
$$d(f,g) = \sum_{k\geq 1} \frac{1}{2^k} \frac{\big | \int_0^T dt h_k(t)g(t) - \int_0^T dt h_k(t)f(t)\big|}{1+ \big| \int_0^T dt h_k(t)g(t) - \int_0^T dt h_k(t)f(t)\big|},$$
where $\{h_k\}_{k\geq 1}$ is a dense sequence of functions in $L^1[0,T]$ (cf. p. 200 in \cite{Kipnis}).

 Note now that $\pi^{N,l}_t \in L^\infty([0, T], M_+(\R^n \times [0, \infty)))$, and accordingly $\{Q^{N,l}\}$ are measures on $L^\infty([0, T], M_+(\R^n \times [0, \infty)))$.  The general strategy, as in \cite{Rezakhanlou}, is to characterize limit points $Q^*$ of $\{Q^{N,l}\}$ in terms of unique `measure weak' solutions to the hydrodynamic equation.
 
 At this point, we remark, in later development, functions $F(s,u,\lambda)=G(s,u)f(\lambda)$ where $f$ is not compactly supported, but bounded $|f(\lambda)|\leq C\lambda$ for all large $\lambda$, will have use.  Although such functions are not part of the topology on $L^\infty([0, T], M_+(\R^n \times [0, \infty)))$, we establish in Subsection \ref{cutoff_sect}, for a subsequence $\{Q^{N',l'}\}$ converging to $Q^*$, that
\begin{eqnarray}
&&\lim_{l'\rightarrow\infty}\lim_{N'\rightarrow\infty} Q^{N',l'} \Big( \big| \int_0^T \langle \pi^{N',l'}_s, G(s,u)f(\lambda)\rangle ds \big|>\epsilon_0\Big) = 0 \ \ {\rm for \ all \ } \epsilon_0>0\nonumber\\
&&\Longleftrightarrow \ Q^*\Big(\big|\int_0^T \langle \pi_s, G(s,u)f(\lambda)\rangle ds\big|>\epsilon_0\Big) = 0 \ \ {\rm for \ all \ }\epsilon_0>0. 
\label{cutoff_stat}
\end{eqnarray}

\medskip
We now define the notion of `measure weak' solution.
Consider the weak formulation of the differential equation in terms of a weak solution $\rho(s,u)$. The measure weak formulation is obtained by replacing  $\rho(s,u)$ where ever it appears by $\lambda$ and then integrating against the measure $\rho(s, u, d\lambda)$ with respect to $\lambda$.  So, $f(\rho(s,u))$ becomes $\int f(\lambda)\rho(s,u,d\lambda)$.   
If $\rho(s, u, d\lambda)$ is a solution of the resulting equation, it is called a measure weak solution. For example,
\begin{eqnarray*}
&&\int_{\R^n} G_0(u) \int_\R \lambda \rho(0,u, d\lambda) du + \int_0^\infty \int_{\R^n} \partial_s G_s(u) \int_\R \lambda \rho(s,u, d\lambda)  du ds\\
&&\ \ \ \ \ \ \ \ \ \ \ \ \ \ \ \  + \gamma_\alpha \int_0^\infty \int_{\R^n} \int_\R \partial_{\textbf{1}(n)}G_s(u)\Phi\Psi(\lambda)\rho(s,u, d\lambda) duds=0
\end{eqnarray*}
is the measure weak formulation of nonnegative solutions of the hydrodynamic equation $\partial_t \rho + \gamma_\alpha \partial_{\textbf{1}(n)} \Phi(\rho)\Psi(\rho)=0$. 

We recall now part of Kru$\breve{\text{z}}$kov's entropy condition \eqref{entropy_solution_defn} on a weak solution of the hydrodynamic equation: For $c\in \R$,
$$\partial_t |\rho-c|+\gamma_\alpha \partial_{\textbf{1}(n)} [\text{sgn}(\rho-c)(\Phi\Psi(\rho)-\Phi\Psi(c))] \leq 0.$$
  It is known that there is a unique bounded weak solution of the hydrodynamic equation which satisfies Kru$\breve{\text{z}}$kov's entropy condition, with bounded initial data $w_0$ (cf. \cite{Kruz}, \cite{DiBen}, \cite{Evans}). 
The corresponding measure weak formulation
is given by 
\begin{eqnarray*}
&&\int_{\R^n} G_0(u) \int_\R |\lambda - c|\rho(0,u,d\lambda)du + \int_0^\infty \int_{\R^n}  \partial_sG_s\left( u\right) \int_\R |\lambda - c| \rho(s,u, d\lambda) du ds\\
&&\ \ \ \ \ \ \ \ + \gamma_\alpha\int_0^\infty \int_{\R^n} \partial_{\textbf{1}(n)} G_s\left(u\right) \int_\R q(\lambda,c)\rho(s,u, d\lambda)duds \geq 0
\end{eqnarray*}
where $q(\lambda,c)= \text{sgn}(\lambda-c)(\Phi\Psi(\lambda)-\Phi\Psi(c))$ and $G$ is a nonnegative test function. 
We will say that $\rho(t, u, d\lambda)$ is a measure weak entropy solution if it is a measure weak solution of the hydrodynamic equation that measure weakly satisfies the entropy condition. 

We are now ready to state DiPerna's uniqueness theorem (cf. Theorem 4.2 in \cite{DiPerna}) for such measure weak solutions.

\begin{thm}
\label{DiPerna_thm}
Suppose $w(t, u, d\lambda)du$ is a measure weak entropy solution of 
$$\partial_t w + \upsilon\cdot \nabla Q(w)=0.$$
Here, $Q\in C^1$, $\upsilon\in \R^n$, and initial condition $w(0,u,d\lambda) = \delta_{w_0(u)}$, where $w_0$ is bounded and integrable.  Suppose also that the following three conditions are satisfied:
\begin{enumerate}
\item Bounded support and probability measure:  The support of $w(t,u,d\lambda)$ is bounded in the interval $A=[a,b]$, for some $a,b\in \R$, uniformly in $(t,u)\in [0,T]\times \R^n$.  Also, for each $(t,u)$, $w(t,u,d\lambda)$ is a probability measure.
\item Initial condition: 
$$\liminf\limits_{t \to 0} \frac{1}{t}\int_0^t \int_{\R^n} \int_A |\lambda- w_0(u)|w(s,u, d\lambda) du ds =0$$
\item Ess sup mass condition:
$${\rm ess\,sup}_t \int_{\R^n} \int_A |\lambda|w(t,u, d\lambda) du < \infty. $$
\end{enumerate}
Then, $w(t,u, d\lambda)$ is the Dirac measure supported on the unique, bounded entropy solution $w(t,u)$ of $\partial_t w+ \upsilon\cdot \nabla Q(w)=0$, with initial condition $w_0(u)$, that is, $w(t,u, d\lambda)= \delta_{w(t,u)}(d\lambda)$.
\end{thm}

\medskip
Given this preamble, we now begin the main part of the proof of Theorem \ref{alpha>1}.
\vskip .2cm

{\it Step 1.} First, we claim that the measures $\{Q^{N,l}\}$ are tight.  This follows the same proof as given in Lemma 1.2, Chapter 8 in \cite{Kipnis}.    Next, as $N\uparrow\infty$ subsequentially, we may obtain a weak limit $Q^l$, and as $l\uparrow\infty$ subsequentially, we obtain a limit point $Q^*$.  We claim that $Q^*$ is supported on measures in the form $\pi(s,du,d\lambda)=\rho(s, u, d\lambda)du$, which are absolutely continuous in $u$.  This also follows the same proof as given for item 1, p. 201 of \cite{Kipnis}.

Also, by the law of large numbers and our initial conditions, we have $\langle \pi_0^{N,l}, F(u,\lambda)\rangle \rightarrow \int_{\R^n} F(u,\rho_0(u))du$ and so the identification $\rho(0,u, d\lambda) = \delta_{\rho_0(u)}(\lambda)$ a.e. $u$.

In addition, $\rho(s,u,d\lambda)$ is supported in a bounded interval, uniformly in $s,u$:  If $M_0<\infty$, that is $h(m)=0$ for some $m$, then there can be at most $M_0$ particles per site in the process.  In particular, $\eta_s^l(x)\leq M_0$ for all $x,s, l$, and so $0\leq \rho(s,u, d\lambda) \leq M_0$ for all $s,u$, without using `attractiveness'.  On the other hand, if $h(m)>0$ for all $m$, 
by the `basic coupling' proof, using `attractiveness', and that the measures $\{\mu_N\}$ are `stochastically bounded' by $\nu_{\rho^\#}$ where $\rho^\# = \|\rho_0\|_\infty$, as given for item (ii) in the proof of Theorem 1.1 of Chapter 8, p. 201-203 in \cite{Kipnis}, we obtain $\rho(s,u,d\lambda)$ is supported in $[0,\rho^{\#}]$ (cf.  related comments, on the `basic coupling', at the beginning of Section \ref{couplingsection}).

We also assert that $\rho(s,u,d\lambda)$ can be identified as probability measures, $\int_0^\infty \rho(s,u,d\lambda)=1$ for $s, u$.   Indeed, $\pi_s^{N,l}(B,[0,\infty))$ is nonrandom, and converges as $N$ and $l\uparrow\infty$ to the Lebesgue measure $m(B)$.  The assertion follows from, say \eqref{cutoff_stat}, and the limit,
\begin{align*}
&\int_0^T f(s) \pi_s^{N,l}(B, [0,\infty))ds= \int_0^T f(s) ds\cdot \frac{1}{N^n}\sum_{x\in \Z^n} 1_B\big(\frac{x}{N}\big) \\
&\ \ \ \ \ \ \ \ \ \rightarrow \int_0^T f(s)\int_B \int_0^\infty \rho(s,u,d\lambda)duds = m(B)\int_0^Tf(s)ds.
\end{align*}

\vskip .2cm

{\it Step 2.}  We will show $Q^*$ a.s. that the density $\rho(t,u,d\lambda)$ satisfies the following four conditions:

\begin{thm}\label{msw_thm}
 $\rho$ is a measure weak solution of $\partial_t \rho + \gamma_\alpha \partial_{\textbf{1}(n)} \Phi\Psi(\rho)=0$. 
 \end{thm}

\begin{thm}
\label{ent_cond_thm}  
The entropy condition holds measure weakly for any $c \in \R$: $$\partial_t |\rho-c|+\gamma_\alpha \partial_{\textbf{1}(n)} [\text{sgn}(\rho-c)(\Phi\Psi(\rho)-\Phi\Psi(c))] \leq 0.$$ 
\end{thm}

\begin{thm}
\label{mass_thm} We have that
$${\rm ess\,sup}_t \int_{\R^n} \int_0^\infty |\lambda - \rho^*|\rho(t,u, d\lambda) du \leq \int_{\R^n}|\rho_0(u)-\rho^*|du< \infty. $$
\end{thm}

\begin{thm}
\label{init_cond_thm}  
The initial condition holds,
 $$\liminf\limits_{t \to 0} \frac{1}{t}\int_0^t \int_{\R^n} \int_0^\infty |\lambda- \rho_0(u)|\rho(s,u, d\lambda) du ds =0.$$
\end{thm}

We prove Theorems \ref{msw_thm}, \ref{ent_cond_thm}, \ref{mass_thm}, \ref{init_cond_thm},  in Subsections \ref{msw_sect}, \ref{ent_cond_section}, \ref{L^1mass_sect}, and \ref{init_sect}, respectively

\vskip .2cm
{\it Step 3.}  Although our initial condition, as $\rho_0(u) = \rho^*$ for $|u|$ large, is not integrable, the function $\hat\rho_0(u) = \rho_0(u)-\rho^*$, is also bounded, and belongs to $L^{1}(\R^n)$.  By considering $\rho^*$-shifted solutions, we will see that the items in Steps 1 and 2 allow to use DiPerna's Theorem \ref{DiPerna_thm} to characterize the limit points $Q^*$.   First, we note the following equivalences.

\vskip .1cm
{\it Equivalence of weak entropy solutions.}  Define $\hat{\rho}(t,u)= \rho(t,u)- \rho^*$ and $\Phi\Psi_{\rho^*}(x)= \Phi\Psi(x+\rho^*)$.  Note as $\Phi\Psi\in C^1$ that also $\Phi\Psi_{\rho^*}\in C^1$ on its domain.  We observe that $\rho(t,u)$ is a weak entropy solution of $\partial_t \rho + \gamma_\alpha \partial_{\textbf{1}(n)} [\Phi\Psi(\rho)]=0$ 
if and only if $\hat{\rho}(t,u)$ is a weak entropy solution of $\partial_t \hat{\rho} + \gamma_\alpha \partial_{\textbf{1}(n)} [\Phi\Psi_{\rho^*}(\hat{\rho})]=0$.

\vskip .1cm
{\it Equivalence of measure weak entropy solutions.}
Similarly, define $\hat{\rho}(t,u, d\lambda)$ through $\hat{\rho}(t,u, F) = \rho(t,u, F+\rho^*)$ for any measurable set $F$. 
Observe, for a function $f$, that
$$\int_a^b f(\lambda) \rho(t,u, d\lambda)=\int_{a-\rho^*}^{b-\rho^*} f_{\rho^*}(\lambda) \hat{\rho}(t,u, d\lambda)$$ where $f_{\rho^*}(\lambda) = f(\lambda+ \rho^*)$.  Note that $\rho(t, u, d\lambda)$ is a probability measure exactly when $\hat{\rho}(t, u, d\lambda)$ is a probability measure.  Also, $\rho(t,u,d\lambda)$ has bounded support in $\lambda$ exactly when $\hat\rho(t,u,d\lambda)$ has bounded support in $\lambda$.

 Hence, 
 $\rho(t,u, d\lambda)$ is a measure weak entropy solution of $\partial_t\rho + \gamma_\alpha \partial_{\textbf{1}(n)} [\Phi\Psi(\rho)]=0$, satisfying the initial condition $\rho_0(u)$ 
if and only if 
$\hat{\rho}(t,u, d\lambda)$ is a measure weak entropy solution of $\partial_t\hat{\rho} + \gamma_\alpha \partial_{\textbf{1}(n)} [\Phi\Psi_{\rho^*}(\hat{\rho})]=0$, satisfying the initial condition $\hat{\rho_0}(u)= \rho_0(u)-\rho^*$, that is, if the following holds:
\begin{enumerate}
\item Measure weakly $\partial_t \hat{\rho} + \gamma_\alpha \partial_{\textbf{1}(n)} \Phi\Psi_{\rho^*}(\hat{\rho})=0$. 
\item Entropy condition holds measure weakly, for $c \in \R$,
$$\partial_t |\hat{\rho}-c|+\gamma_\alpha \partial_{\textbf{1}(n)} [\text{sgn}(\hat{\rho}-c)(\Phi\Psi_{\rho^*}(\hat{\rho})-\Phi\Psi_{\rho^*}(c))] \leq 0.$$ 
\item $L^1$ mass bound holds,
$${\rm ess\,sup}_{t\geq 0} \int_{\R^n}\int_{-\rho^*}^\infty |\lambda|\hat\rho(t,u,d\lambda)du<\infty.$$
\item Initial boundary layer holds, 
$$\lim\limits_{t \downarrow 0} \frac{1}{t}\int_0^t \int_{\R^n} \int_{-\rho^*}^\infty |\lambda- \hat{\rho_0}(u)|\hat{\rho}(s,u, d\lambda) du ds =0.$$
\end{enumerate}

\vskip .2cm

{\it Step 4.}
If now $\rho(t,u, d\lambda)$ satisfies Theorems \ref{msw_thm}, \ref{ent_cond_thm}, \ref{mass_thm}, and \ref{init_cond_thm} 
in Step 2, then $\hat{\rho}(t,u,d\lambda)$ will satisfy the equivalent item versions (1),(2),(3), and (4) in Step 3.

Note also $\hat\rho$ is supported in the bounded interval $A=[-\rho^*, \|\rho_0\|_\infty -\rho^*]$ uniformly in $(t,u)$, as in Step 1 we showed $\rho$ is supported in $[0,\|\rho_0\|_\infty]$.

Then, by Theorem \ref{DiPerna_thm}, we conclude  $\hat{\rho}(t,u,d\lambda)=\delta_{\hat{\rho}(t,u)}(d\lambda)$, where $\hat\rho(t,u)$ is the unique bounded weak entropy solution $\hat{\rho}(t,u)$ of $\partial_t \hat{\rho} + \gamma_\alpha \partial_{\textbf{1}(n)} \Phi\Psi_{\rho^*}(\hat{\rho})=0$ with initial condition $\hat{\rho_0}(u)$.  
It then follows that $\rho(t,u, d\lambda)=\delta_{\rho(t,u)}$, where
$\rho(t,u)$ is the unique bounded entropy solution of $\partial_t \rho + \gamma_\alpha \partial_{\textbf{1}(n)} \Phi\Psi(\rho)=0$ with initial condition $\rho_0(u)$.

Hence, all limit points $Q^*$ of $\{Q^{N,l}\}$ are the same, uniquely characterized in terms of the weak entropy solution of the hydrodynamic equation, $Q^* = \delta_{ \rho(t,u)}$. 
\vskip .2cm

{\it Step 5.} 
We now relate the limit points $\{Q^*\}$ to the limit points of $\{P^N\}$, and thereby prove Theorem \ref{alpha>1}.
We have shown, for test functions $f(s)G(u)$ that
$$Q^*\left(\left|\int_0^t \int_{\R^n} \int_0^\infty f(s)G(u)\lambda \rho(s, u, d\lambda) du  ds- \int_0^t \int_{\R^n} f(s)G(u){\rho}(s,u)ds \right| > \epsilon_0 \right) = 0$$
for all $\epsilon_0 >0$.  Then, as $Q^{N,l}$ on a subequence converges to $Q^*$, by \eqref{cutoff_stat}, 
$$\limsup_{l \to \infty} \limsup_{N \to \infty} Q^{N,l}\left(\left|\int_0^t f(s)\lt \pi^{N,l}_s,  G(u)\lambda \rt ds- \int_0^t \int_{\R^n} f(s)G(u){\rho}(s,u)ds \right| > \epsilon_0\right) = 0,$$
or in other words,
\begin{eqnarray*}
&&\limsup_{l \to \infty} \limsup_{N \to \infty} P_N\Big(\Big|\int_0^t \frac{1}{N^n} \sum_{x \in \Z^n} f(s)G\left(\frac{x}{N}\right)\eta^l_s(x) ds\\
&&\ \ \ \ \ \ \ \ \ - \int_0^t \int_{\R^n} f(s)G(u){\rho}(s,u)ds \Big| > \epsilon_0\Big) = 0.
\end{eqnarray*}

By discrete integration-by-parts, smoothness and compact support of $G$, we may replace $\eta^l_s(x)$ by $\eta_s(x)$ with expected error of order $E^N\int_0^tN^{-(n+1)}\sum_{|x|\leq R'N} \eta_s(x)ds$, which vanishes by say Lemma \ref{particlebound}.

Therefore, 
$$\limsup_{N \to \infty} P^N\left(\left|\int_0^t f(s)\lt \pi^N_s, G\rt ds- \int_0^t \int_{\R^n} f(s)G(u){\rho}(s,u)ds \right| > \epsilon_0\right) = 0.$$
Now, by the assumption FEM, limit points of $\{P^N\}$ are supported on absolutely continuous measures $\pi_s = \bar\rho(s,u)du$; this observation, made in Step 6 in Section \ref{proofsection_alpha<1} for the case $\alpha<1$, also directly applies when $\alpha>1$.  Then, as $\pi\mapsto \int_0^t f(s)\langle \pi^N_s, G\rangle ds$ is continuous, for every limit point $P^*$, we have
$$P^*\left(\left|\int_0^t \int_{\R^n} f(s)G(u) \bar\rho(s,u)du ds- \int_0^t \int_{\R^n} f(s)G(u){\rho}(s,u)ds \right| > \epsilon_0\right) = 0.$$

But, as tightness of $\{P^N_G\}$ was shown with respect to the uniform topology (Proposition \ref{PNtightness}), the limit $\int_{\R^n} G(u)\bar\rho(s,u) du$ is continuous function in time $s$.  One also has that $\int_{\R^n}G(u)\rho(s,u)du$ is continuous in $s$ (cf. Theorem 2.1 \cite{cances}). Therefore, $\int_{\R^n} G(u)\bar\rho(s,u) du= \int_{\R^n} G(u){\rho}(s,u) du$ for all times $s$.

We conclude all limit points $P^*$ are the same, that is, supported on absolutely continuous measures $\pi_t=\rho(t, u)du$ whose density is the unique weak entropy solution of the hydrodynamic equation, and so Theorem \ref{alpha>1} follows. 
\qed

\subsection{Proof of \eqref{cutoff_stat}}
\label{cutoff_sect}
We first note,
for all large $A$ and $\lambda$, by the bound $|f(\lambda)|\leq C\lambda$ and compact support of $G$,
\begin{eqnarray*}
 &&E^{N,l}\int_0^t \lt \pi^{N,l}_s, |G_s(u)||f(\lambda)| {\bf 1}(\lambda > A) \rt ds\\
&&= E^{N}\int_0^t \frac{1}{N^n}\sum_{x \in \Z^n} |G_s\left(\frac{x}{N}\right)||f(\eta^l_s(x))| {\bf 1}(\eta^l_s(x) > A) ds\\
&&\leq C_{G,f} E^{N} \int_0^t \frac{1}{N^n}\sum_{|x|\leq RN} \eta^l_s(x){\bf 1}(\eta^l_s(x) > A) ds
 =C_{G,f} E^{N,l}\int_0^t \langle \pi_s^{N,l}, \lambda {\bf 1}(\lambda >A)\rangle ds.
\end{eqnarray*}
Then, by Lemma \ref{particlebound2} and that $\pi\mapsto \int_0^t \langle \pi_s, \lambda {\bf 1}(\lambda\geq A)\rangle ds$
is a lower semi-continuous function, we have
$$0=\limsup_{A \to \infty} \limsup_{l \to \infty} \limsup_{N \to \infty} E^{N,l}\int_0^t \lt \pi^{N,l}_s, \lambda {\bf 1}(\lambda > A) \rt ds\\
\geq \limsup_{A \to \infty} E_{Q^*}\int_0^t \lt \pi_s, \lambda{\bf 1}(\lambda > A)\rt ds.
$$
In particular, as $\pi\mapsto \int_0^t \lt \pi_s, |G_s(u)||f(\lambda)| {\bf 1}(\lambda \geq A) \rt ds$ is also lower semi-continuous,
\begin{eqnarray}
&&\limsup_{A \to \infty} \limsup_{l \to \infty} \limsup_{N \to \infty} E^{N,l}\int_0^t \lt \pi^{N,l}_s, |G_s(u)||f(\lambda)| {\bf 1}(\lambda > A) \rt ds\nonumber\\
&&\ \ = \ 
\limsup_{A \to \infty} E_{Q^*}\int_0^t \lt \pi_s, |G_s(u)||f(\lambda)|{\bf 1}(\lambda> A)\rt ds
\ =\ 0. 
\label{cutoff1}
\end{eqnarray}

We now argue the left to right equivalence.  In the left-side of \eqref{cutoff_stat}, by \eqref{cutoff1},
we may introduce the indicator function ${\bf 1}(\lambda  \leq A)$.
 Then, as
$\pi\mapsto \int_0^t \lt \pi_s, G_s(u)f(\lambda){\bf 1}(\lambda \leq A)\rt ds$ is continuous, 
we have
\begin{eqnarray*}
&&\limsup_{A \to \infty} \lim_{l' \to \infty} \lim_{N' \to \infty} Q^{N',l'}\Big(\big|\int_0^t \lt \pi^{N',l'}_s, G_s(u)f(\lambda){\bf 1}(\lambda \leq A)\rt\big| > \epsilon_0\Big)\\
&& = 
\limsup_{A \to \infty} Q^*\Big(\big|\int_0^t \lt \pi_s, G_s(u)f(\lambda){\bf 1}(\lambda\leq A)\rt ds\big| > \epsilon_0\Big)=0.
\end{eqnarray*}
The right-side of \eqref{cutoff_stat} follows now by \eqref{cutoff1} applied again.

The right to left equivalence in \eqref{cutoff_stat} follows by similar steps in reverse, given now $Q^*\big(\big|\int_0^t \langle \pi_s, G_s(u)f(\lambda)\rangle ds\big| \geq \epsilon_0\big) = 0$ for all $\epsilon_0$.  Here, without loss of generality we have replaced `$>$' by `$\geq$' to maintain the correct bounds implied by weak convergence.
\qed

\section{Measure weak solutions:  Proof of Theorem \ref{msw_thm}}
\label{msw_sect}

The argument follows some of the initial reasoning given for the proof of Theorem \ref{alpha<1}, in the $\alpha <1$ case, relying however on the $1$-block estimate Lemma \ref{1block_prop_alpha>1}.

{\it Step 1.}  The same estimate as in Step 1 in Section \ref{proofsection_alpha<1}, with respect to the martingale $M^{N,G}_t$, gives that
$$\limsup\limits_{N \to \infty} P^N\Big(\big|\lt \pi^N_0 , G_0\rt + \int_0^t  \lt \pi^N_s, \partial_sG_s \rt ds + \int_0^t N L_N\lt \pi^N_s , G_s\rt ds \big|> \epsilon_0 \Big)=0.$$
 Here, we recall from \eqref{gen_quad},
\begin{equation*}
\int_0^t N L_N\lt \pi^N_s , G_s\rt ds = \int_0^t \frac{N}{N^n}\sum\limits_{x \in \Z^n} \sum\limits_{\|d\|=1}^\infty \frac{1}{\|d\|^{n+\alpha}} gh_d(\eta_s(x)) \big[G_s\big(\frac{x+d}{N}\big)-G_s\big(\frac{x}{N}\big)\big]ds.
\end{equation*}

\vskip .2cm

{\it Step 2.}
We would like to replace $G_s\left(\frac{x+d}{N}\right)-G_s\left(\frac{x}{N}\right)$ by $\nabla G_s(\frac{x}{N}) \cdot d/N$. To this aim, noting $gh_d(\eta_s(x))\leq \kappa\|h\|_\infty\eta_s(x)$, by Lemma \ref{lemma23_cor1}, we may truncate the sum on $d$ to $\|d\|\leq \epsilon N$, in terms of a parameter $\epsilon$ which will vanish after $N$ diverges. 
Next, as $|G_s\left(\frac{x+d}{N}\right)-G_s\left(\frac{x}{N}\right)-\nabla G_s\left(\frac{x}{N}\right)\cdot \frac{d}{N} |\leq  \|\nabla^2 G_s\|
{\bf 1}(|x|<(R+\epsilon)N)\|\frac{d}{N}\|^2$, we have
\begin{eqnarray*}
&&\int_0^t \frac{N}{N^n} \sum_{x \in \Z^n}\sum_{\|d\|=1}^{\epsilon N} \frac{1}{\|d\|^{n+\alpha}} gh_d(\eta_s(x)) \Big|G_s\left(\frac{x+d}{N}\right)-G_s\left(\frac{x}{N}\right)- \nabla G_s\left(\frac{x}{N}\right)\cdot \frac{d}{N}\Big|ds\\
&& \leq 
\frac{\kappa\|h\|_\infty\|\nabla^2 G\|}{N} \sum_{\|d\|=1}^{\epsilon N} \frac{\|d\|^2}{\|d\|^{n+\alpha}} \frac{1}{N^n} \sum_{|x|\leq (R+\epsilon)N} \eta_s(x) \ \leq \ \frac{\epsilon C_{G,\alpha}}{N^n}\sum_{|x| \leq (R+\epsilon)N } \eta_s(x),
\end{eqnarray*}
which vanishes in expected value, noting Lemma \ref{particlebound}, as $N\uparrow\infty$ and $\epsilon\downarrow 0$.
Therefore, 
\begin{eqnarray}
&&\lt \pi^N_0 , G_0\rt + \int_0^t  \lt \pi^N_s, \partial_sG_s \rt ds \nonumber\\
&&\ \ \ + \int_0^t \frac{1}{N^n}\sum\limits_{x \in \Z^n} \sum\limits_{\|d\|=1}^{\epsilon N} \frac{1}{\|d\|^{n-1+\alpha}} gh_d(\eta_s(x)) \nabla G_s\left(\frac{x}{N}\right) \cdot \frac{d}{\|d\|} ds
\label{alpha>1_step2}
\end{eqnarray}
converges to zero in probability after taking the appropriate limits. 
Moreover, with similar reasoning,  
we may further replace the sum on $d$ to a truncated sum over $\|d\|\leq D$, where $D$ will diverge after $N$.

\vskip .2cm

 {\flushleft \it Step 3a.} 
Now, by the method of Step 3a in Section \ref{proofsection_alpha<1} for the $\alpha<1$ case, 
we substitute $gh_d(\eta_s(x))$ with  $(gh_d)^l(\eta_s(x))$ where $l$ will go to infinity after $N$ but before $D$. We will also replace $\eta_s(x)$ by $\eta^l_s(x)$ in the first and second terms in \eqref{alpha>1_step2}. 
Hence,
\begin{eqnarray*}
&& \frac{1}{N^n} \sum_{x \in \Z^n} \eta^l_0(x)G_0\left(\frac{x}{N}\right) +\int_0^t \frac{1}{N^n} \sum_{x \in \Z^n} \eta^l_s(x)\partial_s G_s\left(\frac{x}{N}\right) ds\\
&&\ \ \ \ \ \ \ \ +\int_0^t \frac{1}{N^n}\sum\limits_{x \in \Z^n} \sum\limits_{\|d\|=1}^{D} \frac{1}{\|d\|^{n-1+\alpha}} (gh_d)^l(\eta_s(x)) \nabla G_s\left(\frac{x}{N}\right)\cdot \frac{d}{\|d\|}ds
\end{eqnarray*}
converges to zero in probability as the appropriate limits are taken. 

\vskip .2cm

{\flushleft \it Step 3b.} We now replace 
$(gh_d)^l(\eta_s(x))$ by $\Phi(\eta^l_s(x))\Psi(\eta^l_s(x))$ for $1\leq \|d\|\leq D$, using the $1$-block estimate Proposition \ref{1block_prop_alpha>1} and $\sum_{\|d\|=1}^D \|d\|^{-(n-1+\alpha)}<\infty$.   
After this replacement, we will have shown 
\begin{eqnarray*}
&& \frac{1}{N^n} \sum_{x \in \Z^n} \eta^l_0(x)G_0\big(\frac{x}{N}\big) +\int_0^t \frac{1}{N^n} \sum_{x \in \Z^n} \eta^l_s(x)\partial_s G_s\big(\frac{x}{N}\big) ds\\
&&\ \ \ \ \ \ \ \ \  +\int_0^t \frac{1}{N^n}\sum\limits_{x \in \Z^n} \Phi\Psi(\eta^l_s(x)) \sum\limits_{\|d\|=1}^{D} \frac{d}{\|d\|^{n+\alpha}} \cdot \nabla G_s\big(\frac{x}{N}\big) ds
\end{eqnarray*}
vanishes in probability as $N,l\uparrow\infty$. 

By $\Phi\Psi(\eta(x))\leq \kappa\|h\|_\infty \eta(x)$, compact support of $G$, and Lemma \ref{particlebound}, we can further replace $\sum_{\|d\|=1}^D d/\|d\|^{n+\alpha}$  with $\gamma_\alpha \textbf{1}(n) = \sum_{d=1}^\infty d/\|d\|^{n+\alpha}$, taking $D\uparrow\infty$ at the end, where $\textbf{1}(n)$ is the unit vector in the direction $\lt 1, 1, \dots, 1 \rt$.

\vskip .2cm

{\flushleft \it Step 3c.} Now, in terms of the Young measures $\pi^{N,l}_t$ defined in Section \ref{proof_section_alpha>1}, we have
 $$\lt \pi^{N, l}_0 , G_0(u)\lambda\rt + \int_0^t  \lt \pi^{N, l}_s, \partial_s G_s(u)\lambda \rt ds -\gamma_\alpha \int_0^t \lt \pi^{N, l}_s, \partial_{\textbf{1}(n)} G_s(u)\Phi\Psi(\lambda) \rt ds$$
vanishes in probability.
\vskip .2cm

{\it Step 4.}
Consider a limit point $Q^*$ of the measures $\{Q^{N,l}\}$ governing $\pi^{N,l}_t$. 
Recall that we observed in Step 1 in Section \ref{proof_section_alpha>1} that $Q^*$ is supported on $\pi_s = \pi(s,du,d\lambda) = \rho(s,u,d\lambda)du$, $\rho(0,u,d\lambda) = \delta_{\rho_0(u)}$, and also $\langle \pi^{N,l}_0, G_0(u)\lambda\rangle$ converges to $\int G_0(u)\rho_0(u)du$.

Then, on a subsequence as $N\uparrow\infty$ and $l\uparrow\infty$, we conclude, noting \eqref{cutoff_stat} and $\Phi\Psi(\lambda)\leq C\lambda$, that
a.s. $Q^*$,
\begin{eqnarray*}
&&\int_{\R^n} G_0(u) \int_0^\infty \lambda \rho(0,u, d\lambda) du + \int_0^\infty \int_{\R^n} \partial_s G_s(u) \int_0^\infty \lambda \rho(s,u, d\lambda)  du ds\\
&&\ \ \ \ \ \ \ \ \ \ \ \ \ \ \ 
+ \gamma_\alpha \int_0^\infty \int_{\R^n} \int_0^\infty \partial_{\textbf{1}(n)}G_s(u)\Phi\Psi(\lambda)\rho(s,u, d\lambda) duds=0.
\end{eqnarray*}
Here, we replaced the limit $t$ with $\infty$, noting that $G$ has compact support in $[0,T)\times \R^n$.  Hence, $\rho(s,u,d\lambda)$ is a measure-weak solution of the hydrodynamic equation. \qed

\section{A coupled process} \label{couplingsection}

We introduce the basic coupling for misanthrope processes.
Let $\tilde{P}^N$ denote the distribution of the coupled process $(\eta_t, \xi_t)$ with
generator $\tilde{L}$, given by its action on test functions,
\begin{eqnarray*}
\tilde{L}f(\eta, \xi) & = & \sum\limits_{x,y}p(y-x)\min\nolimits_{x, y}(f(\eta^{x, y}, \xi^{x, y}) -f(\eta, \xi)) \nonumber\\
& + & \sum\limits_{x,y}p(y-x)(b(\eta(x), \eta(y)) -\min\nolimits_{x,y})(f(\eta^{x, y}, \xi) -f(\eta, \xi)) \nonumber\\
& + & \sum\limits_{x,y}p(y-x)(b(\xi(x), \xi(y)) -\min\nolimits_{x, y})(f(\eta, \xi^{x, y}) -f(\eta, \xi)), \nonumber
\end{eqnarray*}
where $\min_{x,y}=\min\{b(\eta(x), \eta(y)), b(\xi(x), \xi(y))\}$.
From the form of the generator, it follows that the marginals are themselves misanthrope processes.  

Suppose now that the process is `attractive', that is when $b(n,m)=g(n)h(m)$, with $g$ increasing and $h$ decreasing in particle numbers. 
Then, if $\eta_s(x) \leq \xi_s(x)$ for all $x\in \Z^n$,
at any later time $t\geq s$, we still have the same ordering.  This observation is the crux of the proof of the `$L^\infty$' bound in \cite{Kipnis}, referred to in Step 1 in Section \ref{proof_section_alpha>1}.
This is the first of the two ways where `attractiveness' is used in the proof of Theorem \ref{alpha>1}.

We will use the following teminology. For any set $\Lambda \subseteq \Z^n$, we write $\eta \geq \xi$ on $\Lambda$ if $\eta(x) \geq \xi(x)$ for all $x \in \Lambda$, and we write $\eta > \xi$ on $\Lambda$ if $\eta \geq \xi$ on $\Lambda$ and $\eta(x) > \xi(x)$ for at least one $x \in \Lambda$. If $\eta \geq \xi$ or $\xi \geq \eta$ on $\Lambda$, we say that $\eta$ and $\xi$ are ordered on $\Lambda$. Otherwise, we say that $\eta$ and $\xi$ are unordered on $\Lambda$.

Let $U_{\Lambda}(\eta, \xi)={\bf 1}(\eta {\rm \ and \ }\xi \ {\rm are \ not \ ordered \ on \ }\Lambda)$.
Let $U_{x,d}(\eta, \xi) = U_{\{x,x+d\} }(\eta, \xi)$. We also define  
$$O_{x,d}(\eta, \xi) =\left\{ \begin{array}{cr} 1 & \text{if } \eta > \xi \text{ on } \{x, x+d\}\\
-1 & \text{if } \xi > \eta \text{ on } \{x, x+d\}\\
 0 & \text{otherwise} \\
\end{array} \right.$$

Define the coupled empirical measure by 
$$\tilde{\pi}^N_t = \frac{1}{N^n}\sum_{x\in \Z^n} |\eta_t(x)-\xi_t(x)| \delta_{{x}/{N}}.$$

We now introduce martingales which will be useful in the sequel.  The first two are the coupled versions of $M^{N,G}_t$ and the associated `variance' martingale:  For test functions $G$ on the coupled space, define the martingale,
$$\tilde{M}^{N,G}_t = \lt \tilde{\pi}^N_t, G_t \rt - \lt \tilde{\pi}^N_0, G_0 \rt - \int_0^t  \lt \tilde{\pi}^N_s, \partial_sG_s \rt +N\tilde{L}_N\lt \tilde{\pi}^N_s, G_s \rt ds.$$
With respect to the quadratic variation,
 $$\lt \tilde{M}^{N,G}\rt_t = \int_0^t N\tilde{L}_N[(\lt \tilde{\pi}^N_s, G_s \rt)^2]-2N \lt \tilde{\pi}^N_s, G_s \rt \tilde{L}_N \lt \tilde{\pi}^N_s, G_s \rt ds,$$
 the process $(\tilde{M}^{N,G}_t)^2 - \lt \tilde{M}^{N,G}\rt_t$ is also a martingale.
 We may compute
\begin{eqnarray}
\label{coupled_gen}
&&N\tilde{L}_N\lt \tilde{\pi}^N_s, G_s \rt \\
&&= \frac{N}{N^n} \sum\limits_{x \in \Z^n} \sum\limits_{\|d\|=1}^\infty \frac{1}{\|d\|^{n+\alpha}} (gh_d(\eta_s(x))-gh_d(\xi_s(x))O_{x,d}(\eta_s, \xi_s) \big[G_s\big(\frac{x+d}{N}\big)-G_s\big(\frac{x}{N}\big) \big]\nonumber\\
&&\ \ \  -\frac{N}{N^n}\sum\limits_{x \in \Z^n} \sum\limits_{\|d\|=1}^\infty \frac{1}{\|d\|^{n+\alpha}} (gh_d(\eta_s(x))-gh_d(\xi_s(x)) U^{\pm}_{x,d}(\eta_s, \xi_s) \big[G_s\big(\frac{x+d}{N}\big)+G_s\big(\frac{x}{N}\big) \big],\nonumber
\end{eqnarray}
where 
$$U^{\pm}_{x,d}(\eta_s, \xi_s) = \left\{\begin{array}{rl}
1& \ {\rm if \ } \eta_s(x)>\xi_s(x) \ {\rm and \ } \eta_s(x+d)<\xi_s(x+d)\\
-1& \ {\rm  if\ } \eta_s(x)<\xi_s(x) \ {\rm and \ } \eta_s(x+d)>\xi_s(x+d)\\
0& \ {\rm otherwise.}
\end{array}\right.
$$ 
Note that $U_{x,d}(\eta_s, \xi_s)= |U^{\pm}_{x,d}(\eta_s, \xi_s)|$.

When the process is `attractive', we have
$$(gh_d(\eta_s(x))-gh_d(\xi_s(x)) U^{\pm}_{x,d}(\eta_s, \xi_s)=|gh_d(\eta_s(x))-gh_d(\xi_s(x)| U_{x,d}(\eta_s, \xi_s).$$
In this case, the second line of the generator computation \eqref{coupled_gen} simplifies to
\begin{equation}
\label{alpha>1_generator}
  -\frac{N}{N^n}\sum\limits_{x \in \Z^n} \sum\limits_{\|d\|=1}^\infty \frac{1}{\|d\|^{n+\alpha}} |gh_d(\eta_s(x))-gh_d(\xi_s(x)| U_{x,d}(\eta_s, \xi_s) \big[G_s\big(\frac{x+d}{N}\big)+G_s\big(\frac{x}{N}\big) \big].
\end{equation}
We remark that this is the second of two places where the `attractiveness' condition is explicitly used, featuring in the proof of the `Ordering Lemma', stated later.

\blem \label{alpha>1martbound} 
When $\alpha>1$ and $G$ is nonnegative,
\begin{eqnarray*}
&&N\tilde{L}_N\lt \tilde{\pi}^N_s, G_s \rt \\
&& \leq \frac{N}{N^n} \sum\limits_{x \in \Z^n} \sum\limits_{\|d\|=1}^\infty \frac{1}{\|d\|^{n+\alpha}} (gh_d(\eta_s(x))-gh_d(\xi_s(x))O_{x,d}(\eta_s, \xi_s) \big[G_s\big(\frac{x+d}{N}\big)-G_s\big(\frac{x}{N}\big)\big].
\end{eqnarray*} \elem

\bpf The bound follows as, in \eqref{alpha>1_generator}, all terms in the second line are nonnegative.
\epf

In the next two results, we will start the coupled process $(\eta_s,\xi_s)$ from an arbitrary initial distribution $\tilde \mu_N$ whose marginals are $\mu_N$ and $\nu_c$, for a $0\leq c\leq M_0$ if $M_0<\infty$, and $c\geq 0$ if otherwise.  The coupled process measure is denoted by $\tilde P_N$ and the associated expectation is given by $\tilde E_N$.

For the quadratic variation, $\lt\tilde{M}^{N,G} \rt_t$, a straightforward computation gives that
\begin{eqnarray*}
&&N\tilde{L}_N[(\lt \tilde{\pi}^N_s, G_s \rt)^2]-2N \lt \tilde{\pi}^N_s, G_s \rt \tilde{L}_N \lt \tilde{\pi}^N_s, G_s \rt \\
&&\ = \frac{N}{N^{2n}}\sum\limits_{x\in \Z^n} \sum\limits_{\|d\|=1}^\infty \frac{1}{\|d\|^{n+\alpha}} (gh_d(\eta_s(x))-\min\{gh_d(\eta_s(x)),gh_d(\xi_s(x))\})\times\\ 
&&\ \ \ \ \ \ \ \ \ \ \ \ \ \ \big[(G_s(\frac{x+d}{N})-G_s(\frac{x}{N}))O_{x,d}(\eta_s, \xi_s)- (G_s(\frac{x+d}{N})+G_s(\frac{x}{N}))U^\pm_{x,d}(\eta_s, \xi_s)\big]^2\\
&& \ \ \ \ + \frac{N}{N^{2n}}\sum\limits_{x\in \Z^n} \sum\limits_{\|d\|=1}^\infty \frac{1}{\|d\|^{n+\alpha}} (gh_d(\xi_s(x))-\min\{gh_d(\eta_s(x)),gh_d(\xi_s(x))\})\times \\
&&\ \ \ \ \ \ \ \ \ \ \ \ \ \ \big[(G_s(\frac{x+d}{N})-G_s(\frac{x}{N}))O_{x,d}(\xi_s, \eta_s)- (G_s(\frac{x+d}{N})+G_s(\frac{x}{N}))U^\pm_{x,d}(\eta_s, \xi_s)\big]^2.
\end{eqnarray*}

\blem \label{alpha>1quadbound}
When $\alpha>1$, we have
\begin{eqnarray*}
&&\tilde E^N|\lt \tilde{M}^{N,G}\rt_t| \leq \frac{C(G,t)}{N^n} \\
&& + \frac{16\|G\|^2N}{N^n} \tilde E^N\int_0^t \frac{1}{N^n}\sum\limits_{\substack{|x| \leq RN\\ |x+d| \leq RN}}\sum\limits_{\|d\|=1}^\infty \frac{1}{\|d\|^{n+\alpha}} 
|gh_d(\eta_s(x))-gh_d(\xi_s(x)| U_{x,d}(\eta_s, \xi_s) ds.
\end{eqnarray*}\elem

\bpf
In the expression for the quadratic variation, we may bound factors 
$(gh_d(\eta_s(x))-\min\{gh_d(\eta_s(x)),gh_d(\xi_s(x))\})$ by $|gh_d(\eta_s(x))-gh_d(\xi_s(x)|$.  Also, we note $|O_{x,d}|\leq 1$ and $(U^\pm_{x,d})^2 = U_{x,d}$.  Using the inequality
 $(a-b)^2\leq 2(a^2+b^2)$, 
 we have that $N\tilde{L}_N[(\lt \tilde{\pi}^N_s, G_s \rt)^2]-2N \lt \tilde{\pi}^N_s, G_s \rt \tilde{L}_N \lt \tilde{\pi}^N_s, G_s \rt$
is bounded above by
\begin{eqnarray*}
\label{coupled_quad}
&&\int_0^t 4\frac{N}{N^{2n}}\sum\limits_{x\in \Z^n} \sum\limits_{\|d\|=1}^\infty  \frac{1}{\|d\|^{n+\alpha}}  |gh_d(\eta_s(x))-gh_d(\xi_s(x)|\times\\ 
&&\ \ \ \ \ \ \ \ \ \ \ \ \  \big[\big(G_s\big(\frac{x+d}{N}\big)-G_s\big(\frac{x}{N}\big)\big)^2+ \big(G_s\big(\frac{x+d}{N}\big)+G_s\big(\frac{x}{N}\big)\big)^2U_{x,d}(\eta_s, \xi_s)\big]ds,\nonumber
\end{eqnarray*}
which we split as $A_1 + A_2$, the term $A_1$ involving $\big(G_s\big(\frac{x+d}{N}\big)-G_s\big(\frac{x}{N}\big)\big)^2$ and $A_2$ involving the other squared quantity.

Since, $gh_d(\eta(x)) \leq \kappa \|h\|_\infty\eta(x)$ by \eqref{LB},
we observe that $|gh_d(\eta_s(x)) - gh_d(\xi_s(x))|\leq 2\|h\|_\infty\kappa(\eta_s(x) + \xi_s(x))$.  Hence, $A_1 \leq A_{11} + A_{12}$, where $A_{11}$ and $A_{12}$ involve each only the $\eta_\cdot$ and $\xi_\cdot$ process respectively.   By the proof of Lemma \ref{alpha<1quadbound}, starting from \eqref{lemma4.6_eq}, $\tilde E_N A_{11}\leq K_G t/N^n$.  
A similar bound and argument holds when $\xi_s(x)$ is present as $\nu_c$ is invariant, and therefore $\xi_s\sim \nu_c$ and $\tilde E_N \sum_{a\leq |x|\leq b}\xi_s(x) =c(b^n - a^n)$.
Hence, $\tilde E_N A_1\leq C(G,t)/N^n$.

The remaining part $\tilde E_N A_2$,
as the the sum of the $G$'s squared is bounded by $4\|G\|^2\big[{\bf 1}(|x|\leq RN)+ {\bf 1}(|x+d|\leq RN)\big]$,
is majorized by
$$16\|G\|^2\frac{N}{N^n}\tilde E^N\int_0^t\frac{1}{N^{n}}\sum\limits_{\substack{|x|\leq RN\\ |x+d|\leq RN} } \sum\limits_{\|d\|=1}^\infty  \frac{1}{\|d\|^{n+\alpha}} |gh_d(\eta_s(x))-gh_d(\xi_s(x)|U_{x,d}(\eta_s, \xi_s)ds.
$$
This finishes the proof.
\epf

We now state an `Ordering Lemma' which,
in essence, tells us that $\eta_t$ and $\xi_t$ are ordered on average, even if they are not initially ordered.
This result is analogous to those in the finite-range setting, Lemma 3.3 in \cite{Rezakhanlou} and Lemma 2.2 on p. 209 of \cite{Kipnis}.

\begin{lem}[Ordering Lemma]
\label{orderinglemma}
For $\alpha >1$ and $a,b\in \R^n$,
$$\limsup_{N \to \infty} \tilde{E}^N\int_0^t \frac{1}{N^n} \sum\limits_{\substack{x \in [a, b]N\\ x+d \in [a,b]N}} \sum\limits_{\|d\|=1}^\infty \frac{1}{\|d\|^{n+\alpha}}
|gh_d(\eta_s(x))-gh_d(\xi_s(x)| U_{x,d}(\eta_s, \xi_s) ds= 0$$ 
where $[a,b]N=[aN, bN]$ and $[a,b]=\prod_{j=1}^n [a_j,b_j]$ denotes the $n$-dimensional hyper-rectangle with diagonal extending from $a$ to $b$. 

We also have, for all $d$ with $\|d\|\geq 1$, that
$$\limsup_{N \to \infty} \tilde{E}^N\int_0^t \frac{1}{N^n} \sum\limits_{|x| \leq RN} U_{x,d}(\eta_s, \xi_s) ds= 0.$$

\end{lem}

We postpone the proof the `Ordering Lemma' to the Appendix.

\section{Entropy condition:  Proof of Theorem \ref{ent_cond_thm}}
\label{ent_cond_section}

 We note, as specified in the definition of the measure weak entropy condition, the test functions $G$ in this section are nonnegative.

\vskip .2cm

{\it Step 1.} Since $\rho \geq  0$ a.e. (cf.  Step 1 of Section \ref{proof_section_alpha>1}),  it is enough to prove Theorem \ref{ent_cond_thm} when $c \geq 0$.   When the max occupation number $M_0<\infty$, it is enough to consider $0\leq c\leq M_0$.

Suppose we may show, for $\epsilon_0>0$ and $t\leq T$, that
\begin{eqnarray}
\label{microentropy}
&&\liminf\limits_{l \to \infty} \liminf\limits_{N \to \infty}P^N\Big(\frac{1}{N^n}\sum_{x\in \Z^n} \big|\rho_0(x/N)- c\big|G_0\big(\frac{x}{N}\big) \\
&&\ \ \ \ \ \ \ + \int_0^t \frac{1}{N^n} \sum\limits_{x \in \Z^n}  \big|\eta_s^l(x)-c\big|\partial_sG_s\big( \frac{x}{N}\big)ds\nonumber\\
&&\ \ \ \ \ \ \ + \gamma_\alpha\int_0^t \frac{1}{N^n} \sum\limits_{x \in \Z^n}  \sgn(\eta_s^l(x)-c)\big(\Phi\Psi(\eta_s^l(x))-\Phi\Psi(c)\big)\partial_{\textbf{1}(n)}G_s\big(\frac{x}{N}\big)ds \geq -\epsilon_0\Big)  =1.\nonumber
\end{eqnarray}
 In terms of Young measures and $Q^{N,l}$, 
\eqref{microentropy} is written
\begin{eqnarray*}
&&\liminf\limits_{l \to \infty} \liminf\limits_{N \to \infty} Q^{N, l} \Big(\frac{1}{N^n}\sum_{x\in \Z^n} \big|\rho_0(x/N)- c\big|G_0\big(\frac{x}{N}\big)
 + \int_0^t \lt \pi^{N, l}_s, \partial_sG_s( u)|\lambda - c|\rt ds \\
&&\ \ \ \ \ \ \ \ \ \ \ \ \ \ \  \ \ \ \ \ \ \ \ \ + \gamma_\alpha\int_0^t \lt \pi^{N,l}_s, \partial_{\textbf{1}(n)}G_s(u)q(\lambda,c)ds \geq -\epsilon_0\Big)  =1,
\end{eqnarray*}
where $q(\lambda, c) = \sgn(\lambda-c)(\Phi\Psi(\lambda)-\Phi\Psi(c))$. 

In this case, the desired measure weak formulation of the entropy condition would follow:  By tightness of $\{Q^{N,l}\}$, let $Q^*$ be a limit point.
Such a $Q^*$ is supported on absolutely continuous measures $\pi_s = \rho(s,u,d\lambda)du$ and $\rho(0,u,d\lambda) = \delta_{\rho_0(u)}$ (cf. Step 1 of Section \ref{proof_section_alpha>1}).  Then, as $\Phi\Psi(\lambda) \leq \kappa \|h\|_\infty\lambda$ (cf. \eqref{LB}),  by the weak convergence statement \eqref{cutoff_stat}, we would have $Q^*$ a.s. that
\begin{eqnarray*}
&&\int_{\R^n} G_0(u) |\rho_0(u)-c|du + \int_0^\infty \int_{\R^n}  \partial_sG_s( u) \int_0^\infty |\lambda - c| \rho(s,u, d\lambda) du ds\\
&&\ \ \ \  + \gamma_\alpha\int_0^\infty \int_{\R^n} \partial_{\textbf{1}(n)} G_s(u) \int_0^\infty q(\lambda,c)\rho(s,u, d\lambda)duds \geq 0.
\end{eqnarray*}

\vskip .2cm

{\it Step 2.} To begin to establish \eqref{microentropy},
consider a coupled process $(\eta_t, \xi_t)$ where the initial distribution is such that $\xi_0$ is the invariant measure $\nu_c$ with density $c$.  We will specify the form of the coupled initial distribution at the beginning of Subsection \ref{coupmicro_sect}, and show there a coupled version of the microscopic entropy inequality:
For $\epsilon >0$, and $t \leq T$, we have
\begin{eqnarray}
&&\liminf\limits_{l \to \infty} \liminf\limits_{N \to \infty}\tilde{P}^N\Big(\frac{1}{N^n}\sum_{x\in \Z^n}\big|\rho_0(x/N) - c\big| 
G_0\big(\frac{x}{N}\big)\nonumber\\
&&+ \int_0^t \frac{1}{N^n} \sum\limits_{x \in \Z^n} \big|\eta_s^l(x)-\xi_s^l(x)\big|\partial_sG_s\big( \frac{x}{N}\big)ds + \gamma_\alpha \int_0^t \frac{1}{N^n} \sum\limits_{x \in \Z^n}   \sgn(\eta_s^l(x)-\xi_s^l(x))\nonumber \\
&&\ \ \ \ \ \ \ \ \times \big(\Phi\Psi(\eta_s^l(x))-\Phi\Psi(\xi_s^l(x))\big)\partial_{\textbf{1}(n)}G_s\big(\frac{x}{N}\big)ds \geq -\epsilon_0\Big)  =1.
\label{coupmicroentro}
\end{eqnarray}

\vskip .2cm

{\it Step 3.}
We now show how the microscopic entropy inequality \eqref{microentropy} can be deduced from the coupled microscopic entropy inequality \eqref{coupmicroentro}.
 It is enough to show that the following terms vanish as $N$ and then $l$ go to infinity:
\begin{eqnarray}
&&\tilde{E}^N \int_0^t \frac{1}{N^n} \sum\limits_{|x|\leq RN } \big|\big|\eta_s^l(x)-\xi_s^l(x)\big|-\big|\eta_s^l(x)-c\big|\big|ds\nonumber\\
&& \ \ {\rm and \ \ \ \ }
\tilde{E}^N\int_0^t \frac{1}{N^n} \sum\limits_{|x|\leq RN } \big|q(\eta_s^l(x),\xi_s^l(x))-q(\eta_s^l(x), c)\big| ds.
\label{terms}\end{eqnarray}

To analyze the second term, we note, as $\Psi$ is bounded by $\|h\|_\infty$ and $\Phi$ is Lipschitz, that $|q(\eta_s^l(x),\xi_s^l(x))-q(\eta_s^l(x), c)|=O\big(\eta^l_s(x)+\xi^l_s(x)+c\big)$. Thus, we can introduce the indicator function ${\bf 1}(\eta^l_s(x)\vee \xi^l_s(x) \leq A)$, the error vanishing by Lemma \ref{particlebound2} and that $\xi_\cdot \sim \nu_c$. Now note that $q(z,w)$ is uniformly continuous on $[-A,A]^2$.

On the other hand, for the first term, by the triangle inequality, $||\eta_s^l(x)-\xi_s^l(x)| - |\eta_s^l(x) -c|| \leq |\xi^l_s(x) -c|$.

   Hence, the terms in \eqref{terms} will vanish, if we show that 
$$\limsup\limits_{l\uparrow\infty}\limsup\limits_{N\rightarrow\infty}
\tilde E^N\int_0^t \frac{1}{N^n} \sum\limits_{|x| \leq RN } |\xi_s^l(x))-c|ds = 0.$$  
But, since the state $\xi_s$ has distribution $\nu_c$,  
it follows that
$$\tilde E^N\int_0^t \frac{1}{N^n} \sum\limits_{|x| \leq RN } |\xi_s^l(x))-c|ds \leq (2R+1)^n t \cdot E_{\nu_c}\big |\xi_0^l(0)-c\big |,$$
which vanishes by the law of large numbers as $l\uparrow\infty$. \qed

\subsection{Proof of \eqref{coupmicroentro}}
\label{coupmicro_sect}

We proceed in some steps, recalling estimates in Section \ref{couplingsection}.
First, we specify the initial coupled distribution in Step 2 above:  We will take $\tilde{\mu}_N$ as a product measure over $x\in \Z^n$ with $x$-marginal given by $\tilde{\mu}_N(\eta_0(x)\geq \xi_0(x))=1$ if $\rho_0(x/N)\geq c$ and $\tilde{\mu}_N(\eta_0(x)\leq \xi_0(x))=1$ if $\rho(x/N)\leq c$.  Such a coupled initial measure may be constructed (cf. \cite{Liggett}) as the $x$-marginals of $\mu^N$ and $\nu_c$ are stochastically ordered, that is $\Theta_{\rho_0(x/N)}$ is stochasically more or less than $\Theta_c$ if $\rho_0(x/N)$ is more or less than $c$ respectively.  Then, $\tilde P_N$ is the coupled process measure starting from $\tilde \mu_N$.
\vskip .2cm

{\it Step 1.}
 By Lemma \ref{alpha>1quadbound} and the Ordering Lemma \ref{orderinglemma}, the expected value $\tilde E^N\langle \tilde M^{N,G}\rangle_t$ vanishes as $N\uparrow\infty$.
Hence, for $\epsilon_0>0$, 
$\lim\limits_{N \to \infty} \tilde{P}^N\big(- \tilde{M}^{N,G}_t  \geq -\epsilon_0 \big)=1$.

Since $G$ has compact support in $[0,T) \times \R^n$, we have
$\lt \tilde{\pi}^N_t, G_t \rt=0$ for $t\geq T$, and so
$-\tilde{M}^{N,G}_t = \lt \tilde{\pi}^N_0, G_0 \rt+ \int_0^t \lt \tilde{\pi}^N_s, \partial_sG_s \rt +N\tilde{L}_N\lt \tilde{\pi}^N_s, G_s \rt ds$. 
It follows, as $G$ is nonnegative, from the bound in Lemma \ref{alpha>1martbound}, that
\begin{eqnarray*}
&&\lim\limits_{N\uparrow\infty}
\tilde{P}^N\Big(\frac{1}{N^n}\sum_{x\in \Z^n}\big|\eta_0(x) - \xi_0(x)\big| G_0\big(\frac{x}{N}\big)+\int_0^t \frac{1}{N^n}\sum\limits_{x \in \Z^n}  |\eta_s(x)-\xi_s(x)|\partial_s G_s\big(\frac{x}{N}\big)ds \\
&&\ \  + \int_0^t \frac{N}{N^n} \sum\limits_{x \in \Z^n} \sum\limits_{\|d\|=1}^\infty \frac{1}{\|d\|^{n+\alpha}} (gh_d(\eta_s(x))-gh_d(\xi_s(x)))O_{x,d}(\eta_s, \xi_s)\\
&&\ \ \ \ \ \ \ \ \ \ \ \ \ \ \ \ \ \ \ \ \ \ \times \big[G_s\big(\frac{x+d}{N}\big)-G_s\big(\frac{x}{N}\big) \big] ds \geq -\epsilon_0\Big) = 1.
\end{eqnarray*}

We now replace the second integral in the last display by one with a nicer form.  We make substitutions following the same reasoning as in Step 2 of Section \ref{msw_sect}, the estimates for the $\xi_\cdot$ process easier as $\xi_\cdot \sim \nu_c$.  First, we limit the sum over $d$ to when $\|d\|$ is at most $\epsilon N$, where $N\uparrow\infty$ and then $\epsilon\downarrow 0$.  Next,  
 $\big[G_s\big(\frac{x+d}{N}\big)-G_s\big(\frac{x}{N}\big)\big]$ is replaced by $\nabla G_s\big(\frac{x}{N}\big)\cdot \frac{d}{N}$. Finally, the sum over $d$ is replaced by that when $\|d\|$ is at most $D$, which tends to infinity after $N$ diverges.
After this replacement, we have with probability tending to $1$ that
\begin{eqnarray*}
&&\frac{1}{N^n}\sum_{x\in \Z^n} \big|\eta_0(x) - \xi_0(x)\big|G_0\big(\frac{x}{N}\big) 
 + \int_0^t \frac{1}{N^n}\sum\limits_{x \in \Z^n}  |\eta_s(x)-\xi_s(x)|\partial_s G_s(\frac{x}{N})ds \\
&&\ \ + \int_0^t \frac{1}{N^n}\sum\limits_{x \in \Z^n} \sum\limits_{\|d\|=1}^{D} \frac{1}{\|d\|^{n-1+\alpha}} (gh_d(\eta_s(x))-gh_d(\xi_s(x)))O_{x,d}(\eta_s, \xi_s)\\
&&\ \ \ \ \ \ \ \ \ \ \ \ \ \ \ \ \ \ \ \ \ \ \ \ \ \ \times \nabla G_s\left(\frac{x}{N}\right) \cdot \frac{d}{\|d\|} ds \geq -\epsilon_0.
\end{eqnarray*}

\vskip .2cm
{\it Step 2.}
As in Step 3a in Section \ref{proofsection_alpha<1}, we may substitute $l$-averages for $|\eta_s(x)-\xi_s(x)|$ and $(gh_d(\eta_s(x))-gh_d(\xi_s(x)))O_{x,d}(\eta_s, \xi_s)$, where $l$ diverges after $N$ but before $D$, through   
a discrete integration-by-parts, the smoothness and compact support of $G$, as well as the particle bound Lemma \ref{particlebound}, and with respect to the $\xi_\cdot$ process that $\xi_\cdot \sim \nu_c$. 
Then, we have
\begin{eqnarray}
&&\frac{1}{N^n}\sum_{x\in \Z^n}|\eta_0(x) - \xi_0(x)| G_0\big(\frac{x}{N}\big) + \int_0^t \frac{1}{N^n}\sum\limits_{x \in \Z^n}  |\eta_s(x)-\xi_s(x)|^l \partial_s G_s\big(\frac{x}{N}\big)ds\nonumber\\
&&\ \ +\int_0^t \frac{1}{N^n}\sum\limits_{x \in \Z^n} \sum\limits_{\|d\|=1}^{D} \frac{1}{\|d\|^{n-1+\alpha}} \big[(gh_d(\eta_s(x))-gh_d(\xi_s(x)))O_{x,d}(\eta_s, \xi_s)\big]^l\nonumber\\
&&\ \ \ \ \ \ \ \ \ \ \ \ \ \ \ \ \ \ \ \ \ \ \ \ \ \ \times \nabla G_s\big(\frac{x}{N}\big) \cdot \frac{d}{\|d\|} ds \geq -\epsilon_0
\label{coup_eq1}
\end{eqnarray}
with high probability as $N, l$, and $D$ go to infinity. 

\vskip .2cm
{\flushleft \it Step 3a.}
We now begin to perform a `$1$-block' replacement in the last display, which will allow us to access the Young measure formulation. 

It is only here that we leverage the full form of the initial coupled distribution in order to treat the first term on the left-side of \eqref{coup_eq1}.  Let $A_1$ and $A_2$ be the set of sites $x$ in $\Z^n$ where
$\rho_0(x/N)\geq c$ and $\rho_0(x/N)<c$ respectively.  
Write, using the coupling, noting that $G$ is nonnegative, that
\begin{align*}
&\frac{1}{N^n}\sum_{x\in \Z^n}|\eta_0(x) - \xi_0(x)| G_0\big(\frac{x}{N}\big)\\
&\ = \frac{1}{N^n}\sum_{j=1}^2 \sum_{\stackrel{x\in A_j}{|x|\leq RN}}|\eta_0(x) - \xi_0(x)| G_0\big(\frac{x}{N}\big) \ = \   \sum_{j=1}^2 \big|\frac{1}{N^n}\sum_{\stackrel{x\in A_j}{|x|\leq RN}}\big[\eta_0(x) - \xi_0(x)\big]G_0\big(\frac{x}{N}\big)\big|.
\end{align*}

We now add and subtract $\rho_0(x/N) - c$ inside the square bracket.  Noting the compact support of $G$, we observe
$$
\tilde E_N| \frac{1}{N^n}\sum_{\stackrel{x\in A_j}{|x|\leq RN}}\big[\eta_0(x) - \rho_0(x/N)\big)]G_0\big(\frac{x}{N}\big)|^2
 \leq \frac{\|G\|}{N^{2n}}\sum_{\stackrel{x\in A_j}{|x|\leq RN}}{\rm Var}_{\nu_{\rho_0(x/N)}}(\eta(x)) = O(N^{-n}).
$$
A similar argument, using that $\xi_0$ has distribution $\nu_c$, works for the difference between $\xi_0(x)-c$.  
Hence, with high probability as $N\uparrow\infty$, we may bound above $\frac{1}{N^n}\sum_{x\in \Z^n}|\eta_0(x) - \xi_0(x)| G_0\big(\frac{x}{N}\big)$
by 
$$\frac{1}{N^n}\sum_{j=1}^2\big|\sum_{x\in A_j}\big(\rho_0(x/N) - c \big)G_0\big(\frac{x}{N}\big)\big|\leq
\frac{1}{N^n}\sum_{x\in \Z^n}|\rho_0(x/N) - c| G_0\big(\frac{x}{N}\big).$$

\vskip .2cm
{\flushleft \it Step 3b.}
Now, we will replace $|\eta_s(x)-\xi_s(x)|^l$ by $|\eta^l_s(x)-\xi^l_s(x)|$ in the the first integral, and $[(gh_d(\eta_s(x))-gh_d(\xi_s(x))) O_{x,d}(\eta_s, \xi_s)]^l$ by $(\Phi\Psi(\eta^l_s(x))-\Phi\Psi(\xi^l_s(x)))\sgn(\eta^l_s(x)-\xi^l_s(x))$ in the second integral of \eqref{coup_eq1}. 

Indeed, by the compact support of $\partial_s G$ and $\nabla G$, and $\sum_{\|d\|=1}^\infty \|d\|^{-(n-1+\alpha)}<\infty$, 
it will be enough to show that the expected integral over time of the quantities,
\begin{eqnarray*}
&&S_1=\frac{1}{N^n}\sum\limits_{|x| \leq RN}  \big||\eta_s(x)-\xi_s(x)|^l - |\eta^l_s(x)-\xi^l_s(x)| \big| \ \ {\rm and \ }\\
&&S_2=\frac{1}{N^n}\sum\limits_{|x| \leq RN} \big|[(gh_d(\eta_s(x))-gh_d(\xi_s(x)))O_{x,d}(\eta_s, \xi_s)]^l\\
&&\ \ \ \ \ \ \ \ \ \ \ \ \ \ \ \ \ \ \ \ \ \ \ \ \ \ \ -(\Phi\Psi(\eta^l_s(x))-\Phi\Psi(\xi^l_s(x))\sgn(\eta^l_s(x)-\xi^l_s(x))\big|,
\end{eqnarray*}
vanish in expectation as $N$ and then $l$ go to infinity, for each $\|d\|\leq D$.  

Divide now each of the sums above into two parts, $S_i = S_i^1 + S_i^2$, where $S_i^1$ is the part where $\eta_s$ and $\xi_s$ are ordered on $x+[-(l+D), l+D]^n$, and $S_i^2$ is the part where they are not. When $\eta_s$ and $\xi_s$ are ordered on the set $x+[-(l+D), l+D]^n$, we have
$|\eta_s(x)-\xi_s(x)|^l = |\eta^l_s(x)-\xi^l_s(x)|$,
and 
\begin{align*}
&\big|[(gh_d(\eta_s(x))-gh_d(\xi_s(x)))O_{x,d}(\eta_s, \xi_s)]^l-(\Phi\Psi(\eta^l_s(x))-\Phi\Psi(\xi^l_s(x))\sgn(\eta^l_s(x)-\xi^l_s(x))\big|\\
&\ \  \leq \big|(gh_d)^l(\eta_s(x))-\Phi\Psi(\eta^l_s(x))\big| + \big|(gh_d)^l(\xi_s(x))-\Phi\Psi(\xi^l_s(x))\big|.
\end{align*}
Therefore, the sum $S^1_1$ vanishes. 
But, by the $1$-block estimate Proposition \ref{1block_prop_alpha>1}, we have
$$\limsup_{l \to \infty} \limsup_{N \to \infty} \tilde{E}^N\int_0^t \frac{1}{N^n}\sum_{|x| \leq RN}\big| (gh_d)^l(\eta_s(x))-\Phi\Psi(\eta^l_s(x))\big|ds= 0,$$
and its counterpart with $\eta_\cdot$ replaced by $\xi_\cdot\sim \nu_c$ also vanishes. Hence, the expectation of the time integral of $S_2^1$ 
vanishes in the limit.

\vskip .2cm
{\flushleft \it Step 3c.}
When $\eta_s$ and $\xi_s$ are not ordered on $x+[-(l+D), l+D]^n$, as $h,\Psi\leq \|h\|_\infty$ are bounded, and $g, \Phi$ are Lipschitz, we have
\begin{align*}
&L_1=\big||\eta_s(x)-\xi_s(x)|^l - |\eta^l_s(x)- \xi^l_s(x)|\big|,\\
&L_2=\big|[(gh_d(\eta_s(x))-gh_d(\xi_s(x))) O_{x,d}(\eta_s, \xi_s)]^l \\
&\ \ \ \ \ \ \ \ - (\Phi\Psi(\eta^l_s(x))-\Phi\Psi(\xi^l_s(x)))\sgn(\eta^l_s(x)-\xi^l_s(x))\big|
\end{align*}
are both bounded by a constant times $(\eta^l_s(x)+\xi^l_s(x))$. Therefore, we may introduce the indicator function ${\bf 1}(\eta^l_s(x) \vee \xi^l_s(x) \leq A)$ when taking expectations, the error vanishing as $N,l\uparrow\infty$ by Lemma \ref{particlebound2}, and that $\xi_\cdot\sim \nu_c$. 

Once this indicator is introduced, both terms $L_1$, $L_2$ are bounded by a constant $C$, which allows further to introduce the indicator function ${\bf 1}(\eta^l_s(x+d) \vee \xi^l_s(x+d) \leq A)$, say by Lemma \ref{particlebound2} and that $\xi_\cdot\sim \nu_c$. 

Now, for $k=1,2$, we have
\begin{eqnarray*}
&&L_k{\bf 1}(\eta^l_s(x) \vee \xi^l_s(x) \leq A){\bf 1}(\eta^l_s(x+d) \vee \xi^l_s(x+d) \leq A)U_{x+[-(l+D), l+D]^n}(\eta_s, \xi_s) \\
&&\ \ \ \ \ \ 
\leq C U_{x+[-(l+D), l+D]^n}(\eta_s, \xi_s).
\end{eqnarray*}
Therefore, to complete the $1$-block replacement, it is enough to show that
\begin{equation}
\label{entropyboundstep4}
\limsup_{N \to \infty} \tilde{E}^N\int_0^t \frac{1}{N^n}\sum\limits_{|x| \leq RN} U_{x+[-(l+D), l+D]^n}(\eta_s, \xi_s) ds =0.
\end{equation}

\vskip .2cm

{\flushleft  \it Step 3d.}
Recall that $U_\Lambda(\eta, \xi)$ indicates when $\eta$ and $\xi$ are not ordered on $\Lambda$, and also that $U_{x, d}(\eta, \xi)= U_{\{x,x+d\}}(\eta, \xi)$. We then have the bound, 
$$U_{x+[-(l+D), l+D]^n}(\eta_s, \xi_s) \leq \sum_{|y|\leq l+D}\sum\limits_{|d|\leq 2(l+D)} U_{x+y, d}(\eta_s, \xi_s),$$
from which it follows for each $l$ that
$$\sum_{|x|\leq RN} U_{x+[-(l+D), l+D]^n}(\eta_s, \xi_s) \leq \sum_{|y|\leq l+D}\sum\limits_{|d|\leq 2(l+D)}\sum_{|x|\leq R^+N} U_{x, d}(\eta_s, \xi_s), $$
for large enough $N$ where $R^+ > R$.
However, by the Ordering Lemma \ref{orderinglemma}, we have for each $l$ and $d$ that
$\tilde{E}^N\int_0^t \frac{1}{N^n}\sum_{|x| \leq R^+N} U_{x, d}(\eta_s, \xi_s) ds$ vanishes as $N\uparrow\infty$.  Hence,
\eqref{entropyboundstep4} holds.
and the $1$-block replacement follows.

In particular, we have
\begin{eqnarray}
&&\frac{1}{N^n}\sum_{x\in \Z^n} \big|\rho_0(x/N) - c\big|G_0\big(\frac{x}{N}\big) + \int_0^t \frac{1}{N^n}\sum\limits_{x \in \Z^n}  |\eta^l_s(x)-\xi^l_s(x)|\partial_s G_s\left(\frac{x}{N}\right)ds\nonumber\\
&&\ \  + \int_0^t \frac{1}{N^n}\sum\limits_{x \in \Z^n} \sum\limits_{\|d\|=1}^{D} \frac{1}{\|d\|^{n-1+\alpha}} (\Phi\Psi(\eta^l_s(x))-\Phi\Psi(\xi^l_s(x)))\nonumber\\
&&\ \ \ \ \ \ \ \ \ \ \ \ \ \ \ \ \ \times \sgn(\eta^l_s(x)- \xi^l_s(x))\nabla G_s\left(\frac{x}{N}\right) \cdot \frac{d}{\|d\|} ds \geq -\epsilon_0
\label{last11}
\end{eqnarray}
with high probability as $N, l$, and $D$ go to infinity. 
To recover \eqref{coupmicroentro} from \eqref{last11}, we may group together the terms involving $d$, and remove the bound $\|d\|\leq D$, by appealing to the argument in Step 3a in Subsection \ref{msw_sect}. Then, the sum on $d$ is replaced by $\gamma_\alpha \partial_{\textbf{1}(n)} G_s(\frac{x}{N})$. \qed

\section{$L^1$ mass bound: Proof of Theorem \ref{mass_thm}}
\label{L^1mass_sect}
We leverage the weak formulation of the entropy condition.

\vskip .2cm

{\it Step 1.} Consider a test function in form $G(s,u)=H(s)G(u)$ for $G$ nonnegative, and $c= \rho^*$.  Define $V_G(s)=\int_{\R^n} \int_0^\infty G(u)|\lambda - \rho^*| \rho(s,u, d\lambda) du$. By the `Mass Bounding' Lemma \ref{VGbound} shown below, $V_G$ is finite.  Moreover, by the measure weak entropy condition inequality,
\begin{eqnarray*}
&&-\int_0^\infty   \partial_sH(s)V_G(s) ds\\
&&\ \ \ \ \ \  \leq H(
0)V_G(0) + \gamma_\alpha\int_0^\infty H(s) \int_{\R^n} \int_0^\infty \partial_{\textbf{1}(n)} G\left(u\right)  q(\lambda,\rho^*)\rho(s,u, d\lambda)duds,
\end{eqnarray*}
where we recall $q(\lambda,c) = \sgn(\lambda -c)(\Phi\Psi(\lambda) - \Phi\Psi(c))$.
Since $\Phi,\Psi$ are Lipschitz, and also $\Psi\leq \|h\|_\infty$, we have
$|q(\lambda, c)|\leq |(\Phi(\lambda)-\Phi(c))\Psi(\lambda) + \Phi(c)(\Psi(\lambda)-\Psi(c))|\leq C|\lambda -c|$ where the constant depends on $c$.  Hence, $|q(\lambda,\rho^*)|\leq  
 C|\lambda-\rho^*|$. 
 Then,
$$-\int_0^\infty   \partial_sH(s)V_G(s) ds \leq H(0)V
_G(0) + \gamma_\alpha C(\rho^*) \int_0^\infty |H(s)| V_{|\partial_{\textbf{1}(n)} G|}(s)ds.$$
\vskip .2cm

{\it Step 2.} We now define a sequence of test functions $\{H_i\}$ on $\R$. 
With respect to $0\leq t< T$ and $0\leq \delta\leq T-t$, let $H_0(s)=1$ up to $s=t-\delta$, then decreasing to $0$ by $s=t+\delta$. We may do this in such a way that
$-\partial_s H_0(s)$ is positive on $(t-\delta,t+\delta)$ and weakly approaches a delta function at $t$.  For instance, we can take $-\partial_sH_0(s)$ as the linear interpolation between $(t-\delta, 0)$, $(t-\delta+\gamma, L)$, $(t+\delta-\gamma, L)$ and $(t+\delta, 0)$ where
$L= (2\delta -\gamma)^{-1}$ for $\gamma<\delta/2$.

For each $i\geq 1$, define $H_{i+1}(s)= \int_s^{t+\delta} H_i(u)du$. Therefore, each $H_i(s)$ is nonnegative, vanishing for $s>t+\delta$. We then have $-\partial_s H_{i+1}(s)= |H_i(s)|$.  Note also
 $H_i(0) \leq (t+\delta)^i/i!$.

 Define $G_i(u)=1$ on $[-(i+1)R, (i+1)R]^n$, decreasing to zero within $[-(i+2)R, (i+2)R]^n$, so that $|\partial_{\textbf{1}(n)}G_i(u)| \leq 2/R$. Then, $|\partial_{\textbf{1}(n)}G_i(u)| \leq (2/R) G_{i+1}(u)$.   Here, the limit of these functions, $G_\infty(u)\equiv 1$.
  
 With respect to $H=H_i$ and $G=G_i$, we have
$$\int_0^\infty  - \partial_sH_i(s)V_{G_i}(s) ds \leq H_i(0)V
_{G_i}(0) + \frac{2\gamma_\alpha C(\rho^*)}{R} \int_0^\infty - \partial_sH_{i+1}(s) V_{G_{i+1}}(s)ds.$$

\vskip .2cm
{\it Step 3.}
Iterating the above inequality $k$ times, 
starting with $i=0$, gives
$$-\int_0^\infty  \partial_sH_0(s)V_{G_0}(s) ds \leq \sum_{i=0}^{k-1} \big(\frac{2\gamma_\alpha C}{R}\big)^iH_i(0)V_{G_i}(0) + \big(\frac{2\gamma_\alpha C}{R}\big)^{k} \int_0^\infty - \partial_sH_{k}(s) V_{G_{k}}(s)ds.$$

Since $V_{G_k}(s) \leq {\rm ess\,sup}_{0\leq t \leq T}V_{G_k}(t)$,  
 $V_{G_i}(0)\leq V_{G_\infty}(0)$, and $H_i(0) \leq (t+\delta)^i/i!$, we obtain
\begin{eqnarray}
&&-\int_0^\infty   \partial_sH_0(s)V_{G_0}(s) ds \leq \sum_{i=0}^{k-1} \frac{(2\gamma_\alpha CT/R)^i}{i!}V_{G_\infty}(0) \nonumber\\
&&\ \ \ \ \ \ \ \ \ \ \ \ \ \ \ \ \ \ \ \ \ \ + \frac{(2\gamma_\alpha CT/R)^k}{k!} {\rm ess\,sup}_{0\leq t \leq T}V_{G_k}(t).
\label{step3_equation}
\end{eqnarray}

\vskip .2cm
{\it Step 4.}
Choose now $k=n+2$.
The supremum over $0<\gamma<\delta/2$, $0<\delta\leq T-t $ and $0\leq t< T$ of the left-side of \eqref{step3_equation}, $Q^*$ a.s., increases by monotone convergence as $R\uparrow\infty$ to 
$$\sup_{0<\gamma<\delta/2}\sup_{0<\delta\leq T-t  }\sup_{0\leq t<T}\int_0^\infty  -\partial_sH_0(s)V_{G_\infty}(s) ds.$$

To capture the limit of the right-side,
note
$$V_{G_\infty}(0) = \int_{\R^n} \int_0^\infty |\lambda-\rho^*|\rho(0,u, d\lambda) du = \int_{\R^n} |\rho_0(u)-\rho^*|du < \infty.$$
Then, the first term on the right-side of \eqref{step3_equation} converges to $V_{G_\infty}(0)$ as $R\uparrow\infty$.

However, by the
`Mass Bounding' Lemma \ref{VGbound}, we have $E_{Q^*} {\rm ess\,sup}_{0\leq t\leq T}V_{G_{n+2}}(t) = O(R^n)$, and so
\begin{eqnarray*} 
 \frac{(2\gamma_\alpha CT/R)^{n+2}}{(n+2)!}E_{Q^*}\big[{\rm ess\,sup}_{0\leq t \leq T}V_{G_{n+2}}(t)\big] = O(R^{-2}).
\end{eqnarray*}
Hence, $Q^*$ a.s., by Borel-Cantelli lemma, as $R\uparrow\infty$,
the second term on the right-side of \eqref{step3_equation} vanishes.

\vskip .2cm

{\it Step 5.} 
Therefore, we have
$\sup_{0<\gamma<\delta/2}\sup_{0<\delta\leq T-t}\sup_{0\leq t< T}
\int_0^\infty  - \partial_sH_0(s)V_{G_\infty}(s) ds\leq V_{G_\infty}(0)$,
with respect to a $Q^*$ probability $1$ set.
Moreover, on this set, as $-\partial_s H_0$ is positive on $(t-\delta,t+\delta)$, we have that
$V_{G_\infty}$ is locally integrable on $[0,T]$.  Also, by Fatou's lemma, for each $t$ and small enough $\delta>0$, we have
$$\frac{1}{2\delta}\int_{t-\delta}^{t+\delta} V_{G_\infty}(s)ds \leq \liminf_{\gamma\downarrow 0} \int_0^\infty  - \partial_sH_0(s)V_{G_\infty}(s) ds \leq V_{G_\infty}(0).$$
In fact, for each Lebesgue point $t$ of $V_{G_\infty}$,  
as $\delta\downarrow 0$, we have
$V_{G_\infty}(t) \le V_{G_\infty}(0)$.
 We conclude, as Lebesgue points are dense, that $Q^*$ a.s.
\begin{equation*}
{\rm ess\,sup}_t \int_{\R^n} \int_0^\infty |\lambda-\rho^*|\rho(t,u, d\lambda) du \leq V_{G_\infty}(0),
\end{equation*}
finishing the argument. \qed

\subsection{Mass bounding lemma}
\label{mass_sect}

The following result bounds the mass in finite regions.

\blem \label{VGbound} Let $G$ be a nonnegative function with support in $[-R, R]^n$ such that $|G|\leq 1$. For every limit point $Q^*$ and $c\in \R$, we have 
$$E_{Q^*}\Big[{\rm ess\,sup}_{0\leq t \leq T} \int_{\R^n} G(u) \int_0^\infty |\lambda-c| \rho(t, u, d\lambda)du\Big] = O(R^n).$$ 
\elem

\bpf First we bound $|\lambda-c|$ by $\lambda+ |c|$.  Since $\rho(t,u,d\lambda)$ 
 is a probability measure, we have
$\int_{\R^n} G(u) \int_0^\infty |c| \rho(t, u, d\lambda)du = O(R^n)$. 
Therefore, we only need to prove 
\begin{equation}
\label{mass_0}
E_{Q^*}\Big[{\rm ess\,sup}_{0\leq t\leq T} \int_{\R^n} G(u) \int_0^\infty \lambda \rho(t,u, d\lambda) du\Big] = O(R^n).
\end{equation}

To this end, for $RN \geq l$, note 
\begin{equation*}
\label{mass_1}
\lt \pi^{N,l}_t, G(u)\lambda \rt \leq \frac{1}{N^n} \sum_{|x| \leq RN} \eta^l_t(x) \leq \frac{1}{N^n} \sum_{|x| \leq 2RN} \eta_t(x) \leq \lt \pi^N_t, G_1\rt,
\end{equation*}
where $G_1$ equals $1$ on $[-2R, 2R]^n$, and decreases to zero within $[-3R, 3R]^n$.
By the definition of the martingale $M^{N,G_1}_t$, we have
$$E^N \sup_{0 \leq t\leq T} \lt  \pi^{N,l}_t, G_1\rt \leq E^N \lt \pi^N_0, G_1\rt+ E^N\int_0^T |N L_N \lt \pi^N_s, G_1\rt| ds + E^N \sup_{0 \leq t \leq T} |M^{N,G_1}_t|.$$
By our initial conditions, $E^N\lt\pi^{N}_0, G_1 \rt \leq (4R+1)^n \|\rho_0\|_{\infty}$, and 
by  Lemma \ref{alpha<1martbound}, we have $E^N\int_0^T |N L_N \lt \pi^N_s, G_1\rt| ds =O(R^n)$, independent of $N$ and $l$. 
Also, by Doob's inequality and Lemma \ref{alpha<1quadbound}, $E^N\sup_{0 \leq t \leq T} |M^{N,G_1}_t| \leq 4
E^N \lt M^{N, G_1} \rt_T =O(N^{-n})$.
Therefore, for all large $N$, 
we have
$E^{N,l}{\rm ess\,sup}_{0 \leq t\leq T} \lt \pi^{N,l}_t, G(u)\lambda \rt = O(R^n)$.

Finally, as ${\rm ess\,sup}_{0 \leq t\leq T} \lt \pi^{N,l}_t, G(u)\lambda \rt$ is a lower semi-continuous function of $\pi^{N,l}$, we may take subsequential limits as $N,l\uparrow\infty$, for which $Q^{N,l}\Rightarrow Q^*$, 
to obtain \eqref{mass_0}.  \epf

\section{Initial conditions: Proof of Theorem \ref{init_cond_thm}}
\label{init_sect}

The strategy is to approximate the initial density $\rho_0$ in compact sets via the weak form of the entropy inequality.

\vskip .2cm

{\it Step 1.}
Since $\rho_0$ is a continuous function that equals a constant $\rho^*$ outside of a compact set $[-R, R]^n$, it is uniformly continuous. 
Fix a $\delta= (\delta_0, \ldots, \delta_0)$ with $0<\delta_0<1$. Divide $\R^n$, regularly, into countably many overlapping hyper-rectangles $[a_i-\delta, b_i+\delta]=\prod_{j=1}^n [a_{i,j} - \delta_0, b_{i,j}+\delta_0]$ such that $\cup_{i=1}^\infty [a_i-\delta, b_i+\delta] = \R^n$, and the $[a_i, b_i]$ are disjoint. Finitely many of these hyper-rectangles cover $[-R, R]^n$.  The parameter $\delta$ may be chosen so that
$\rho_0$ varies at most $\epsilon_0>0$ on each hyper-rectangle.

For each hyper-rectangle, we construct a nonnegative smooth bump function, $G_i(u)$, that is $1$ on $[a_i,b_i]$ and decreases to $0$ outside of $[a_i-\delta, b_i+\delta]$. The $\{G_i\}$ may be constructed such that $\sum_{i=1}^\infty G_i(u)=1$ for all $u\in \R^n$, and $\max_i\|\partial_{{\bf 1}(1)} G_i\|_\infty$ is bounded. We also choose constants $c_i=\min \{\rho_0(u): u \in [a_i-\delta,b_i+\delta]$ so that $|\rho_0(u)-c_i|\leq\epsilon$ for all $u \in [a_i-\delta,b_i+\delta]\}$. By the triangle inequality, $|\lambda- \rho_0(u)| \leq |\lambda- c_i|- |\rho_0(u)-c_i|+ 2\epsilon$ on any hyper-rectangle that intersects $[-R, R]^n$; on the other hyper-rectangles, as $\rho_0 = \rho^*$, we have $|\lambda- \rho_0(u)|= |\lambda- c_i|-|\rho_0-c_i|$. 

Note that $(2R)^n$ is the volume of $[-R, R]^n$, $\rho(s,u,d\lambda)$ is a probability measure (cf. Step 1 in Section \ref{proof_section_alpha>1}), and $|\lambda- c_i|-|\rho_0-c_i| = |\lambda - c_i|\geq 0$ for all but finitely many hyper-rectangles.  Then, a Fubini-Tonelli theorem may be applied, so that
\begin{eqnarray*}
&& E_{Q^*}\Big[\frac{1}{t}\int_0^t \int_{\R^n} \int_0^\infty |\lambda- \rho_0(u)|\rho(s, u, d\lambda) duds\Big]\\
&&\ \ \leq  E_{Q^*}\Big[\frac{1}{t}\int_0^t \sum_{i=1}^\infty \int_{\R^n} G_i(u) \int_0^\infty (|\lambda- c_i|-|\rho_0(u)-c_i|)\rho(s, u, d\lambda) duds\Big]  +2\epsilon (2R)^n\\
&&\ \ \ = \sum_{i=1}^\infty E_{Q^*} \Big[\frac{1}{t}\int_0^t \int_{\R^n} G_i(u) \int_0^\infty (|\lambda- c_i|-|\rho_0(u)-c_i|)\rho(s, u, d\lambda) duds\Big]  +2\epsilon (2R)^n.
\end{eqnarray*}

\vskip .2cm

{\it Step 2.}
Suppose, for all $i$, that
\begin{eqnarray}
\label{initcondlem}
&&\limsup\limits_{t \downarrow 0}  E_{Q^*} \Big| \frac{1}{t} \int_0^t \int_{\R^n} G_i(u) \int_0^\infty  (|\lambda - c_i|-|\rho_0(u)- c_i|)\rho(s,u,d\lambda)duds \Big|=0,\\
&& \ {\rm and \ } \ \ \sum_{i=1}^\infty \sup_{0<t<T} E_{Q^*}  \Big|\frac{1}{t} \int_0^t \int_{\R^n} G_i(u) \int_0^\infty  (|\lambda - c_i|-|\rho_0(u)- c_i|)\rho(s,u,d\lambda)duds\Big| <\infty.\nonumber
\end{eqnarray}
Then, by Fatou-Lebesgue lemma , we would have
\begin{eqnarray*}
&&E_{Q^*} \Big[\liminf\limits_{t \downarrow 0} \frac{1}{t} \int_0^t \int_{\R^n}\int_0^\infty  |\lambda - \rho_0(u)|\rho(s,u,d\lambda)duds\Big]\\
&& \leq \liminf\limits_{t \downarrow 0} E_{Q^*}\frac{1}{t}\int_0^t \int_{\R^n} \int_0^\infty |\lambda- \rho_0(u)|\rho(s, u, d\lambda) du ds\leq 2\epsilon_0 (2R)^n,
\end{eqnarray*}
from which Theorem \ref{init_cond_thm} would follow as $\epsilon_0>0$ is arbitrary.

\vskip .2cm
{\it Step 3.}
To finish the proof, we establish \eqref{initcondlem}.  As discussed in Step 1 in Section \ref{proof_section_alpha>1}, with respect to $Q^*$, initially $\rho(0, u,d\lambda) = \delta_{\rho_0(u)}$, and $\rho(s,u,d\lambda)$ is a probability measure.  Then, with respect to a test function $G(s,u) = H(s)G(u)$ with $G$ nonnegative, we have
$H(0)\int_{\R^n}G(u)\int_0^\infty |\rho_0(u)-c|\rho(s,u,d\lambda)du = H(0)\int_{\R^n}\int_0^\infty G(u)|\lambda -c|\rho(0,u,d\lambda)du$.
Hence, we can write the measure weak formulation of the entropy condition in Section \ref{proof_section_alpha>1} as
\begin{eqnarray*}
&&\int_0^\infty \partial_s H(s) \int_{\R^n} G(u) \int_0^\infty  (|\lambda - c|-|\rho_0(u)- c|)\rho(s,u,d\lambda)duds\\
&&\ \ \  + \gamma_\alpha\int_0^\infty H(s) \int_{\R^n} \int_0^\infty \partial_{\textbf{1}(n)} G\left(u\right)  q(\lambda, c)\rho(s,u, d\lambda)duds \geq 0,
\end{eqnarray*}
recalling 
$q(\lambda, c)= \sgn(\lambda-c)(\Phi\Psi(\lambda)- \Phi\Psi(c))$.  
As in Step 1 of Section \ref{L^1mass_sect},
we have $|q(\lambda, c)|\leq C_2|\lambda -c|$ where the constant $C_2$ depends on $\|\rho_0\|_\infty$.  
Therefore, putting the first term on the other side of the inequality,
\begin{eqnarray}
&&-\int_0^\infty \partial_s H(s) \int_{\R^n} G(u) \int_0^\infty  (|\lambda - c|-|\rho_0(u)- c|)\rho(s,u,d\lambda)duds\nonumber\\
&& \leq \gamma_\alpha C_2\int_0^t |H(s)| \int_{\R^n} \int_0^\infty |\partial_{\textbf{1}(n)} G\left(u\right)|  |\lambda-c|\rho(s,u, d\lambda)duds.
\label{initcond_step3}
\end{eqnarray}

\vskip .2cm
{\it Step 4.}  Recall, that
$E_{Q^*} {\rm ess\,sup}_{0\leq t\leq T} \int_{\R^n} \int_0^\infty G(u) |\lambda -c|\rho(t,u,d\lambda)du <\infty$, by the `Mass Bounding Lemma' \ref{VGbound}, for all nonnegative $G$'s with compact support.
Consider, with respect to a small $\delta>0$, a smooth $H$ such that $\partial_s H(s) = -t^{-1}$ for $0\leq s\leq t-\delta$, linearly interpolates from $(t-\delta, -t^{-1})$ to $(t,0)$, and vanishes for $s\geq t$.  Taking $H(0)= 1-\delta/(2t)$, we have that $H$ vanishes for $s\geq t$.  

Since, $|H(s)|\leq 1$, the right-side of the inequality \eqref{initcond_step3} is bounded by taking a supremum over time, and the left-side, by dominated convergence, as $\delta\downarrow 0$, tends to
\begin{eqnarray*}
&&\frac{1}{t}\int_0^t \int_{\R^n} G(u) \int_0^\infty  (|\lambda - c|-|\rho_0(u)- c|)\rho(s,u,d\lambda)duds\\
&&\leq \gamma_\alpha C_2t \ {\rm ess\,sup}_{0 \leq s \leq T} \int_{\R^n} \int_0^\infty |\partial_{\textbf{1}(n)} G\left(u\right)|  |\lambda-c|\rho(s,u, d\lambda)du.
\end{eqnarray*}

As $|\partial_{\textbf{1}(n)} G|$ is compactly supported,
noting the `Mass Bounding Lemma' \ref{VGbound} again, the expected value of the right-side of the above display vanishes as $t$ goes to zero. Therefore, the first line of \eqref{initcondlem} holds.

To obtain the second line of \eqref{initcondlem}, instead of bounding the right-side of \eqref{initcond_step3} by a supremum, bound it by
increasing the time integration to $[0,T]$.  Then, we have
\begin{eqnarray*}
&&\sup_{0<t< T}\Big|\frac{1}{t}\int_0^t \int_{\R^n} G_i(u) \int_0^\infty  \big(|\lambda - c_i|-|\rho_0(u)- c_i|\big)\rho(s,u,d\lambda)duds\Big|\\
&&\leq \gamma_\alpha C_2\int_0^T \int_{\R^n} \int_0^\infty |\partial_{\textbf{1}(n)} G_i\left(u\right)|  |\lambda-c_i|\rho(s,u, d\lambda)duds.
\end{eqnarray*}
For only finitely many $i$ does $c_i$ differ from $\rho^*$ and $|\rho_0(u)-c_i|>0$.  Also, by the comment in the previous paragraph, for each $i$, the expected value of the right-side of the above display is bounded. Note now, by the regular division, that the support of each $\partial_{\textbf{1}(n)}G_i$ is overlapped by the support of at most an uniformly bounded number, in terms of the covering, of other $\{\partial_{\textbf{1}(n)}G_j\}$.  Note also, from construction, that $\partial_{\textbf{1}(n)}G_j$ is uniformly bounded in $j$.   Also, from Theorem \ref{mass_thm}, we have that
$E_{Q*}\big[\int_0^T \int_{\R^n}\int_0^\infty |\lambda - \rho^*|\rho(s,u,d\lambda)duds\big] \leq T \int_{\R^n} |\rho_0(u)-\rho^*|du<\infty$.  Hence, summability in \eqref{initcondlem} follows, and the proof of Theorem \ref{init_cond_thm} is complete. \qed

\vskip .1cm
In passing, we remark that this proof, making use of the weak formulation of the entropy condition, seems new and more direct than proofs in \cite{Rezakhanlou} and \cite{Kipnis} which introduce types of particle couplings in the finite-range setting, without going to the continuum equation.  We note, in the PhD thesis \cite{DS}, an alternate argument for the first line of \eqref{initcondlem} through a simpler and different particle coupling will be found.

\appendix
\section{Proof of the Ordering Lemma \ref{orderinglemma}}
\label{ordering_section}

{\it Step 1.} We now show the first part of the lemma.
Let $G_s(u)$ be a nonnegative smooth function that is $1$ on hyper-rectangle $[a,b]=\prod_{j=1}^n[a_j,b_j]$ and decreases to $0$ outside of $[a-\delta, b+\delta]=\prod_{j=1}^n [a_j-\delta_j,b_j+\delta_j]$ where $\delta= (\delta_1, \ldots, \delta_n)$ with $\delta_i > 0$ and $\|\delta\|<1$.

Then, noting the computation of $N\tilde{L}_N\lt \tilde{\pi}^N_s, G_s \rt$ in \eqref{coupled_gen} and \eqref{alpha>1_generator}, we have
\begin{eqnarray*}
&&\frac{N}{N^n} \sum\limits_{\|d\|=1}^\infty \frac{1}{\|d\|^{n+\alpha}} \sum\limits_{\substack{x \in [a, b]N\\ x+d \in [a,b]N}} 
|gh_d(\eta_s(x))-gh_d(\xi_s(x)| U_{x,d}(\eta_s, \xi_s)\\
&&\leq \frac{N}{N^n} \sum\limits_{x \in \Z^n} \sum\limits_{\|d\|=1}^\infty \frac{1}{\|d\|^{n+\alpha}} (gh_d(\eta_s(x))-gh_d(\xi_s(x))O_{x,d}(\eta_s, \xi_s) \big[G_s\big(\frac{x+d}{N}\big)-G_s\big(\frac{x}{N}\big) \big]\\
&&\ \ \ \ \ \ \ \ - N\tilde{L}_N\lt \tilde{\pi}^N_s, G_s \rt \ =: J_1 - N\tilde{L}_N\lt \tilde{\pi}^N_s, G_s \rt.
\end{eqnarray*}
\vskip .2cm

{\it Step 2.}
Note, the expression $|gh_d(\eta_s(x))-gh_d(\xi_s(x)|\leq \kappa \|h\|_\infty (\eta_s(x)+\xi_s(x))$ by \eqref{LB},
and $|O_{x,d}(\eta_s, \xi_s)|\leq 1$.  Hence,
$$J_1\leq \kappa \|h\|_\infty \frac{N}{N^n} \sum\limits_{x \in \Z^n} \sum\limits_{\|d\|=1}^\infty \frac{1}{\|d\|^{n+\alpha}}   (\eta_s(x)+ \xi_s(x))\big|G_s\big(\frac{x+d}{N}\big)-G_s\big(\frac{x}{N}\big) \big|.$$

We now split the sum over $d$ into two parts, namely when $\|d\|\leq N$ and $\|d\|> N$, and write $J_1 = J_{11} + J_{12}$ accordingly. When $\|d\|\leq N$, we bound $|G_s\left(\frac{x+d}{N}\right)-G_s\left(\frac{x}{N}\right)| \leq \|\nabla G\| \cdot \|d\|/N{\bf 1}(x \in [(a-\delta-1)N, (b+\delta+1)N])$. 
Then,
$$J_{11} \leq \kappa \|h\|_\infty \sum\limits_{\|d\|=1}^{N} \frac{1}{\|d\|^{n-1+\alpha}} \|\nabla G\|  \frac{1}{N^n} \sum' (\eta_s(x)+ \xi_s(x))$$
where $\sum'$ refers to a sum over $O(N^n)
$ values of $x$.
By Lemma \ref{particlebound}, and that the $\xi_s$ process starts in the invariant measure $\nu_c$, we have
$\tilde{E}^N\sum' (\eta_s(x)+ \xi_s(x))\leq 2K_0N^n$ say.  Therefore, 
$\tilde E^N J_{11} \leq \tilde{C}_1 := \kappa \|h\|_\infty\sum_{\|d\|=1}^\infty \frac{1}{\|d\|^{n-1+\alpha}} \|\nabla G\| 2K_0$, finite when $\alpha>1$.

On the other hand, when $\|d\|> N$, we have $|G_s\left(\frac{x+d}{N}\right)-G_s\left(\frac{x}{N}\right)| \leq {\bf 1}(x\in [(a-\delta)N, (b+\delta)N] \cup x+d \in [(a-\delta)N, (b+\delta)N])$.
Then, in terms of a sum $\sum''$ over $O(N^n)$ sites,
$$J_{12}\leq \kappa \|h\|_\infty \sum\limits_{\|d\|> N} \frac{1}{\|d\|^{n-1+\alpha}} \frac{1}{\|d\|N^{n-1}} \sum''(\eta_s(x)+ \xi_s(x)).
$$
As $\|d\|\geq N$, we have
$\tilde{E}^N\left[\frac{1}{\|d\|N^{n-1}}\sum''(\eta_s(x)+ \xi_s(x))\right] \leq 2K_0'$ say, uniformly in $d$,
by Lemma \ref{particlebound} and that $\xi_\cdot\sim \nu_c$. Then, $\tilde E^N J_{12} \leq 
\tilde{C}_2:=\kappa \|h\|_\infty \sum\limits_{\|d\|=1}^\infty \frac{1}{\|d\|^{n-1+\alpha}}2K_0'$.

\vskip .2cm

{\it Step 3.} It follows that $\tilde E^N J_1 \leq \tilde C = \tilde C_1 + \tilde C_2$. Therefore,
\begin{eqnarray}
&&\frac{1}{N^n}\tilde E^N \int_0^t \sum\limits_{\|d\|=1}^\infty \frac{1}{\|d\|^{n+\alpha}} \sum\limits_{\substack{x \in [a, b]N\\ x+d \in [a,b]N}}
|gh_d(\eta_s(x))-gh_d(\xi_s(x)| U_{x,d}(\eta_s, \xi_s) ds\nonumber\\
&& \ \ \ \ \ \leq 
\frac{t\tilde{C}}{N}- \tilde{E}^N\int_0^t \tilde{L}_N\lt \tilde{\pi}^N_s, G_s \rt ds. 
\label{ordered_step3}
\end{eqnarray}

Consider the mean-zero martingale $\tilde M^{N,G}_t$ where $G_s(u) =G(u)\geq 0$ for $s\leq t$, vanishing before time $T$.  As $\partial_sG =0$ for $s\leq t$, and $\lt \tilde{\pi}^N_t, G \rt \geq 0$, we have
$$- \int_0^t N\tilde{L}_N\lt \tilde{\pi}^N_s, G_s \rt ds\leq  \tilde{M}^{N,G}_t +\lt \tilde{\pi}^N_0, G_0 \rt.$$
Therefore, \eqref{ordered_step3} is bounded by
$$\frac{t\tilde{C}}{N} + \frac{1}{N^{n+1}}\tilde{E}^N\sum_{x\in [a-\delta,b+\delta]N} (\eta_0(x)+\xi_0(x)).$$
As the expectation is of order $O(N^n)$ by Lemma \ref{particlebound} and that $\xi_\cdot\sim \nu_c$, the last display vanishes as $N\uparrow\infty$.  This completes the proof of the first part of Lemma \ref{orderinglemma}.
\vskip .2cm

{\it Step 4.}
We now show the second part of Lemma \ref{orderinglemma}.  In general, $gh_d(\eta_s(x)) - gh_d(\xi_s(x))$ may not vanish, and so the first part is not coercive.  To work around this issue, we would like to introduce the indicator function ${\bf 1}(\eta_s(x) \vee \xi_s(x) \vee \eta_s(x+d) \vee \xi_s(x+d) < A)$ into the associated expectation.   This is justified if we show
that 
$$\lim_{A\uparrow\infty}\lim_{N\uparrow\infty}\tilde{E}^N\int_0^t \frac{1}{N^n} \sum\limits_{|x|\leq RN} U_{x,d}(\eta_s, \xi_s){\bf 1}(\eta_s(x) \vee \xi_s(x) \vee \eta_s(x+d) \vee \xi_s(x+d) > A) ds=0.$$
The expectation above is bounded by the sum of
$E^N\int_0^t \frac{1}{N^n} \sum_{|x|\leq RN} {\bf 1}(\eta_s(x) > A)  ds$,
and three other expectations containing the
 indicator functions  ${\bf 1}(\xi_s(x)>A)$,  ${\bf 1}(\eta_s(x+d)>A)$, and ${\bf 1}(\xi_s(x+d) > A))$. 

By the entropy inequality \eqref{entropy_ineq}, the first expectation is bounded by
$$\frac{1}{\gamma N^n} \Big(O(N^n) \ln E_{\nu_{\rho^*}} \big[ e^{\gamma {\bf 1}(\eta(0)>A)} \big] + CN^n\Big).$$
After $N$ and $A$ go to infinity, the limit is $C/\gamma$ which vanishes as $\gamma\uparrow\infty$.  The other three terms are similarly analyzed, using $\xi_\cdot\sim \nu_c$ when $\xi_\cdot$ is involved.

Therefore, it will be enough to prove, for each $A$, that
$$\limsup_{N \to \infty} \tilde{E}^N\int_0^t \frac{1}{N^n} \sum\limits_{|x|\leq RN} U_{x,d}(\eta_s, \xi_s){\bf 1}(\eta_s(x) \vee \xi_s(x) \vee \eta_s(x+d) \vee \xi_s(x+d) < A) ds= 0.$$

\vskip .2cm
{\it Step 5.}
Let $I^{x,d}_{m_1, m_2, k_1, k_2}(\eta, \xi)= {\bf 1}(\eta(x)=m_1, \xi(x)=m_2, \eta(x+d)=k_1, \xi(x+d)=k_2)$. Then, 
$$U_{x,d}(\eta, \xi){\bf 1}(\eta(x) \vee \xi(x) \vee \eta(x+d) \vee \xi(x+d) < A) \leq \widetilde\sum I^{x,d}_{m_1, m_2, k_1, k_2}(\eta, \xi)$$
where the sum $\widetilde\sum$ is over all $m_1, m_2, k_1, k_2$ less than $A$ such that $\eta$ and $\xi$ will be not ordered on the sites $x$ and $x+d$. Since this is a finite sum, it will be enough to prove that 
\begin{equation}
\label{I_eq}
\limsup_{N \to \infty} \tilde{E}^N\int_0^t \frac{1}{N^n} \sum\limits_{|x|\leq RN} I^{x,d}_{m_1, m_2, k_1, k_2}(\eta_s, \xi_s) ds = 0,
\end{equation}
for each $m_1, m_2, k_1, k_2$ indexed in $\widetilde\sum$.

From the `proven first part', in other words that \eqref{ordered_step3} vanishes, we note if $g(m_1)h(k_1)\neq g(m_2)h(k_2)$ then \eqref{I_eq} holds.

Recall $M_0= \min\{k: h(k)=0\}$ is the maximum possible number of particles at a site, with the convention that $M_0=\infty$ if $h(k)$ is never zero. 
If $k_1$ or $k_2$ is greater than $M_0$, then $I^{x,d}_{m_1, m_2, k_1, k_2}(\eta_s, \xi_s)$ identically vanishes, and \eqref{I_eq} holds trivially.

We note \eqref{I_eq} holds also if $k_1$ or $k_2$ equal $M_0 < \infty$: Indeed, without loss of generality, suppose $M_0 = k_1 > k_2$. Then, as the sites are unordered, $0 \leq m_1 < m_2$. It follows that $g(m_2)\neq 0 \neq h(k_2)$. So, as $h(k_1)=0$, we have $g(m_1)h(k_1) = 0 \neq g(m_2)h(k_2)$.  Then, by our earlier comment, the `proven first part' applies, and \eqref{I_eq} holds.

Similarly, \eqref{I_eq} holds if $m_1$ or $m_2$ equal $0$:  Indeed, without loss of generality, suppose $0 = m_1 < m_2$. Then, $ k_1 > k_2$. It follows that $g(m_2)\neq 0=g(m_1)$. If $h(k_2)=0$, then $k_2 \geq M_0$ and we have already shown \eqref{I_eq}. If $h(k_2) \neq 0$, then $g(m_1)h(k_1) = 0 \neq g(m_2)h(k_2)$, and \eqref{I_eq} holds by the `proven first part'.
 
 \vskip .2cm
 
 {\it Step 6.}
We now establish \eqref{I_eq} by induction for all other cases. Without loss of generality, suppose $m_1 < m_2$ and $k_1 > k_2$. Assume for our induction step that \eqref{I_eq} holds for a fixed $m_1 \geq 0$ and for all $m_2 > m_1$ and for all $k_1, k_2$ such that $k_1 > k_2$. Our base case, when $m_1=0$, has already been shown.
Suppose that we can show
\begin{eqnarray}
&&\limsup_{N \to \infty} \frac{1}{\|d\|^{n+\alpha}}g(m_1+1)h(k_1-1)\tilde{E}^N\int_0^t \frac{1}{N^n}\sum_{|x|\leq RN} I^{x,d}_{m_1+1, m_2+1, k_1-1, k_2-1}(\eta_s, \xi_s) ds \nonumber\\
&&\leq \limsup_{A \to \infty} \limsup_{N \to \infty} \nonumber\\
&&\ \ \ \ \ \ \ \ \ 8A\kappa \|h\|_\infty \sum\limits_{\|d'\|=1}^\infty \frac{1}{\|d'\|^{n+\alpha}} \tilde{E}^N \int_0^t \frac{1}{N^n}\sum_{|x|\leq RN} I^{x,d}_{m_1, m_2, k_1, k_2}(\eta_s, \xi_s) ds.
\label{step7_eq}
\end{eqnarray}
Then, 
by the induction assumption,
we would have
$$\limsup_{N \to \infty} \frac{1}{\|d\|^{n+\alpha}}g(m_1+1)h(k_1-1)\tilde{E}^N\int_0^t \frac{1}{N^n}\sum_{|x|\leq RN} I^{x,d}_{m_1+1, m_2+1, k_1-1, k_2-1}(\eta_s, \xi_s) ds = 0.$$
As $g(m_1+1)\geq g(1)>0$ and $h(k_1-1) >0$ if $k_1-1< M_0$, the expectation in the above display vanishes as $N\uparrow\infty$ 
for all $m_2+1 > m_1+1$ and for all $k_1-1, k_2-1$ such that $k_1-1 > k_2-1$ provided that $k_1-1< M_0$. However, we have already shown that the expectation vanishes if $k_1-1 \geq M_0$. 

Thus, the induction step and therefore the second part of Lemma \ref{orderinglemma} would be proved.

\vskip .2cm

{\it Step 7.}
To show \eqref{step7_eq}, we recall $\min_{y,y+d'}=\min\{gh_{d'}(\eta(y)), gh_{d'}(\xi(y))\}$ and write
$\tilde{L}I^{x,d}_{m_1, m_2, k_1, k_2}(\eta, \xi)$ equal to
\begin{eqnarray*}
&&\sum\limits_{y \in \Z^n} \sum\limits_{\|d'\|=1}^\infty \frac{1}{\|d'\|^{n+\alpha}}\min\nolimits_{y, y+d'}\big(I^{x,d}_{m_1, m_2, k_1, k_2}(\eta^{y, y+d'}, \xi^{y, y+d'}) -I^{x,d}_{m_1, m_2, k_1, k_2}(\eta, \xi)\big) \\
&&+\sum\limits_{y \in \Z^n} \sum\limits_{\|d'\|=1}^\infty \frac{1}{\|d'\|^{n+\alpha}}(gh_{d'}(\eta(y))-\min\nolimits_{y, y+d'})\\
&&\ \ \ \ \ \ \ \ \ \ \ \ \ \ \ \ \ \ \ \ \ \times \big(I^{x,d}_{m_1, m_2, k_1, k_2}(\eta^{y, y+d'}, \xi) -I^{x,d}_{m_1, m_2, k_1, k_2}(\eta, \xi)\big)\\
&&+ \sum\limits_{y \in \Z^n} \sum\limits_{\|d'\|=1}^\infty \frac{1}{\|d'\|^{n+\alpha}}( gh_{d'}(\xi(y))-\min\nolimits_{y, y+d'})\\
&&\ \ \ \ \ \ \ \ \ \ \ \ \ \ \ \ \ \ \ \ \ \times \big(I^{x,d}_{m_1, m_2, k_1, k_2}(\eta, \xi^{y, y+d'}) -I^{x,d}_{m_1, m_2, k_1, k_2}(\eta, \xi)\big).\end{eqnarray*}

 All of the terms above vanish except those when $y=x$, $y=x+d$, $y+d'=x$, and $y+d'=x+d$. In making a bound, of the positive terms, we shall keep the term $I^{x,d}_{m_1, m_2, k_1, k_2}(\eta^{y, y+d'}, \xi^{y, y+d'}) =I^{x,d}_{m_1+1, m_2+1, k_1-1, k_2-1}(\eta, \xi)$ when $y=x$ and $d'=d$. For the negative terms, we shall double count the terms where $y=x$ and $d'=d$. 
Note also that the total aggregate rate for all the negative terms is $|gh_{d'}(\eta(y))-gh_{d'}(\xi(y)|$. Therefore,
\begin{eqnarray*}
&&\tilde{L}I^{x,d}_{m_1, m_2, k_1, k_2}(\eta, \xi) \geq \frac{1}{\|d\|^{n+\alpha}}\min\nolimits_{x, x+d} I^{x,d}_{m_1+1, m_2+1, k_1-1, k_2-1}(\eta, \xi)\\
&&- \sum\limits_{\|d'\|=1}^\infty \frac{1}{\|d'\|^{n+\alpha}} \sum_{\substack{y=x, x+d,\\ x-d', x+d-d'}}|gh_{d'}(\eta(y))-gh_{d'}(\xi(y)|I^{x,d}_{m_1, m_2, k_1, k_2}(\eta, \xi).
\end{eqnarray*}
Note that in the above display, $\min\nolimits_{x,x+d} = g(m_1+1)h(k_1-1)$.
In the `negative' terms, we bound $|gh_{d'}(\eta(y))-gh_{d'}(\xi(y)|\leq \kappa \|h\|_\infty (\eta(y)+ \xi(y))$ (cf. \eqref{LB}), and split into terms involving only $\eta$ and only $\xi$.  We may introduce indicator functions ${\bf 1}(\eta(y) \leq A)$ and ${\bf 1}(\eta(y) > A)$ onto the `$\eta$' terms and `$\xi$' terms. Items $\eta(y){\bf 1}(\eta(y) \leq A) \leq A$ and also $I^{x,d}_{m_1,m_2,k_1,k_2}(\eta,\xi)\eta(y){\bf 1}(\eta(y) > A) \leq \eta(y){\bf 1}(\eta(y)>A)$, with similar bounds for the items with $\xi$.

We thus obtain, moving the negative terms to the other side of the inequality,
\begin{eqnarray*}
&&\frac{1}{\|d\|^{n+\alpha}}g(m_1+1)h(k_1-1)\tilde{E}^N\int_0^t \frac{1}{N^n}\sum_{|x|\leq RN} I^{x,d}_{m_1+1, m_2+1, k_1-1, k_2-1}(\eta_s, \xi_s) ds \\
&&\ \ \ \ \ \ \ \ \  \leq S_1+S_2+S_3 \ \ \ {\rm where}
\end{eqnarray*}
\begin{eqnarray*}
 S_1&=& \tilde{E}^N\int_0^t \frac{1}{N^n}\sum_{|x|\leq RN}\tilde{L}I^{x,d}_{m_1, m_2, k_1, k_2}(\eta_s, \xi_s) ds,
\end{eqnarray*}
\begin{eqnarray*}
  S_2&=&8A\kappa \|h\|_\infty \sum\limits_{\|d'\|=1}^\infty \frac{1}{\|d'\|^{n+\alpha}} \tilde{E}^N \int_0^t \frac{1}{N^n}\sum_{|x|\leq RN} I^{x,d}_{m_1, m_2, k_1, k_2}(\eta_s, \xi_s) ds \ \ \ \ \ \ {\rm and \ }\\
 S_3 &=& \kappa \|h\|_\infty \sum\limits_{\|d'\|=1}^\infty \frac{1}{\|d'\|^{n+\alpha}}\\
 &&\ \ \ \ \ \times \int_0^t \tilde{E}^N\frac{1}{N^n} \sum_{|x|\leq RN} \sum_{\substack{ y=x, x+d,\\ x-d', x+d-d'}} \big(\eta_s(y){\bf 1}(\eta_s(y)> A)+\xi_s(y){\bf 1}(\xi_s(y) > A)\big)  ds.
  \end{eqnarray*}

We now show the first and third terms vanish as $N, A\uparrow$ to finish, the term $S_2$ being what we would like to keep.   By Lemma \ref{particlebound2} and that $\xi_\cdot\sim\nu_c$, noting $\sum_{\|d'\|\geq 1} \|d'\|^{-(n+\alpha)}<\infty$, the term $S_3$ goes to zero. For the term $S_1$, consider the mean-zero martingale
$\tilde{M}^{x,d}_{m_1, m_2, k_1, k_2}(t)= I^{x,d}_{m_1, m_2, k_1, k_2}(\eta_t, \xi_t) - I^{x,d}_{m_1, m_2, k_1, k_2}(\eta_0, \xi_0) - \int_0^t N\tilde{L}I^{x,d}_{m_1, m_2, k_1, k_2}(\eta_s, \xi_s) ds$.
As $0\leq I^{x,d}_\cdot \leq 1$, it follows that
$$\int_0^t \frac{1}{N^n}\sum_{|x|\leq RN}\tilde{L}I^{x,d}_{m_1, m_2, k_1, k_2}(\eta_s, \xi_s) ds  \leq \frac{1}{N^n}\sum_{|x|\leq RN}\frac{1}{N}\big(1- \tilde{M}^{x,d}_{m_1, m_2, k_1, k_2}(t)\big).$$
Hence, the expected value $S_1\leq C(R)
/N\rightarrow 0$, as $N\uparrow\infty$, completing the proof of \eqref{step7_eq}. \qed

\medskip
\medskip
{\bf Acknowledgement.}  This work was partially supported by ARO grant W911NF-14-1-
0179, and a Daniel Bartlett graduate fellowship.  Thanks to Jianfei Xue for helpful comments on a preliminary version of the article.

\end{document}